%% file: main.tex
\documentclass[10pt, a4paper, twoside, titlepage]{article}
\usepackage[a4paper, inner=1.15in, outer=1.35in, bottom=1.25in, marginparwidth=1in, marginparsep=0.125in]{geometry}
\pdfoutput=1

\usepackage[ngerman, english]{babel} 
\usepackage{amsmath}
\usepackage{amssymb}
\usepackage{amsthm}
\usepackage{thmtools}
\usepackage{algorithm, algpseudocode}
\usepackage{bbm}
\usepackage{blkarray}
\usepackage{pgfplots}
\usetikzlibrary{graphs,arrows.meta,matrix,fit,patterns,shapes}
\pgfplotsset{compat=1.8}
\usepackage{makecell}
\usepackage{svg}
\usepackage{colonequals}
\usepackage{soul}
\usepackage{color}
\usepackage{setspace}
\usepackage{mathtools}
\usepackage{hyperref}
\usepackage[nameinlink,capitalize]{cleveref}
\usepackage[toc,acronym,symbols]{glossaries}
\usepackage{stmaryrd}
\usepackage{eso-pic}
\usepackage{transparent}
\usepackage{graphicx}
\usepackage[framemethod=TikZ]{mdframed}
\usepackage[shortlabels]{enumitem}
\usepackage[nottoc,numbib]{tocbibind}
\usepackage{ifoddpage}
\usepackage{lmodern} 
\usepackage[T1]{fontenc}
\usepackage{FiraSans}
\usepackage{merriweather} 
\usepackage{FiraMono}
\usepackage{mathpazo} 

\usepackage{fancyhdr}
\usepackage{titlesec} 
\numberwithin{equation}{section} 

\usepackage{tabularx}
\usepackage{multicol}
\newcolumntype{Y}{>{\centering\arraybackslash}X}
\usepackage[column=C]{cellspace}
\addparagraphcolumntypes{X}
\definecolor{codegreen}{RGB}{140,153,0}
\definecolor{codegray}{RGB}{128,128,128}
\definecolor{codepurple}{RGB}{173,20,87}
\definecolor{codeorange}{RGB}{240,107,0}
\definecolor{backcolour}{RGB}{237,237,237}
\definecolor{backcolour2}{RGB}{26,26,26}
\definecolor{solarizedbg}{RGB}{253,246,227}
\definecolor{solarizedtext}{RGB}{101,123,131}
\definecolor{codeblue}{RGB}{3,135,209}

\definecolor{tumBlue}{RGB}{0,101,189} 
\definecolor{tumDarkBlue}{RGB}{0,82,147} 
\definecolor{tumLightBlue}{RGB}{100,160,200} 
\definecolor{tumLighterBlue}{RGB}{152,198,234} 
\definecolor{tumOrange}{RGB}{227,114,34} 
\definecolor{tumGreen}{RGB}{162,173,0} 
\definecolor{tumGray}{RGB}{153,153,153} 
\definecolor{tumLight}{RGB}{218,215,203} 

\colorlet{sectionblue}{tumBlue}
\definecolor{linkred}{RGB}{127,0,0} 
\definecolor{darklinkred}{RGB}{50,0,0} 
\colorlet{headingcolor}{sectionblue}
\colorlet{headingcolormuted}[RGB]{headingcolor!20!white}
\colorlet{linkcolor}{linkred}

\makeatletter
\DeclareRobustCommand{\xhl}[1]{%
 \protected@edef\tmp{#1}%
 \expandafter\hl\expandafter{\tmp}%
}
\makeatother
\newcommand{\hlyellow}[1]{{\sethlcolor{solarizedbg}\protect\xhl{#1}}}
\newcommand*{\yellowemph}[1]{%
  \ifmmode\mathchoice{\tikz[baseline=(text.base)]\node(text)[rectangle, fill=solarizedbg, inner sep=0pt]{$\displaystyle #1$};}{\tikz[baseline=(text.base)]\node(text)[rectangle, fill=solarizedbg, inner sep=0pt]{$\textstyle #1$};}{\tikz[baseline=(text.base)]\node(text)[rectangle, fill=solarizedbg, inner sep=0pt]{$\scriptstyle #1$};}{\tikz[baseline=(text.base)]\node(text)[rectangle, fill=solarizedbg, inner sep=0pt]{$\scriptscriptstyle #1$};}\else\hlyellow{#1}\fi
}

\usepackage[labelfont={color=headingcolor,bf}]{caption}

\definecolor{darkblue}{RGB}{0,102,204}
\colorlet{colorlow}[hsb]{darkblue}
\definecolor{darkred}{RGB}{204,26,0}
\colorlet{colorhigh}[hsb]{darkred}
\colorlet{color1}[RGB]{colorlow!100!colorhigh}
\colorlet{color2}[RGB]{colorlow!85!colorhigh}
\colorlet{color3}[RGB]{colorlow!70!colorhigh}
\colorlet{color4}[RGB]{colorlow!50!colorhigh}
\colorlet{color5}[RGB]{colorlow!30!colorhigh}
\colorlet{color6}[RGB]{colorlow!15!colorhigh}
\colorlet{color7}[RGB]{colorlow!0!colorhigh}

\makeatletter
\def\moverlay{\mathpalette\mov@rlay}
\def\mov@rlay#1#2{\leavevmode\vtop{%
   \baselineskip\z@skip \lineskiplimit-\maxdimen
   \ialign{\hfil$\m@th#1##$\hfil\cr#2\crcr}}}
\newcommand{\charfusion}[3][\mathord]{
    #1{\ifx#1\mathop\vphantom{#2}\fi
        \mathpalette\mov@rlay{#2\cr#3}
      }
    \ifx#1\mathop\expandafter\displaylimits\fi}
\makeatother

\newcommand{\midd}[1]{\mathrel{}\middle#1\mathrel{}}

\mdfsetup{ skipbelow=-0.2em }

\newlength{\spaceblength}
\declaretheoremstyle[
    headfont=\bfseries,
    notefont=\bfseries,
    notebraces={}{\\[\parskip]}, 
    bodyfont=\normalfont\upshape,
    headpunct={},
    postheadspace=\spaceblength,
    spacebelow=\parskip,
    spaceabove=\parskip,
    headformat={%
        \checkoddpage\ifoddpage\rlap{\hskip\textwidth\hskip10pt\color{headingcolor}\ \ \ \NAME\ \NUMBER}\hskip-\spaceblength{\NOTE}%
        \else\makebox[0pt][r]{\color{headingcolor}\NAME\ \NUMBER\hskip10pt\ \ \ \ }\hskip-\spaceblength{\NOTE}\fi%
    },
    mdframed={
        nobreak=true,
        linecolor=headingcolor!20,
        innertopmargin=8pt,
        innerbottommargin=6pt,
        innerleftmargin=8pt,
        innerrightmargin=8pt,
        skipabove=0.7em,
        linewidth=2pt } 
]{boxstyle}
\declaretheoremstyle[
    headfont=\bfseries\itshape,
    notefont=\normalfont\bfseries,
    notebraces={}{\\[\parskip]}, 
    bodyfont=\normalfont,
    headpunct={\indent},
    postheadspace=\spaceblength,
    spacebelow=\parskip,
    spaceabove=\parskip,
    headformat={%
        \checkoddpage\ifoddpage\rlap{\hskip\textwidth\color{headingcolor}\ \ \ \NAME}\hskip-\spaceblength{\NOTE}%
        \else\makebox[0pt][r]{\color{headingcolor}\NAME\ \ \ \ }\hskip-\spaceblength{\NOTE}\fi%
    },
    qed=\qedsymbol
]{proofstyle}
\declaretheoremstyle[
    headfont=\bfseries,
    notefont=\bfseries,
    notebraces={}{\\[\parskip]}, 
    bodyfont=\normalfont,
    headpunct={\indent},
    postheadspace=\spaceblength,
    spacebelow=\parskip,
    spaceabove=\parskip,
    headformat={%
        \checkoddpage\ifoddpage\rlap{\hskip\textwidth\color{headingcolor}\ \ \ \NAME\ \NUMBER}\hskip-\spaceblength{\NOTE}%
        \else\makebox[0pt][r]{\color{headingcolor}\NAME\ \NUMBER\ \ \ \ }\hskip-\spaceblength{\NOTE}\fi%
    },
    qed={$\begingroup\color{headingcolormuted}\blacktriangleleft\endgroup$}
]{examplestyle}
\declaretheoremstyle[
    headfont=\normalfont,
    notefont=\bfseries,
    notebraces={}{\\[\parskip]}, 
    bodyfont=\normalfont,
    headpunct={\indent},
    postheadspace=\spaceblength,
    spacebelow=\parskip,
    spaceabove=\parskip,
    headformat={%
        \checkoddpage\ifoddpage\rlap{\hskip\textwidth\color{headingcolor}\ \ \ \NAME}\hskip-\spaceblength{\NOTE}%
        \else\makebox[0pt][r]{\color{headingcolor}\NAME\ \ \ \ }\hskip-\spaceblength{\NOTE}\fi%
    },
    qed={$\begingroup\color{headingcolormuted}\blacktriangleleft\endgroup$}
]{remarkstyle}

\declaretheorem[style=boxstyle,numberwithin=section]{definition}
\declaretheorem[style=boxstyle,sibling=definition]{theorem}

\declaretheorem[style=boxstyle,sibling=definition]{lemma}
\declaretheorem[style=boxstyle,sibling=definition]{corollary}

\declaretheorem[style=boxstyle,sibling=definition]{conjecture}

\declaretheorem[style=proofstyle,numbered=no,name=Proof]{tproof}
\declaretheorem[style=proofstyle,numbered=no,name=Proof Sketch]{sproof}

\addto\extrasenglish{

}

\newlength{\secskip}
\titleformat{\section}
  {\LARGE\bfseries\color{headingcolor}\sffamily}
  {\makebox[0pt][r]{\thesection\hspace*{\secskip}}}{0em}{}
\titleformat{\subsection}
  {\large\bfseries\color{headingcolor}\sffamily}
  {\makebox[0pt][r]{\thesubsection\hspace*{\secskip}}}{0em}{}
\titleformat{\subsubsection}
  {\bfseries\color{headingcolor}\sffamily}
  {\makebox[0pt][r]{\thesubsubsection\hspace*{\secskip}}}{0em}{}

\hypersetup{
    colorlinks=true,
    allcolors=.,
    urlcolor=linkcolor,
    citecolor=linkcolor,
    linkcolor=linkcolor,
}

\newcommand{\spectitle}[3]{
	\sffamily
	\def\svgwidth{59.99796366911307pt}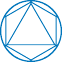\hfill\raisebox{5pt}[0pt][0pt]{\def\svgwidth{96.05263069425209pt}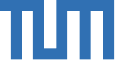}\\[8em]
	{\large\merriweathersanslight Technical University of Munich\\[2em]
	Department of Mathematics\\[2em]
	Master's Thesis in Mathematics}\\[8em]
	{\color{headingcolor}\Huge\bfseries #1}\\[4em]
	{\LARGE #2}\\[10em]
	{\large\merriweathersanslight {\sffamily Supervisor:} Prof.~Dr.~rer.~nat.~habil.~Nina Gantert\\[1.5em]
	{\sffamily Submission Date:} #3}
}

\newtagform{roman}[]()

\DeclareMathOperator{\GEO}{Geo}

\DeclareMathOperator{\EXPD}{Exp}

\DeclareMathOperator{\BETAD}{B}

\DeclareMathOperator{\BR}{br}

\DeclareMathOperator{\DIST}{dist}
\DeclareMathOperator{\POLYLOG}{Li}

\newcommand{\Prb}[1] {\mathbb{P}\left[#1\right]}
\newcommand{\Prbc}[2] {\mathbb{P}_{#1}\left[#2\right]}

\newcommand{\Prbcc}[3] {\mathbb{P}_{#1}^{#2}\left[#3\right]}
\newcommand{\bPrb}[1] {\mathbf{P}\left[#1\right]}

\newcommand{\Ex}[1] {\mathbb{E}\left[#1\right]}
\newcommand{\Exc}[2] {\mathbb{E}_{#1}\left[#2\right]}
\newcommand{\Excc}[3] {\mathbb{E}_{#1}^{#2}\left[#3\right]}

\newcommand{\Geo}[1] {\GEO_{\geq 0}\left(#1\right)}
\newcommand{\Geob}[1] {\GEO_{\geq 1}\left(#1\right)}
\newcommand{\Betad}[2] {\BETAD\left(#1, #2\right)}

\newcommand{\Expd}[1] {\EXPD\left(#1\right)}
\newcommand{\Polylog}[2] {\POLYLOG_{#1}\left(#2\right)}
\newcommand{\br}[1] {\BR\left(#1\right)}
\newcommand{\dx}[1] {\;\mathrm{d}#1}
\newcommand{\mudx}[2] {\;#1\left(\mathrm{d}#2\right)}

\newcommand{\bbP} {\mathbb{P}}
\newcommand{\bbE} {\mathbb{E}}

\newcommand{\bbQ} {\mathbb{Q}}
\newcommand{\bbR} {\mathbb{R}}

\newcommand{\bbN} {\mathbb{N}}
\newcommand{\bbZ} {\mathbb{Z}}

\newcommand{\frP} {\mathfrak{P}}

\newcommand{\calG} {\mathcal{G}}

\newcommand{\calP} {\mathcal{P}}

\newcommand{\bP} {\mathbf{P}}
\newcommand{\bQ} {\mathbf{Q}}

\newcommand{\dist}[2] {\DIST\left(#1, #2\right)}
\newcommand{\indic}[1] {\mathbbm{1}_{#1}}
\DeclareMathOperator*{\argmax}{arg\,max}
\DeclareMathOperator*{\argmin}{arg\,min}

\newcommand{\mmod}{\mathrel{}\operatorname{mod}\mathrel{}}
\newcommand{\stickman}[2] {\begin{scope}[shift={#1}]
  \draw [#2] (-0.2, -0.5) -- (0,-0.25) -- (0,0.25);
  \draw [#2] (0.2, -0.5) -- (0,-0.25);
  \draw [#2] (-0.2, 0.15) -- (0,-0.1);
  \draw [#2] (0.2, 0.15) -- (0,-0.1);
  \fill [#2] (0,0.25) circle (0.125);
\end{scope}}
\newcommand*\circledchar[1]{\tikz[baseline=(char.base)]{
  \node[shape=circle,draw,inner sep=1pt,outer sep=1pt] (char) {#1};}}

\tikzset{
  diagonal fill a/.style n args=2{path picture={%
  \draw[fill=#1, draw=none] (path picture bounding box.south west) --
              (path picture bounding box.north east) -- (path picture bounding box.south east) -- cycle; %
  \draw[fill=#2, draw=none] (path picture bounding box.south west) --
              (path picture bounding box.north east) -- (path picture bounding box.north west) -- cycle;}},
  diagonal fill b/.style n args=2{path picture={%
  \draw[fill=#1, draw=none] (path picture bounding box.south west) --
              (path picture bounding box.north east) -- (path picture bounding box.south east) -- cycle; %
  \draw[pattern=north west lines, pattern color=#2, draw=none] (path picture bounding box.south west) --
              (path picture bounding box.north east) -- (path picture bounding box.north west) -- cycle;}},
  diagonal fill c/.style n args=2{path picture={%
  \draw[pattern=north west lines, pattern color=#1, draw=none] (path picture bounding box.south west) --
              (path picture bounding box.north east) -- (path picture bounding box.south east) -- cycle; %
  \draw[fill=#2, draw=none] (path picture bounding box.south west) --
              (path picture bounding box.north east) -- (path picture bounding box.north west) -- cycle;}},
  diagonal fill d/.style n args=2{path picture={%
  \draw[pattern=north west lines, pattern color=#1, draw=none] (path picture bounding box.south west) --
              (path picture bounding box.north east) -- (path picture bounding box.south east) -- cycle; %
  \draw[pattern=north west lines, pattern color=#2, draw=none] (path picture bounding box.south west) --
              (path picture bounding box.north east) -- (path picture bounding box.north west) -- cycle;}},
  table nodes/.style={
    rectangle,
    draw=none,
    align=center,
    minimum height=7mm,
    text depth=0.5ex,
    text height=2ex,
    inner xsep=0pt,
    outer sep=0pt
  },      
  table/.style={
    matrix of nodes,
    row sep=-\pgflinewidth,
    column sep=-\pgflinewidth,
    nodes={
        table nodes
    }
  }
}

\usepackage{ifthen}
\newboolean{print}
\setboolean{print}{false}

\title{Variations on Reinforced Random Walks}
\author{Fabian Michel}
\date{16/09/2022}

\begin{document}

\begin{titlepage}
	\makeatletter
	
	\sffamily
	\def\svgwidth{59.99796366911307pt}\input{graphics/tum_mathematik_svg-tex.pdf_tex}\hfill\raisebox{5pt}[0pt][0pt]{\def\svgwidth{96.05263069425209pt}\input{graphics/TUM_svg-tex.pdf_tex}}\\[8em]
	{\large\merriweathersanslight Technical University of Munich\\[2em]
	Department of Mathematics\\[2em]
	Master's Thesis in Mathematics}\\[8em]
	{\color{headingcolor}\Huge\bfseries Variations on\\[0.3em]Reinforced Random Walks}\\[4em]
	{\LARGE \@author}\\[10em]
	{\large\merriweathersanslight {\sffamily Supervisor:} Prof.~Dr.~rer.~nat.~habil.~Nina Gantert\\[1.5em]
	{\sffamily Submission Date:} \@date}

	\makeatother
\end{titlepage}

\ifthenelse{\boolean{print}}{\thispagestyle{empty} { ~ } \newpage}{}

\thispagestyle{empty}
{
	\makeatletter
	
	\sffamily
	\def\svgwidth{59.99796366911307pt}\input{graphics/tum_mathematik_svg-tex.pdf_tex}\hfill\raisebox{5pt}[0pt][0pt]{\def\svgwidth{96.05263069425209pt}\input{graphics/TUM_svg-tex.pdf_tex}}\\[6em]
	{\large\merriweathersanslight Technical University of Munich\\[1.5em]
	Department of Mathematics\\[1.5em]
	Master's Thesis in Mathematics}\\[6em]
	{\color{headingcolor}\Huge\bfseries Variations on\\[0.3em]Reinforced Random Walks}\\[2em]
	{\color{headingcolor}\huge Abwandlungen\\[0.3em]selbstverst\"arkender Irrfahrten}\\[3em]
	{\LARGE \@author}\\[8em]
	{\large\merriweathersanslight {\sffamily Supervisor:} Prof.~Dr.~rer.~nat.~habil.~Nina Gantert\\[1.5em]
	{\sffamily Submission Date:} \@date}

	\makeatother
	\ifthenelse{\boolean{print}}{}{~\\[5em]\noindent {\large\merriweathersanslight final digital version}}
}
\newpage

\ifthenelse{\boolean{print}}{\thispagestyle{empty} { ~ } \newpage}{}

\thispagestyle{empty}
{
	\makeatletter
	~\\[20em]
	\noindent I confirm that this master's thesis is my own work and I have documented all sources and material used.\\[1em]
	\noindent Ich erkl\"are hiermit, dass ich diese Arbeit selbst\"andig und nur mit den angegebenen Hilfsmitteln angefertigt habe.
	
	\noindent ~\\[3em]
	\rule{8cm}{.4pt}\\
	\@author, M\"unchen, \@date
	\makeatother
}
\newpage

{
  \thispagestyle{fancy}
  ~\\[4em]
  \begin{center}
    \LARGE\bfseries\color{headingcolor}\sffamily Abstract
  \end{center}

  This thesis examines edge-reinforced random walks with some modifications to the standard definition.
  An overview of known results relating to the standard model is given and
  the proof of recurrence for the standard linearly edge-reinforced random walk on bounded degree graphs
  with small initial edge weights is repeated. Then, the edge-reinforced random walk with multiple
  walkers influencing each other is considered. The following new results are shown: on a segment of
  three nodes, the edge weights resemble a P\'olya urn and the fraction of the edge weights divided by
  the total weight forms a converging martingale. On $\bbZ$, the behavior is the same as for a single
  walker~-- either all walkers have finite range or all walkers are recurrent. Finally, edge-reinforced
  random walks with a bias in a certain direction are analysed, in particular on $\bbZ$. It is shown
  that the bias can introduce a phase transition between recurrence and transience, depending on the
  strength of the bias, thus fundamentally altering the behavior in comparison to the standard linearly
  reinforced random walk.

  ~\\[2em]
  \begin{center}
    \LARGE\bfseries\color{headingcolor}\sffamily Zusammenfassung
  \end{center}

  Diese Masterarbeit betrachtet (kanten-)selbstverst\"arkende Irrfahrten mit einigen Ver\"anderungen
  im Vergleich zur g\"angigen Definition. Es wird eine \"Ubersicht \"uber bekannte Ergebnisse zum g\"angigen Modell
  gegeben und der Beweis, dass die linear selbstverst\"arkende Irrfahrt auf Graphen mit beschr\"anktem
  Grad und hinreichend kleinen anf\"anglichen Kantengewichten rekurrent ist, wird wiederholt. Danach
  werden selbstverst\"arkende Irrfahrten mit mehreren Walkern, die sich gegenseitig beeinflussen,
  untersucht. Die folgenden neuen Resultate werden bewiesen: Auf einer Strecke mit drei Knoten \"ahneln
  die Kantengewichte einer P\'olya-Urne und der Anteil der Kantengewichte am Gesamtgewicht bildet ein
  konvergierendes Martingal. Auf $\bbZ$ ist das Verhalten das gleiche wie f\"ur einen einzelnen Walker~--
  entweder besuchen alle Walker nur einen endlichen Teil des Graphen oder alle sind rekurrent.
  Zum Schluss wird die selbstverst\"arkende Irrfahrt mit Bias in eine bestimmte Richtung betrachtet,
  vor allem auf $\bbZ$. Es wird gezeigt, dass der Bias einen Phasen\"ubergang zwischen Rekurrenz und
  Transienz verursachen kann, der von der St\"arke des Bias abh\"angt, und damit das Verhalten im
  Vergleich zur normalen linear selbstverst\"arkenden Irrfahrt grundlegend \"andert.
}
\newpage

{
\hypersetup{
  hidelinks
}
\tableofcontents
}

\newpage

{
\hypersetup{
  hidelinks
}
\listoffigures
\listoftables

~\\[2em]
\subsection*{A note on cross-references in this thesis}
Cross-references and citations are marked in \textcolor{linkred}{red}, and they are clickable and
directly link to the referenced object in the digital version. Titles, definitions, theorems,
etc.~are colored in \textcolor{sectionblue}{blue}.
}

\newpage

\section{Introduction}
\label{sec:intro}

The central topic of this thesis is the edge-reinforced random walk (ERRW), a special type of random
walk on a graph. The edges in the graph are weighted, and the probability to leave a node via one of
the incident edges is proportional to the respective edge weight compared to the weights of the other
incident edges. Each time an edge is crossed, its weight is increased according to some reinforcement
scheme. Thus, it becomes more likely to visit parts of the graph which have already been visited before
again. Next to reviewing results on this model, two modifications are considered: introducing multiple
random walkers which influence each other and introducing a bias in a certain direction.

The ERRW is harder to analyze compared to random processes which possess the Markov property.
In contrast to Markov chains (MCs), where the transition probabilities from one state to another only depend
on the current state, but not on the past, the transition probabilities of the ERRW change over
time as the edge weights are reinforced. The Markov property is lost. However, for a certain class of ERRWs,
namely the linearly edge-reinforced random walk (LERRW), where the edge weights are increased by a constant increment upon every traversal,
the reinforced random walk is equal in law to a mixture of MCs, i.e.~a MC with
random transition probabilities. This equivalence was the basis for many of the known results for
ERRWs, since it allows us to use tools developed for MCs on the ERRW as well.

The new variations introduced in this thesis, multiple walkers and a bias, break this connection to
MCs as well (at least in most of the considered cases, if not in all of them). The walk can
no longer be represented as a mixture of MCs, and has to be analyzed using different techniques.
The results presented here indicate that multiple walkers do not fundamentally change the behavior of the
LERRW, while a bias is strong enough to do so.

\subsection{Literature}

The ERRW, as well as its counterpart, the vertex-reinforced random walk (VRRW), have been studied extensively,
the first papers dating back to 1987, when the model was introduced by Coppersmith and Diaconis. Even before,
\cite{definettimarkovchain} showed that the LERRW (if certain assumption are satisfied) has a representation
as a mixture of MCs. Much later, \cite{rwforerrw} showed that this representation can be used for the
LERRW on any graph, and \cite{magicformula} even gave a formula for the so-called mixing measure on finite
graphs. The mixing measure simply is the distribution of the random transition probabilities in the mixture of
MCs.

Relatively early, results for the LERRW on trees were obtained. \cite{errwpemantle,rwrelyonspemantle,rwrelyonspemantlecorr}
showed that there is a phase transition between recurrence and transience in the initial edge weights. The interesting
case of $\bbZ^d$ remained open much longer until \cite{localizationerrw,errwandvrjp,transienceerrw} showed that for $d \geq 3$,
there is again a transition from recurrence to transience in the initial weights. \cite{monphasetrans} proved that this transition
is sharp, i.e.~there is a certain critical initial edge weight such that for smaller initial weights, the random walk is
recurrent, and transient for larger weights.

In parallel, the VRRW was analyzed, but with completely different tools (there is no representation as
a mixture of MCs). \cite{vrrwfinitegraphs} analyzed the behavior of the VRRW on finite graphs
with the help of a so-called stochastic approximation: the evolution of the vertex weights is approximated by a differential
equation. \cite{vrrwzfiniterange,vrrwzstuckonfive} proved that the VRRW (with linear reinforcement) gets
stuck on $5$ nodes on $\bbZ$, a result which could in part be generalized to arbitrary graphs: \cite{vrrwarbitrary}
showed that the VRRW gets stuck with positive probability on certain finite subgraphs for almost any graph.
The behavior of the VRRW is thus largely different from the behavior of the ERRW, and it was therefore
unclear how the variations of the ERRW considered here would affect the random walk.

An overview of results on reinforced processes in general can be found in \cite{surveyreinf,rrwsurveykozma}. These surveys
also show that reinforced random walks are closely related to urn processes, which have a very similar reinforcement
component to the linearly reinforced walks: in most urn models, when a ball of a certain color is drawn, a fixed number of
balls of the same color is added to the urn. This is also a type of linear reinforcement, and urns have been used on many
occasions to analyze reinforced walks (see, for example, \cite{errwpemantle,lerrwbsc}).

Reinforced walks with multiple walkers or with a bias have (to the best of our knowledge) not been studied yet.
There is however, another relevant part of the literature which deals with Markovian random walks. Markov chains
are still a vital tool to understand the more complicated reinforced walks. \cite{rwandust,probtreenet} are two
very good resources for results and techniques which can be applied to MCs. Finally, the references
list additional sources which look at various aspects of random walks, including more and different models for
reinforced walks.

\subsection{Main Results}

The following results are presented in this thesis. After covering preliminary definitions in \autoref{sec:prelim},
some of the most important results on reinforced random walks, in particular in relation to this thesis, are presented
in \autoref{sec:selected_results}, without proof. \autoref{sec:recurrence_errw} repeats the proof of one such result:
the recurrence of the LERRW on bounded-degree graphs for small initial weights (see \autoref{thm:rec_bndd_deg}).

\autoref{sec:multiple_walkers} covers the first new variation to the LERRW: multiple walkers influencing each
other. \autoref{cor:two_pl_urn_conv_rand} shows that edge weights of the LERRW for $2$ walkers on a segment of
$3$ nodes behave similarly to the proportion of balls in a P\'olya urn, and it is conjectured in \autoref{conj:two_pl_urn_limit_rand}
that the limit of the edge weight proportions will have similar properties to the limit of the proportion of balls in the urn.
\autoref{thm:allRecurrentOrAllFiniteRange} considers the ERRW with more general reinforcement and multiple walkers
on $\bbZ$ and shows that its behavior is similar to the case with a single walker: either all walkers are recurrent, or all
walkers have finite range.

In \autoref{sec:biased}, the LERRW on $\bbZ$, but with an additional bias is considered. \autoref{lem:asconv} proves
that, for a multiplicative bias, the probability to move in the direction of bias converges to $1$ for nodes which are visited
infinitely often. This might hint at transience, as noted in \autoref{conj:trans}. For an additive bias, the representation as
a mixture of MCs is recovered, allowing us to show a phase transition between recurrence and transience in the bias
in \autoref{thm:add_bias_rectrans}, and it is even possible to show positive speed if the bias is strong enough in
\autoref{thm:add_bias_speed}. Finally, \autoref{ssec:trans_env} considers the LERRW with a strong bias in the
initial edge weights: the initial edge weight to the right of a node $z \in \bbZ$ is set to $\lambda^z$ for some parameter $\lambda$.
The unfinished calculation in \autoref{ssec:trans_env} seems to indicate that the reinforced walk will be transient, i.e.~that
the initially biased environment dominates the reinforcement. This is similar in spirit to \cite[Theorem 4]{monphasetrans},
where it was shown that if the initial weights are set to the conductances of a recurrent MC, then the
LERRW is also recurrent. In \autoref{ssec:trans_env}, we have an initially transient environment. However, there
cannot be a theorem like \cite[Theorem 4]{monphasetrans} for transience, as can already be seen by looking at
\autoref{thm:rec_bndd_deg}.

Finally, \autoref{sec:conclusion} reviews the results on the new variations of the ERRW obtained in this thesis, and
looks at the many open questions which still remain.

\clearpage
\newpage

\section{Preliminaries}
\label{sec:prelim}

Throughout this thesis, random walks on finite and infinite connected graphs will be considered. All
considered graphs will be locally finite, that is, the number of edges incident to a node is always
finite. The following notation is used in connection with graphs:

\begin{definition}[Graph]
  \label{def:graph}
  A (possibly infinite) \textbf{graph} $G$ is a tuple $\left(V, E\right)$ of vertices (or nodes) and
  edges, with  $E \subseteq \left\{ \left\{u, v\right\} \midd| u,v \in V \land u \neq v \right\}$.
  Sometimes, it is useful to consider directed edges, hence we also define the set
  $\ensuremath{\protect\overrightarrow{E}} := \left\{ \left(u, v\right) \midd| u,v \in V \land \left\{u, v\right\} \in E \right\}$
  of directed edges. For $e = \left(u, v\right) \in \overrightarrow{E}$, we call $\ensuremath{\protect\overleftarrow{e}} := \left(v, u\right)$
  the corresponding reversed edge, and we write $\check{e} = u, \hat{e} = v$, so $e = \left(\check{e}, \hat{e}\right)$.
  We write $\ensuremath{u \sim v}$ if $\left\{u, v\right\} \in E$. $\ensuremath{\dist{u}{v}}$ is the edge length of the shortest
  path between the two nodes and $\ensuremath{\dist{v}{e}}$ is the minimum of the distances of $v$ to the two
  endpoints of $e$.
\end{definition}

\subsection{Random Walks}
\label{ssec:rw}

One type of random walk on a graph is the RWRE, which is used as a tool to
analyze reinforced random walks. The RWRE is basically a Markov chain, but with random
transition probabilities, which are also referred to as the environment. The formal definition is as follows:

\begin{definition}[Random Walk in Random Environment]
  \label{def:rwre}
  Consider a graph $G = \left(V, E\right)$. Let $\ensuremath{\mathbf{c}} := \left(c_e\right)_{e \in E}$
  be a collection of random variables over some probability space with $c_e > 0$.
  For $v \in V$, we further denote by $\ensuremath{c_v} := \sum_{u \in V: \ensuremath{u \sim v}} c_{\left\{u,v\right\}}$
  the sum of the conductances of all adjacent edges, and assume $c_v > 0$.
  The \textbf{environment} $\ensuremath{P_{\mathbf{c}}} : V \times V \to \left[0, 1\right]$ on the graph $G$ is defined by
  \begin{align*}
    P_{\mathbf{c}}\left(u, v\right) := \frac{c_{\left\{u,v\right\}}}{c_u}
  \end{align*}
  The Markov chain (MC) $\left(X_n\right)_{n \geq 0}$, conditioned on the knowledge of $\mathbf{c}$, and starting at $v \in V$, is now
  defined by
  \begin{align*}
    \bbP_v^\mathbf{c}\left[X_0 = v\right] = 1
    \qquad
    \Prbcc{v}{\mathbf{c}}{X_{n+1} = t \midd| X_n = u} = P_{\mathbf{c}}\left(u, t\right)
  \end{align*}
  $X_n$ is a random walk taking values in $V$.

  The distribution of $\mathbf{c}$ is called $\ensuremath{\bP}$, and is called the mixing measure,
  and the probability measure $\bbP_v$ is defined by setting
  \begin{align*}
    \Prbc{v}{\cdot} := \int \Prbcc{v}{\mathbf{c}}{\cdot} \mudx{\bP}{\mathbf{c}}
  \end{align*}
  $\left(X_n\right)_{n \geq 0}$, distributed according to $\bbP_v$, is called \textbf{random walk in random environment}
  (RWRE), or mixture of MCs.
\end{definition}

For more background on conductances as well as electrical networks and their relation to
random walks, \cite{probtreenet} gives an excellent introduction (in Chapters 1 and 2).

We shall also use the following terminology:
\begin{enumerate}[(1)]
  \item \label{term:statmeas} $\ensuremath{\pi_{\mathbf{c}}}: V \to \left[0, \infty\right)$ is called a \textbf{stationary
  measure} for the MC given by the conductances $\mathbf{c}$ if it satisfies
  $\pi_{\mathbf{c}}\left(v\right) = \sum_{u \in V: \ensuremath{u \sim v}} \pi_{\mathbf{c}}\left(u\right) P_{\mathbf{c}}\left(u, v\right)$
  for all $v \in V$.
  $\pi_{\mathbf{c}}$ is further called a stationary distribution if it is a probability measure.
  Setting $\pi_{\mathbf{c}}\left(v\right) := c_v$ defines a stationary (and reversible) measure.
  \item \label{term:rectrans} The MC given by $\mathbf{c}$ is called
  \textbf{recurrent} if $\Prbcc{v}{\mathbf{c}}{X_n = v \textrm{ for infinitely many } n} = 1$, and
  \textbf{transient} otherwise. In the transient case, it holds that
  $\Prbcc{v}{\mathbf{c}}{X_n = v \textrm{ for infinitely many } n} = 0$, so the starting node $v$ is only visited
  finitely often a.s., and $\Prbcc{v}{\mathbf{c}}{\exists n > 0: X_n = v} < 1$. In the recurrent
  case, it holds that $\Prbcc{v}{\mathbf{c}}{\exists n > 0: X_n = v} = 1$ (obviously). Since the considered
  MCs are a.s.~irreducible (every node can be reached from every other node
  with positive probability), this corresponds to the standard notion of recurrence and transience.

  The MC is further called \textbf{positive recurrent} if the expected return time to the
  root is finite, i.e.~$\Excc{v}{\mathbf{c}}{\min\left\{ n > 0 : X_n = v \right\}} < \infty$ (this implies recurrence).
  Positive recurrence is equivalent to the existence of a (unique) stationary distribution. Note that the
  existence of a stationary measure does not imply the existence of a stationary distribution since
  the stationary measure may be infinite (if it is finite, it can always be normalized). Indeed,
  in the null recurrent case (recurrent, but not positive recurrent), there exists a unique
  stationary measure which is infinite.
\end{enumerate}

We next define the ERRW, which is one of the main topics of this thesis. The ERRW
is not a MC since the transition probabilities change over time.

\begin{definition}[Edge-Reinforced Random Walk]
  \label{def:errw}
  Consider a graph $G = \left(V, E\right)$ and a counting function $z : \bbN_{\geq 0} \times E \to \bbN_{\geq 0}$
  where we set $z\left(0, e\right) = 0$ for all $e \in E$. Let $W_e : \bbN_{\geq 0} \to \left(0, \infty\right)$ be
  weight functions for $e \in E$. We now define the evolution of $z$ as well as
  the \textbf{edge-reinforced random walk} (ERRW) $\left(X_n\right)_{n \geq 0}$ on $G$
  as follows:
  \begin{align*}
    \textrm{for } n \geq 1: \textrm{ }
    \ensuremath{z\left(n, \left\{u,v\right\}\right)} = \sum_{i=1}^n \mathbbm{1}_{\left(X_{i-1} = u \land X_i = v\right) \lor \left(X_{i-1} = v \land X_i = u\right)}
  \end{align*}
  i.e.~$z$ counts the number of edge traversals, and
  \begin{align*}
    \Prbc{v}{X_0 = v} &= 1 \\
    \Prbc{v}{X_{n+1} = t \midd| X_n = u, \ldots, X_0 = x_0} &= \frac{W_{\left\{u, t\right\}}\left( z\left(n, \left\{u, t\right\}\right) \right)}{\sum_{s \in V: s \sim u} W_{\left\{u, s\right\}}\left( z\left(n, \left\{u, s\right\}\right) \right)}
  \end{align*}
  i.e.~the probability to make a transition from $u$ to $t$ at time $n$ is proportional to the weight $W_{\left\{u, t\right\}}\left( z\left(n, \left\{u, t\right\}\right) \right)$ associated
  to the edge $\left\{u, t\right\}$, which depends on the number of traversals of the edge. We set
  \begin{align*}
    w\left(n, \left\{u, t\right\}\right) := W_{\left\{u, t\right\}}\left( z\left(n, \left\{u, t\right\}\right) \right)
  \end{align*}
  The random variables $w$ thus describe the evolution of the edge weights to which the transition
  probabilities are proportional.
\end{definition}

It is also possible to make the edge weight depend on more than just the number of edge traversals and
the edge, but we restrict ourselves to this case in this thesis which is already quite general. If
$W_e$ is the same for every $e \in E$, then we will write $W$ instead of $W_e$. One important choice
of $W$ is $W\left(n\right) = 1 + n$, which corresponds to starting with initial edge weights $1$ everywhere
and incrementing the weight of the traversed edge by $1$ at every step. We define the
\textbf{linearly edge-reinforced random walk} (LERRW) in a way which is a little more general: the LERRW
is the reinforced random walk where $W_e\left(n\right) = a_e + n$ for some choice of initial weights
$\left(a_e\right)_{e \in E} \in \left(0, \infty\right)^E$.

Later, it will be useful to count the number of directed edge traversals. We therefore also define the
following counting function:
\begin{align*}
  \ensuremath{\protect\overrightarrow{z}\left(n, \left(u,v\right)\right)} = \sum_{i=1}^n \mathbbm{1}_{X_{i-1} = u \land X_i = v}
\end{align*}

Another type of random walk is obtained by using vertex weights instead of edge weights for reinforcing.

\begin{definition}[Vertex-Reinforced Random Walk]
  \label{def:vrrw}
  Consider a graph $G = \left(V, E\right)$ and a counting function $z : \bbN_{\geq 0} \times V \to \bbN_{\geq 0}$
  where we set $z\left(0, v\right) = 0$ for all $v \in V$. Let $W_v : \bbN_{\geq 0} \to \left(0, \infty\right)$ be
  weight functions for $v \in V$. We now define the evolution of $z$ as well as
  the \textbf{vertex-reinforced random walk} (VRRW) $\left(X_n\right)_{n \geq 0}$ on $G$
  as follows:
  \begin{align*}
    \textrm{for } n \geq 1: \textrm{ }
    \ensuremath{z\left(n, v\right)} = \sum_{i=1}^n \mathbbm{1}_{X_i = v}
  \end{align*}
  i.e.~$z$ counts the number of vertex visits, and
  \begin{align*}
    \Prbc{v}{X_0 = v} = 1 \qquad
    \Prbc{v}{X_{n+1} = t \midd| X_n = u, \ldots, X_0 = x_0} = \frac{W_t\left( z\left(n, t\right) \right)}{\sum_{s \in V: s \sim u} W_s\left( z\left(n, s\right) \right)}
  \end{align*}
  i.e.~the probability to go from $u$ to $t$ at time $n$ is proportional to the weight $W_t\left( z\left(n, t\right) \right)$ associated
  to the node $t$, which depends on the number of visits to the node. We set
  \begin{align*}
    w\left(n, t\right) := W_t\left( z\left(n, t\right) \right)
  \end{align*}
\end{definition}

\clearpage
\newpage

\section{Selected Results}
\label{sec:selected_results}

Here, we present known results on reinforced random walks to put the following sections into context.

\subsection{Reinforced Random Walks}

\subsubsection{Linearly Edge-Reinforced Random Walk}

The following relation between LERRW and RWRE is a very important tool to analyze the
LERRW. Many previous results on the LERRW have been proved using the fact that the
LERRW is equal in distribution to a mixture of Markov chains, i.e.~a MC with
random transition probabilities. MCs are well understood, thus this theorem allows us to use the
tools we have for MCs, and apply them to the non-Markovian LERRW.

\begin{theorem}[Reinforced Walk as Mixture of Markov Chains \textnormal{(\cite[Theorem 4]{localizationerrw}, \cite[Theorem 2.2]{rwforerrw}, \cite{definettimarkovchain})}]
  \label{thm:lerrw_is_rwre}
  For a finite graph $G$ with initial weights $\left(a_e\right)_{e \in E}$, the LERRW
  is equal in law to a RWRE with a uniquely determined mixing measure $\ensuremath{\bP}$.
  For any $e \in E$, $\bPrb{c_e > 0} = 1$. In particular, if the probability measure
  corresponding to the LERRW is denoted by $\ensuremath{\bbP_v}$, then
  \begin{align*}
    \Prbc{v}{\cdot} = \int \bbP_v^\mathbf{c}\left[\cdot\right] \mudx{\bP}{\mathbf{c}}
  \end{align*}
  with the notation from \autoref{def:rwre}.
\end{theorem}

We next turn to the LERRW on trees. Using \autoref{thm:lerrw_is_rwre}, one can show that

\begin{theorem}[Phase Transition of Reinforced Walk on Trees \textnormal{(\cite{errwpemantle,lerrwbsc})}]
  \label{thm:errw_trees}
  Let $T$ be a an infinite tree. Consider the LERRW on
  $T$ with initial weights $a > 0$. Let $\Delta_0$ be the solution to the equation
  \begin{align*}
    \frac{\Gamma\left(\frac{2 + \Delta}{4\Delta}\right)^2}{\Gamma\left(\frac{1}{2\Delta}\right)\Gamma\left(\frac{1 + \Delta}{2\Delta}\right)} = \frac{1}{\ensuremath{\br{T}}}
  \end{align*}
  where $\br{T}$ is the branching number of the tree $T$ (a measure of the number of children per node, the branching number is equal to the number
  of children for regular trees; see \cite[Definition 2.6]{lerrwbsc} for the exact definition). Then the LERRW on $T$ is
  \begin{enumerate}[(i)]
    \item a.s.~transient if $a > \frac{1}{\Delta_0}$
    \item a.s.~recurrent if $a < \frac{1}{\Delta_0}$.
  \end{enumerate}
  The expected return time to the root is always infinite.

  The same holds for Galton-Watson trees with $\br{T}$ replaced by the mean number of children, if the mean is
  $>1$ and conditioned upon non-extinction.

  On $\bbZ$, the LERRW is recurrent for any choice of initial weight $a$.
\end{theorem}

In this thesis, a recurrence result for the LERRW on general graphs will be shown. It is a part of the
following theorem, for which we also need to define the Cheeger constant. A graph is called nonamenable if there is
a constant $h > 0$ such that for any finite set of nodes $A \subseteq V$, $\left| \partial A \right| \geq h \left|A\right|$
where $\partial A = \left\{ v \in V \setminus A : \dist{v}{A} = 1 \right\}$, i.e.~$\partial A$ is the set of nodes which
are outside of $A$, but direct neighbors to a node within $A$. The Cheeger constant is the largest such constant $h$, and can be defined as follows:
\begin{align*}
  h\left(G\right) := \inf\left\{ \frac{\left| \partial A \right|}{\left|A\right|} : A \subseteq V \textrm{ finite} \right\}
\end{align*}

\begin{theorem}[LERRW on General Graphs \textnormal{(\cite{localizationerrw,errwandvrjp,transienceerrw})}]
  Let $G$ be any graph. Consider the LERRW with initial edge weight $a$ for every edge. Then:
  \begin{enumerate}[(i)]
    \item \label{thm:lerrw_general_rec} for every constant $K$, there exists some $a_0 > 0$ such that if all degrees of $G$ are bounded by $K$,
    the LERRW on $G$ is recurrent if $a \in \left(0, a_0\right)$.
    \item for every constant $K$, and any $c_0 > 0$, there exists some $a_0 > 0$ such that if all degrees of $G$ are bounded by $K$,
    and if $h\left(G\right) \geq c_0$, then the LERRW on $G$ is transient if $a > a_0$.
    \item for $G = \bbZ^3$ and $G = \bbZ^d$ with any $d \geq 3$, there exists some $a_0 > 0$ such that
    the LERRW on $G$ is transient if $a > a_0$.
  \end{enumerate}
\end{theorem}

Part \ref{thm:lerrw_general_rec} of the theorem above is shown in this thesis in \autoref{thm:rec_bndd_deg}.
Next to the above, and next to other results which are not listed here, the following is known about the LERRW:
\begin{itemize}
  \item On $\bbZ^3$ and $\bbZ^d$ with any $d \geq 3$, there is a sharp phase transition between recurrence and transience
  of the LERRW in the initial edge weights. If the initial edge weights are all set to a constant $a > 0$, then
  there is a critical threshold $a_0$ such that the LERRW on $\bbZ^d$ is recurrent if $a < a_0$ and transient if
  $a > a_0$. See \cite[Theorem 3]{monphasetrans}.
  \item On any locally finite graph, if a random walk with given conductances on the edges is recurrent,
  then so is the LERRW if these conductances are used as initial edge weights. See \cite[Theorem 4]{monphasetrans}.
  \item The VRJP (vertex-reinforced jump process) is a continuous-time version of the LERRW.
  The process $\left(Y_t\right)_{t \geq 0}$ takes values in the vertices of a graph with conductances $W_e$ on the edges.
  If we are at vertex $v$ at time $t$, then, conditional on $\left(Y_s, s \leq t\right)$, the process jumps to a neighbor
  $u$ of $v$ at rate $W_{\left\{v,u\right\}} \cdot L_u\left(t\right)$ with
  \begin{align*}
    L_u\left(t\right) = 1 + \int_0^t \mathbbm{1}_{Y_s = u} \dx{s}
  \end{align*}
  Up to an exponential time change, and using independent Gamma-distributed conductances, this has the same
  distribution as the following model of a continuous-time LERRW when looked at at jump times.
  In order to make the LERRW continuous, we can add clocks on every edge.
  The time of an edge runs only when the process is in a node adjacent to it, and the alarm of an edge rings at
  exponential times with the rate of the edge weight (i.e.~the rates are increasing). When an alarm rings, then
  the process immediately traverses the corresponding edge.

  The VRJP is distributed as a mixture of time-changed Markov jump processes, just as the LERRW is
  distributed as a mixture of Markov chains.

  The results on the LERRW basically also apply to the VRJP (or the other way round). See \cite{errwandvrjp}
  for more details.
  \item For the LERRW on finite graphs (with arbitrary initial weights), there is a  ``magic formula'':
  The LERRW is a mixture of Markov chains where the mixing measure can be given explicitly: if the edge weights
  are normalized to sum to one, they converge to a random limit on the unit simplex. The law can be given as a density
  w.r.t.~the Lebesgue measure on the unit simplex, and the LERRW is distributed as a mixture of (reversible) MCs
  with the weights chosen according to this same measure (the weights are not independent across
  edges). See \cite{magicformula}.

  The main reason why this is true and also the basis of the proof is partial exchangeability, i.e.~two finite paths
  have the same probability if the edges are traversed the same number of times (this is easy to see: the probability of
  a path is a product of fractions; for nominators, consider the edges, for denominators, the vertices, to see that the
  order of traversals is irrelevant).

  For recurrent LERRWs in general, we can use a de Finetti theorem on MCs, by considering blocks of
  excursions from the starting vertex (these are iid conditioned on some parameter) instead of balls drawn as for
  the most basic de Finetti theorem. See \cite{definettimarkovchain}.
  \item The RWRE representation for the LERRW also exists for infinite graphs and in
  particular transient LERRW (no explicit formula for the mixing measure is known). See \cite{rwforerrw}.
\end{itemize}
Just in order to mention that there is more than just linear reinforcement: there is stronger than linear reinforcement which can lead to being stuck on a single edge, and there
is the once-reinforced random walk (where the edge weight is only increased on the first visit and then stays
constant) which is even more difficult to analyze, little is known about it.

This thesis adapts the LERRW and adds multiple walkers or a bias to the model. \autoref{sec:multiple_walkers}
and \autoref{sec:biased} add some results to the list above.

\subsubsection{Vertex-Reinforced Random Walk}

There are some changes in behavior if vertex weights instead of edge weights are reinforced. A short overview of known results:
\begin{itemize}
  \item On any complete graph, the linearly vertex-reinforced random walk converges to uniform vertex occupation
  (meaning that the time spent at each vertex, divided by the total time, will converge to a uniform distribution,
  so the vertex weights will all be of the same order).
  On any finite graph fulfilling a non-degeneracy condition, vertex occupation converges to a rest point
  of a differential equation approximating the evolution of the vertex weights. See \cite{vrrwfinitegraphs}.
  \item On $\bbZ$, the VRRW (linear reinforcement) is eventually trapped on an interval of $5$ vertices,
  the middle vertex is visited with frequency $\frac{1}{2}$, its neighbors with positive frequency, and the two
  outer vertices are visited infinitely often but at frequency $0$. See \cite{vrrwzstuckonfive}.
  \item On $\bbZ^2$, the VRRW (linear reinforcement) is trapped with positive probability on $13$ vertices
  (a diamond with $\frac{1}{2}$-frequency core and with $0$-frequency boundary) and with positive probability on $12$
  vertices (square with $0$-frequency boundary). See \cite{vrrwarbitrary}. It still remains an open question whether the
  VRRW is trapped with probability $1$.
  \item On any graph of bounded degree, the VRRW (linear reinforcement) is trapped with positive probability
  on finitely many vertices which form a complete $n$-partite graph with outer boundary such that every vertex in the
  outer boundary is not connected to one of the partitions plus one extra node. See \cite{vrrwarbitrary}.
  \item On any tree of bounded degree, the VRRW (linear reinforcement) is trapped with probability $1$ on a finite subgraph. See \cite{vrrwarbitrary}.
  \item On a tree with $K_n$ children for every node on level $n$, where $\sum K_n^{-1} < \infty$, there is a positive
  probability for the VRRW (linear reinforcement) to move away from the root on \emph{every single} step (transience), but also a positive probability
  to get trapped (finite range), so there is no $0$-$1$-law. See \cite{vrrwarbitrary}.
\end{itemize}

These results show that an apparently small change in the model can dramatically affect the behavior of the random walk. Instead
of being recurrent or transient, as the LERRW, the VRRW with linear reinforcement has finite range (gets stuck
on a finite subgraph) with positive probability. This is one of the reasons why it was interesting to study variations to the model of
the ERRW, in order to better understand which changes in the model cause which effects.

\subsection{Urn Models}

Urn models are closely related to reinforced random walks since they are also an example of a reinforced
random process: usually, balls are drawn uniformly at random, and more balls of the drawn color are added.
Hence, the fraction of balls of a certain color is reinforced upon drawing a ball of the same color.

\subsubsection{P\'olya Urn}

Consider a P\'olya urn, starting with $w$ white and $b$ black balls. At every time step, draw a ball from the urn,
and put it back along with $\Delta$ more balls of the same color. Call $w_n$ and $b_n$ the number of white
and black balls after the $n$-th draw, so $w_0 = w$ and $b_0 = b$. Then:

\begin{theorem}[P\'olya Urn \textnormal{(\cite[Lemma 1]{errwpemantle})}]
  \label{thm:polya}
  The fraction of white balls $\frac{w_n}{w_n + b_n}$ is a martingale, and it thus converges to a limit random variable
  $L$ which is distributed according to a Beta distribution:
  \begin{align*}
    L \sim \Betad{\frac{w}{\Delta}}{\frac{b}{\Delta}}
  \end{align*}
  The sequence of balls drawn is equal in law to a mixture of iid sequences. Indeed, conditioned on $L$,
  the sequence of balls drawn is distributed as an iid sequence where a white ball is drawn with probability $L$.
\end{theorem}

The fact that the sequence of draws is distributed as a mixture of iid sequences follows from de Finetti's
theorem, which is applicable since the draws are exchangeable: the probability of a certain sequence of draws occurring
depends only on the number of white and black balls in the sequence, but not on their order. This is very similar to what
we already saw in \autoref{thm:lerrw_is_rwre}, and indeed, the same arguments are used in the proof.
Another thing to note is that \autoref{thm:polya} does not only hold when $w, b$ and $\Delta$ are integers, but also when they are replaced by any positive
real numbers. In other words, it is also possible to have only a fraction of a ball in the urn.

In \autoref{sec:multiple_walkers}, \autoref{cor:two_pl_urn_conv_rand} and \autoref{conj:two_pl_urn_limit_rand} show that something similar to a two-player
P\'olya urn (where two players draw balls and do not put them back immediately) exhibits very similar properties to that
of the standard one-player P\'olya urn.

\subsubsection{Randomly Reinforced Urn}

Consider the following urn model: initially, the urn contains $b_0$ black and $w_0$ white balls.
The numbers $b_n$ and $w_n$ of balls contained in the urn in the $n$-th step is determined as follows.
Draw a ball uniformly at random from the urn with $b_{n-1}$ black balls and $w_{n-1}$ white balls. If the ball is black, replace it with $M_n$ black
balls. If the ball is white, replace it with $N_n$ white balls. $M_n$ and $N_n$ are series of
independent random variables on the nonnegative natural numbers and all $M_n$ have the common distribution
$\mu$, all $N_n$ have the common distribution $\nu$. It is assumed that $\mu$ and $\nu$ both have finite support.
Call $z_n = \frac{b_n}{b_n + w_n}$ the fraction of black balls in the urn. Denote by $B_n, W_n, Z_n$
the corresponding random variables.

\begin{theorem}[Randomly Reinforced Urn \textnormal{(\cite[Theorem 4.1]{rreinfurn})}]
  \label{thm:urnconv}
  If the following assumptions hold (and \textbf{$\mu$ and $\nu$ both have finite support}, as stated above)
  \begin{enumerate}[(i)]
    \item $b_0 > 0 \quad$ (at least one black ball in the urn at the beginning)
    \item $\mu\left(\left\{0\right\}\right) = 0 \quad$ (if a black ball is drawn, it is replaced by at least one black ball)
    \item $\mu \geq_{\textrm{st}} \nu \quad$ ($\mu$ stochastically dominates $\nu$)
    \item $\Ex{M_n} = \int_{\bbN} x \;\mu\left(\textrm{d}x\right) > \int_{\bbN} x \;\nu\left(\textrm{d}x\right) = \Ex{N_n} \quad$
    (the expected number of black balls added is bigger than the expected number of white balls added)
  \end{enumerate}
  Then $\lim_{n \to \infty} Z_n = \lim_{n \to \infty} \frac{B_n}{B_n + W_n} = 1$ a.s. $\quad$
  (the fraction of black balls converges almost surely to $1$).
\end{theorem}

\autoref{thm:urnconv} covers a large class of urns, including urns with deterministic reinforcement
where a different number of balls is added depending on the color.

\clearpage
\newpage

\section{Recurrence of Edge-Reinforced Random Walk}
\label{sec:recurrence_errw}

In this section, we repeat the proof of the known result that the LERRW is
recurrent on any bounded-degree graph for small initial weights. We already know that the LERRW
has an equivalent representation as a RWRE with random conductances $c_e$ on the edges. If
we further normalize these conductances by setting $c_{v_0} = 1$ where $v_0$ is
the starting vertex of the LERRW (we will assume this normalization in the following sections),
then we can give a more precise statement of what we want
to show. We consider the weight function $W\left(n\right) = a + n$ and are
interested in the following result, given in \cite[Theorem 2]{localizationerrw}:
if $G$ is a graph with degree at most $K$, if the LERRW is recurrent on $G$ for any choice
of initial weights, and if $s \in \left(0, \frac{1}{4}\right)$,
then there is some $a_0 > 0$, depending on $s$ and $K$, such that for initial weights $a \in \left(0, a_0\right)$
we have $\Exc{v_0}{c_e^s} \leq 2K \left( C\left(s, K\right)\sqrt{a} \right)^{\dist{v_0}{e}}$
where $C\left(s, K\right)$ solely depends on $K$ and $s$. In particular, for $a$ small enough, we have exponential
decay in the expected conductances $c_e$ of the random environment corresponding to the LERRW
(exponential decay with increasing distance, and with the conductances raised to the $s$-th power).
As a consequence, we can prove under almost the same assumptions, but without requiring a priori
that the random walk is recurrent, that the LERRW is recurrent for $a$ small enough.

The proof given here will follow the proof given in \cite[Section 2]{localizationerrw} with only minor
modifications, which hopefully make a proof which was already really well presented in \cite{localizationerrw}
even easier to understand. So the mathematics presented here was really done by \cite{localizationerrw}.
The main idea is to estimate the ratios $\frac{c_e}{c_f}$ for edges $e$ and $f$ which
are incident on the same vertex by comparing with a corresponding ratio of numbers of edge traversals
of the LERRW, under the assumption that the random walk is recurrent. We then bound the
estimated ratio as well as the error caused by the estimation to conclude.

\subsection{Estimating Conductance Ratios}

The first step is to construct an estimate of conductance ratios. Our final goal will be to bound
the expectation of $c_e^s$ for any edge $e$. We therefore choose an arbitrary directed edge $e \in \ensuremath{\protect\overrightarrow{E}}$
to begin with. Next, we will define a path from the starting vertex $v_0$ to $e$ and estimate the conductance
ratios along this path, which will result in an estimate for $c_e$. The construction of the path is as follows:
for the edge $e = \left(\check{e}, \hat{e}\right)$, we can find the directed edge $e'$ through which
$\check{e}$ was first reached by the LERRW (or the RWRE, depending on your point of view).
If $e$ is traversed before its corresponding inverse $\overleftarrow{e}$, then $e' \neq e$ as undirected edges. Since
$e'$ must have been traversed before its inverse, by its definition, it is possible to iterate this
construction until the starting vertex $v_0$ is reached. This results in a (random) path $\gamma_e = \left( \ldots, e'', e', e \right)$
starting from $v_0$ and ending with $e$ (the path can be a loop if $\hat{e} = v_0$). For a deterministic
path $\gamma$, we call $D_{\gamma}$ the event that for the last edge $e$ of $\gamma$, we have $\gamma_e = \gamma$,
and that $e$ is traversed before $\overleftarrow{e}$.

For $e$ such that $e' \neq \overleftarrow{e}$, we now estimate the ratio $\frac{c_e}{c_{e'}}$ (i.e.~the
ratio of the conductances of two edges along one of the paths constructed above) by
\begin{align}
  \label{eq:est_ratios}
  Q\left(e\right) &:= \frac{M_e}{M_f} \qquad \textrm{where } f := \overleftarrow{e'},
  \qquad \textrm{an estimate for } R\left(e\right) := \frac{c_e}{c_{e'}}
\end{align}
where $M_e$ counts the number of times the directed edge $e$ was crossed in the right direction before
the edge $f$ was crossed in the right direction by the random walk, if $e$ is crossed before $f$. In this
case, $M_f$ is set to $1$. Otherwise, $M_f$ counts the number of times $f$ was crossed before $e$, and
$M_e$ is $1$. In other words, $M_e$ and $M_f$ count the numbers of directed edge crossings of $e$ and $f$ until both
edges have been traversed by the random walk. Since both edges lead out of the same vertex, this can also
be seen as counting the number of departures from $\check{e}$ until both edges were used.
\begin{figure}[H]
  \begin{center}
    \begin{tikzpicture}
      \draw[tumGray] (0, 0) -- (1.4, 1.4);
      \draw[tumGray] (2, 0) -- (1.4, 1.4);
      \draw[tumGray] (2, 0) -- (3.5, 1);
      \draw[tumGray] (-3.9, -0.2) -- (-2, 0);
      \node[circle,fill=black,inner sep=0.7mm] (N0) at (0, 0) {};
      \node[circle,fill=black,inner sep=0.7mm] (N1) at (2, 0) {};
      \node[circle,fill=black,inner sep=0.7mm] (Nm1) at (-2, 0) {};
      \node[circle,fill=tumGray,inner sep=0.7mm] (Nm2) at (-3.5, 1) {};
      \node[circle,fill=tumGray,inner sep=0.7mm] (N1b) at (1.4, 1.4) {};
      \node[circle,fill=tumGray,inner sep=0.7mm] (N2) at (3.5, 1) {};
      \node[circle,fill=tumGray,inner sep=0.7mm] (Nm2b) at (-3.9, -0.2) {};
      \node[above=1mm,tumGray] at (Nm2) {$v_0$};
      \node[above=1mm] at (N0) {$\check{e}$};
      \node[above=1mm] at (N1) {$\hat{e}$};
      \draw (N0) edge[-{Latex[length=2mm,width=2mm]},sectionblue,line width=0.5mm] node[above,sectionblue] {$e$} (N1);
      \draw (Nm1.north) edge[-{Latex[length=2mm,width=2mm]},line width=0.5mm] node[above] {$e'$} (N0.north);
      \draw (N0.south) edge[-{Latex[length=2mm,width=2mm]},tumOrange] node[below,tumOrange] {$f$} (Nm1.south);
      \draw (Nm2) edge[-{Latex[length=2mm,width=2mm]},tumGray,line width=0.5mm] node[above=1mm,tumGray] {$e''$} (Nm1);
      \node at (-2.2, -1) {$\textcolor{tumOrange}{M_f}\;=$ number of};
      \node at (-2.2, -1.5) {departures from $\check{e}$ along $\textcolor{tumOrange}{f}$};
      \node at (-2.2, -2) {until both $e$ and $f$ used};
      \node at (2.2, -1) {$\textcolor{sectionblue}{M_e}\;=$ number of};
      \node at (2.2, -1.5) {departures from $\check{e}$ along $\textcolor{sectionblue}{e}$};
      \node at (2.2, -2) {until both $e$ and $f$ used};
      \node at (5.5, 0) {$\displaystyle{Q\left(e\right) = \frac{\textcolor{sectionblue}{M_e}}{\textcolor{tumOrange}{M_f}}}$, an estimate for $\displaystyle{R\left(e\right) = \frac{c_e}{c_{e'}}}$};
      \node at (0, 2) {$\boldsymbol{\gamma}_{\boldsymbol{e}} = \left(e'', e', e\right)$};
    \end{tikzpicture}
    \caption{Estimating conductance ratios}
    \label{fig:ratio_estimate_q}
  \end{center}
\end{figure}
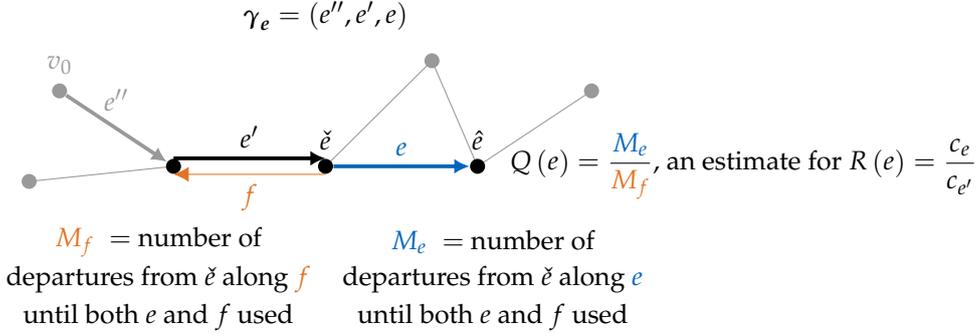

We next define the deterministic set of simple paths (or loops, if one of the endpoints of $e$ is $v_0$)
$\Gamma_e$ as the set of paths which end with $e$ or $\overleftarrow{e}$. So, $\Gamma_e$ is the set
of possible values which the random path $\gamma_e$ (or $\gamma_{\overleftarrow{e}}$, if $\overleftarrow{e}$
is traversed before $e$) can take. We now have
\begin{align*}
  \Exc{v_0}{c_e^s} = \sum_{\gamma \in \Gamma_e} \Exc{v_0}{c_e^s \cdot \indic{D_{\gamma}}}
  \leq \sum_{\gamma \in \Gamma_e} \Exc{v_0}{\left(\frac{c_e}{c_{\gamma_1}}\right)^s \cdot \indic{D_{\gamma}}}
\end{align*}
where $\gamma_1$ is the first edge of $\gamma$. Recall that we normalized the environment such that
$c_{v_0} = 1$, so in particular, $c_{\gamma_1} \leq 1$ since $\gamma_1$ leads out of $v_0$.

One of the assumptions we started with was that we already know that the LERRW is recurrent.
This implies that $e$ will be visited a.s.~and therefore the partition into possible paths
$\gamma \in \Gamma_e$ does indeed cover the whole probability space. Furthermore, by \autoref{thm:lerrw_is_rwre} and \cite{definettimarkovchain},
the conductances in the random environment corresponding to the LERRW are uniquely defined
and we can therefore look at $\frac{c_e}{c_{\gamma_1}}$.

We now want to compare the estimated conductance ratios along the path $\gamma_e$ with the real
conductance ratios. If we fix a path $\gamma$ ending in $e$ and $D_{\gamma}$ occurs, then
\begin{align*}
  \frac{c_e}{c_{\gamma_1}} = \prod_{f \in \gamma, f \neq \gamma_1} R\left(f\right)
  = \underbrace{\prod_{f \in \gamma, f \neq \gamma_1} \frac{R\left(f\right)}{Q\left(f\right)}}_{\textrm{error}} \underbrace{\prod_{f \in \gamma, f \neq \gamma_1} Q\left(f\right)}_{\textrm{estimate}}
\end{align*}
By the Cauchy-Schwarz inequality, we have
\begin{align*}
  \Exc{v_0}{\left(\frac{c_e}{c_{\gamma_1}}\right)^s \cdot \indic{D_{\gamma}}} \leq
  \Exc{v_0}{\prod_{f \in \gamma, f \neq \gamma_1} \left(\frac{R\left(f\right)}{Q\left(f\right)}\right)^{2s} \cdot \indic{D_{\gamma}}}^{\frac{1}{2}}
  \Exc{v_0}{\prod_{f \in \gamma, f \neq \gamma_1} Q\left(f\right)^{2s} \cdot \indic{D_{\gamma}}}^{\frac{1}{2}}
\end{align*}

The above equation allows us to bound the expectation of the error introduced by the estimate separately from the
expectation of the estimate. The following two sections will bound the two terms on the right hand side. Then, we will
be able to show exponential decay of the edge conductances in the distance to the starting vertex for graphs on which
the LERRW is recurrent, a result which can subsequently actually be used to show recurrence of the
LERRW itself.

\subsection{Bounding the Error}

\begin{lemma}[Error Bound \textnormal{(\cite[Lemma 7]{localizationerrw})}]
  \label{lem:error_bound}
  Let $G$ be a graph, $v_0 \in V$ the starting vertex, and $a \in \left(0, \infty\right)$ such that the
  LERRW on $G$ with initial weights equal to $a$ is recurrent. Then, for any $e \in E$,
  any $\gamma \in \Gamma_e$, and any $s \in \left(0, 1\right)$:
  \begin{align*}
    \Exc{v_0}{\prod_{f \in \gamma, f \neq \gamma_1} \left(\frac{R\left(f\right)}{Q\left(f\right)}\right)^{s} \cdot \indic{D_{\gamma}}} \leq C\left(s\right)^{\left|\gamma\right| - 1}
  \end{align*}
  where $C\left(s\right)$ is a constant depending solely on $s$ and $\left|\gamma\right|$ is the length of the path $\gamma$.
\end{lemma}

Note that the bound of the error given above shows that we have at most an exponential increase of
our error in the length of the path $\gamma$. The corresponding bound of the estimated ratios will
balance this potentially exponentially increasing error.

\begin{tproof}
  $R$ and $Q$ depend in general on the random path $\gamma_e$ (which was defined recursively as the
  path of first entrance edges in the section above). For a deterministic path $\gamma$, we can define
  $R_{\gamma}$ and $Q_{\gamma}$ for edges $f \in \gamma, f \neq \gamma_1$ as we defined $R$ and $Q$
  for the random path: if $f \in \gamma$, and if $f' \in \gamma$ is the preceding edge in the path, then
  $R_{\gamma}\left(f\right) = c_f c_{f'}^{-1}$, and $Q_{\gamma}\left(f\right)$ is the ratio of exits from
  $\check{f}$ along $f$ and $\overleftarrow{f'}$ until both edges have been visited (see \autoref{eq:est_ratios}).
  We have, for fixed $\gamma$
  \begin{align*}
    \Exc{v_0}{\prod_{f \in \gamma, f \neq \gamma_1} \left(\frac{R\left(f\right)}{Q\left(f\right)}\right)^{s} \cdot \indic{D_{\gamma}}}
    \leq
    \Exc{v_0}{\prod_{f \in \gamma, f \neq \gamma_1} \left(\frac{R_{\gamma}\left(f\right)}{Q_{\gamma}\left(f\right)}\right)^{s}}
  \end{align*}
  This inequality holds since $R_{\gamma}\left(f\right) = R\left(f\right)$ and $Q_{\gamma}\left(f\right) = Q\left(f\right)$
  if $D_{\gamma}$ occurs. The inequality is a crude bound but sufficient for our purposes.

  We want to show a bound on the error of estimating the conductance ratios. It is therefore necessary
  to look at the RWRE description of the random walk, and we will now condition on a particular
  vector $\ensuremath{\mathbf{c}}$ of conductances being chosen. We then show that for every such vector, the
  bound above holds, and this will give the desired result. So from now on, assume we have conductances
  $\mathbf{c}$. $R_{\gamma}\left(f\right)$ is then also a fixed number, since it only depends on the conductances.
  The $Q_{\gamma}\left(f\right)$, on the other hand, are independent for any two different edges if
  we condition on the conductances $\mathbf{c}$. This is because we just count directed exits from
  $\check{f}$, and we can couple these with iid sequences of random variables $\left(Z_n^v\right)_{n \geq 1}$
  for every node $v$ where $Z_n^v$ takes values in the edges incident to node $v$, and is distributed
  such that an edge $e$ appears with probability $c_e \left( \sum_{f \ni v} c_f \right)^{-1}$. These
  $Z_n^v$ are iid for every node, and independent for different nodes.

  By the independence of the $Q_{\gamma}\left(f\right)$ (conditioned on $\mathbf{c}$), we have
  \begin{align*}
    \Exc{v_0}{\prod_{f \in \gamma, f \neq \gamma_1} \left(\frac{R_{\gamma}\left(f\right)}{Q_{\gamma}\left(f\right)}\right)^{s} \midd| \mathbf{c}}
    &= \prod_{f \in \gamma, f \neq \gamma_1} \Exc{v_0}{\left(\frac{R_{\gamma}\left(f\right)}{Q_{\gamma}\left(f\right)}\right)^{s} \midd| \mathbf{c}} \\
    &= \prod_{f \in \gamma, f \neq \gamma_1} \Exc{v_0}{\left(\frac{c_f}{c_{f'}} \cdot \frac{M_{\overleftarrow{f'},\gamma}}{M_{f,\gamma}}\right)^{s} \midd| \mathbf{c}}
  \end{align*}
  where $f'$ is the edge preceding $f$ in $\gamma$ and $M_{f,\gamma}$ is defined as $M_f$ in \autoref{eq:est_ratios}
  but with the deterministic path $\gamma$ instead of the random path $\gamma_e$. To finish the proof,
  it is therefore sufficient to show that, for any edge $f \in \gamma$ preceded by $f'$,
  \begin{align*}
    \Exc{v_0}{\left(\frac{c_f}{c_{f'}} \cdot \frac{M_{\overleftarrow{f'},\gamma}}{M_{f,\gamma}}\right)^{s} \midd| \mathbf{c}} \leq C\left(s\right)
  \end{align*}
  where $C\left(s\right)$ does not depend on the conductance vector $\mathbf{c}$ we condition on.

  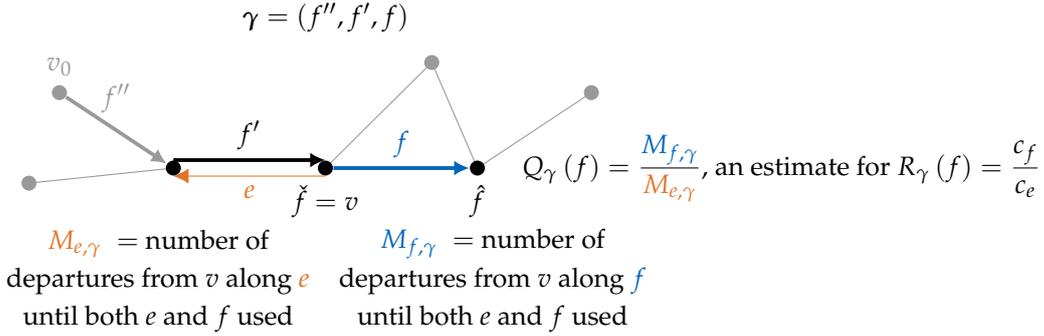
\begin{figure}[H]
    \begin{center}
      \begin{tikzpicture}
        \draw[tumGray] (0, 0) -- (1.4, 1.4);
        \draw[tumGray] (2, 0) -- (1.4, 1.4);
        \draw[tumGray] (2, 0) -- (3.5, 1);
        \draw[tumGray] (-3.9, -0.2) -- (-2, 0);
        \node[circle,fill=black,inner sep=0.7mm] (N0) at (0, 0) {};
        \node[circle,fill=black,inner sep=0.7mm] (N1) at (2, 0) {};
        \node[circle,fill=black,inner sep=0.7mm] (Nm1) at (-2, 0) {};
        \node[circle,fill=tumGray,inner sep=0.7mm] (Nm2) at (-3.5, 1) {};
        \node[circle,fill=tumGray,inner sep=0.7mm] (N1b) at (1.4, 1.4) {};
        \node[circle,fill=tumGray,inner sep=0.7mm] (N2) at (3.5, 1) {};
        \node[circle,fill=tumGray,inner sep=0.7mm] (Nm2b) at (-3.9, -0.2) {};
        \node[above=1mm,tumGray] at (Nm2) {$v_0$};
        \node[below=1mm] at (N0) {$\check{f} = v$};
        \node[below=1mm] at (N1) {$\hat{f}$};
        \draw (N0) edge[-{Latex[length=2mm,width=2mm]},sectionblue,line width=0.5mm] node[above,sectionblue] {$f$} (N1);
        \draw (Nm1.north) edge[-{Latex[length=2mm,width=2mm]},line width=0.5mm] node[above] {$f'$} (N0.north);
        \draw (N0.south) edge[-{Latex[length=2mm,width=2mm]},tumOrange] node[below,tumOrange] {$e$} (Nm1.south);
        \draw (Nm2) edge[-{Latex[length=2mm,width=2mm]},tumGray,line width=0.5mm] node[above=1mm,tumGray] {$f''$} (Nm1);
        \node at (-2.2, -1) {$\textcolor{tumOrange}{M_{e,\gamma}}\;=$ number of};
        \node at (-2.2, -1.5) {departures from $v$ along $\textcolor{tumOrange}{e}$};
        \node at (-2.2, -2) {until both $e$ and $f$ used};
        \node at (2.2, -1) {$\textcolor{sectionblue}{M_{f,\gamma}}\;=$ number of};
        \node at (2.2, -1.5) {departures from $v$ along $\textcolor{sectionblue}{f}$};
        \node at (2.2, -2) {until both $e$ and $f$ used};
        \node at (6, 0) {$\displaystyle{Q_{\gamma}\left(f\right) = \frac{\textcolor{sectionblue}{M_{f,\gamma}}}{\textcolor{tumOrange}{M_{e,\gamma}}}}$, an estimate for $\displaystyle{R_{\gamma}\left(f\right) = \frac{c_f}{c_{e}}}$};
        \node at (0, 2) {$\boldsymbol{\gamma} = \left(f'', f', f\right)$};
      \end{tikzpicture}
      \caption{Estimating conductance ratios along a deterministic path $\gamma$}
      \label{fig:ratio_estimate_lemproof}
    \end{center}
  \end{figure}
  Write $v = \check{f}$, $e = \overleftarrow{f'}$ (see \autoref{fig:ratio_estimate_lemproof} for reference).
  The law of $M_{e,\gamma}$ and $M_{f,\gamma}$ is determined by the iid sequence $\left(Z_n^v\right)_{n \geq 1}$
  (since we condition on $\mathbf{c}$), and does not depend on other edges incident to $v$: it is only
  relevant with which probability the two edges appear relative to each other in the sequence. We can
  therefore assume w.l.o.g.~that only $e$ and $f$ appear in the sequence $Z_n^v$ with probabilities
  $p$ and $q = 1 - p$, respectively. Then
  \begin{align*}
    \Exc{v_0}{\left(\frac{c_f}{c_{e}} \cdot \frac{M_{e,\gamma}}{M_{f,\gamma}}\right)^{s} \midd| \mathbf{c}}
    &= \left(\frac{c_f}{c_e}\right)^s \left( \sum_{k \geq 1} k^s \cdot \Prbc{v_0}{k \textrm{ exits from } v \textrm{ along } e, \textrm{ then exit } v \textrm{ along } f} \right. \\
    &\hphantom{= \left(\frac{c_f}{c_e}\right)^s (}\left. + \sum_{k \geq 1} k^{-s} \cdot \Prbc{v_0}{k \textrm{ exits from } v \textrm{ along } f, \textrm{ then exit } v \textrm{ along } e} \right) \\
    &= \left(\frac{q}{p}\right)^s \left( \sum_{k \geq 1} k^s p^k q + \sum_{k \geq 1} k^{-s} q^k p \right) \\
    &\overset{\circledast}{\leq} \left(\frac{q}{p}\right)^s \left( q^{-s}p + C\left(s\right) qp^s \right) = p^{1-s} + C\left(s\right) q^{1+s} \leq 2 C\left(s\right)
  \end{align*}
  Since the bound is independent of $p$ and $q$ and thereby independent of the choice of the conductances $\mathbf{c}$,
  this completes the proof, up to showing $\circledast$. Note that the bound can be improved, dependent on the value of $s$.
  We show $\circledast$ using the following observations.
  If $X \sim \Geo{q}$ is a random variable with geometric distribution of parameter $q$, i.e.~$\Prb{X = k} = p^k q$, then we can write
  \begin{align}
    \label{eq:first_sum_bound}
    \Ex{X^s} = \sum_{k \geq 1} k^s p^k q \leq pq^{-s}
    &\iff \Ex{\left(qX\right)^s} \leq p
  \end{align}
  Now, proving that the inequality in \autoref{eq:first_sum_bound} holds will show the bound on the first sum in $\circledast$.
  Note that the function $f\left(s\right) = x^s$ is convex for arbitrary choice of $x > 0$.
  Therefore, $g\left(s\right) = \Ex{\left(qX\right)^s}$ is convex in $s$. As $X$ takes values in the integers, we have that
  $\Ex{\left(qX\right)^s} = q^s \Ex{X^s} \leq \Ex{X} < \infty$ for $s \in \left(0, 1\right)$.
  By dominated convergence, $g\left(s\right)$ is therefore continuous in $s$, and can be continued continuously
  to the interval $s \in \left[0, 1\right]$. For $s = 1$, we have $g\left(1\right) = q \Ex{X} = q \cdot \frac{p}{q} = p$,
  and for $s = 0$, we have
  \begin{align*}
    \lim_{s \to 0} g\left(s\right) = \lim_{s \to 0} \Ex{\left(qX\right)^s} = \Ex{\lim_{s \to 0} \left(qX\right)^s}
    = \Ex{\indic{\left\{X > 0\right\}}} = p
  \end{align*}
  By convexity of $g$, we conclude that $g\left(s\right) \leq p$ for all $s \in \left(0, 1\right)$, which proves \autoref{eq:first_sum_bound}.

  For bounding the second sum, let $X \sim \Geob{p}$ be a random variable with geometric distribution of parameter $p$, but which takes values
  in $\bbN \setminus \left\{0\right\}$, so $\Prb{X = k} = q^{k-1} p$ where $k \geq 1$. Then we can write
  \begin{align}
    \label{eq:scnd_sum_bound}
    q \cdot \Ex{X^{-s}} = q \cdot \sum_{k \geq 1} k^{-s} q^{k - 1} p = \sum_{k \geq 1} k^{-s} q^k p \leq C\left(s\right) qp^s
    &\iff \Ex{\left(pX\right)^{-s}} \leq C\left(s\right)
  \end{align}
  Now, $g\left(p\right) = \Ex{\left(pX\right)^{-s}}$ is decreasing in $p$ (to see this, note that
  $\Ex{\left(pX\right)^{-s}} = \frac{p^{1 - s}}{1 - p} \cdot \Polylog{s}{1 - p}$ where $\POLYLOG$
  denotes the polylogarithm, derive: \href{https://www.wolframalpha.com/input?i=d%2Fdx+of+x%5E%281-s%29+%2F+%281-x%29+*+PolyLog%5Bs%2C+1-x%5D}{https://bit.ly/3BOmLjP},
  and analyze the derivative: \href{https://www.wolframalpha.com/input?i=plot+%28s+*+%28x+-+1%29+%2B+1%29+*+PolyLog%5Bs%2C+1-x%5D+-+x+*+PolyLog%5Bs+-+1%2C+1-x%5D+for+x+from+0+to+1+and+s+from+0+to+1}{https://bit.ly/3RVn5CA}), so it is sufficient to show that
  $\lim_{p \to 0} g\left(p\right) =: C\left(s\right)$ exists and is finite. Note that for $p \to 0$,
  $pX$ converges in distribution to a random variable $Y$ which is $\Expd{1}$-distributed
  (we have that $\Prb{pX > t} \to e^{-t}$ for $p \to 0$). Therefore, $\lim_{p \to 0} g\left(p\right) = \Ex{Y^{-s}} = \Gamma\left(1 - s\right) < \infty$
  for $s \in \left(0, 1\right)$. Thus, we can choose $C\left(s\right) = \Gamma\left(1 - s\right)$ in \autoref{eq:scnd_sum_bound}.
\end{tproof}

\subsection{Bounding the Estimate}

\begin{lemma}[Estimate Bound \textnormal{(\cite[Lemma 8]{localizationerrw})}]
  \label{lem:estimate_bound}
  Let $G$ be a graph with degree at most $K$, $v_0 \in V$ the starting vertex, and $a \in \left(0, \infty\right)$ such that the
  LERRW on $G$ with initial weights equal to $a$ is recurrent. Then, for any $e \in E$,
  any $\gamma \in \Gamma_e$, and any $s \in \left(0, \frac{1}{2}\right)$:
  \begin{align*}
    \Exc{v_0}{\prod_{f \in \gamma, f \neq \gamma_1} Q\left(f\right)^{s} \cdot \indic{D_{\gamma}}} \leq \left(C\left(s, K\right)a\right)^{\left|\gamma\right| - 1}
  \end{align*}
  where $C\left(s\right)$ is a constant depending solely on $s$ and $K$.
\end{lemma}

Choosing $a$ small enough, the bound above shows that the estimated ratios decrease exponentially
and fast enough to balance the potential error caused by the estimation given in \autoref{lem:error_bound}.

\begin{tproof}
  The idea of the proof is to find iid random variables $\overline{Q}\left(f\right)$ which
  stochastically dominate the variables $Q\left(f\right)$ on all edges of the deterministic path $\gamma$,
  provided that $D_{\gamma}$ occurs. It is then much easier to bound the expectation of the random
  variables $\overline{Q}\left(f\right)$. We start by fixing $\gamma$ and now want to define
  $\overline{Q}\left(f\right)$ such that $Q\left(f\right) \leq \overline{Q}\left(f\right)$ for all
  $f \in \gamma, f \neq \gamma_1$ on the event $D_{\gamma}$.

  For every such edge $f$, we define two independent series of Bernoulli random variables (also
  independent of the corresponding series for different edges) $\left(Y_j\right)_{j \geq 0}$ and $\left(Y_j'\right)_{j \geq 0}$:
  \begin{align*}
    \Prb{Y_j = 1} = \frac{a}{j + 1 + 2a} = 1 - \Prb{Y_j = 0} \qquad
    \Prb{Y_j' = 1} = \frac{1 + a}{2j + 1 + Ka} = 1 - \Prb{Y_j' = 0}
  \end{align*}
  Intuitively, and up to some technical details, these variables are used as follows. Recall that
  $f'$ was the edge through which $\check{f} = v$ was first reached (which is the edge preceding $f$ in $\gamma$
  on the event $D_{\gamma}$) and we set again $e = \overleftarrow{f'}$.
  We want to bound $Q\left(f\right) = \frac{M_f}{M_e}$, hence we want to upper bound $M_f$, the number
  of exits along $f$ until both edges were used, and we want to lower bound $M_e$. First, $Y_0' = 1$
  represents the event that on the first visit to $v$, we depart along $e$, and sometimes also if
  $Y_0' = 0$. $Y_0' = 1$ is thus a lower bound for the departure along $e$. Now,
  \begin{itemize}
    \item if $Y_0' = 0$, i.e.~if the first departure was not necessarily along $e$, and $e$ was not used yet: then,
    for $n \geq 1$, on the $n$-th visit to $v$, we will couple the
    random walk with the variables $Y'$ such that we depart along $e$ if $Y_{n-1}' = 1$, and sometimes
    also if $Y_{n-1}' = 0$. $Y_{n-1}' = 1$ lower bounds the event that $e$ is crossed.
    \item if $Y_0' = 1$, i.e.~if the first departure was along $e$, and if $f$ was not used yet: then,
    for $n \geq 1$, on the $n$-th visit to $v$ where $v$ is left along $f$ or $e$, we will couple the
    random walk with the variables $Y$ such that we depart along $e$ if $Y_{n-1} = 0$, and sometimes
    also if $Y_{n-1} = 1$. $Y_{n-1} = 1$ upper bounds the event that $f$ is crossed.
  \end{itemize}
  We now set
  \begin{align*}
    \overline{Q}\left(f\right) &= \frac{\overline{M}_f}{\overline{M}_e} \quad \textrm{where} \\
    \overline{M}_f &= \min\left\{ j \geq 1 : Y_j' = 1 \right\} \;\; \textrm{ and } \;\; \overline{M}_e = 1 \;\; \textrm{ if } \;\; Y_0' = 0 \\
    \overline{M}_e &= \min\left\{ j \geq 1 : Y_j = 1 \right\} \;\; \textrm{ and } \;\; \overline{M}_f = 1 \;\; \textrm{ if } \;\; Y_0' = 1
  \end{align*}

  \begin{figure}[H]
    \begin{center}
      \begin{tikzpicture}
        \node at (0, 1.8) {First exit not necessarily along $e$, $Y_0'=0$};
        \draw[tumGray] (0, 0) -- (1.4, 1.4);
        \draw[tumGray] (0, 0) -- (-1.4, 1.4);
        \draw[tumGray] (0, 0) -- (0, -2);
        \node[circle,fill=black,inner sep=0.7mm] (N0) at (0, 0) {};
        \node[circle,fill=black,inner sep=0.7mm] (N1) at (2, 0) {};
        \node[circle,fill=black,inner sep=0.7mm] (Nm1) at (-2, 0) {};
        \node[circle,fill=tumGray,inner sep=0.7mm] (N1b) at (1.4, 1.4) {};
        \node[circle,fill=tumGray,inner sep=0.7mm] (N1c) at (-1.4, 1.4) {};
        \node[circle,fill=tumGray,inner sep=0.7mm] (N1d) at (0, -2) {};
        \node[below=1mm] at (N0) {$\check{f} = v$};
        \node[below=1mm] at (N1) {$\hat{f}$};
        \draw (N0) edge[-{Latex[length=2mm,width=2mm]},sectionblue] node[above] {$f$} (N1);
        \draw (Nm1.north) edge[-{Latex[length=2mm,width=2mm]},tumOrange] node[above] {$f'$} (N0.north);
        \draw (N0.south) edge[-{Latex[length=2mm,width=2mm]},tumOrange] node[below] {$e$} (Nm1.south);
        \node at (0, -2.4) {Weights at the $n$-th visit to $v$:};
        \node[tumOrange] at (-1, -2.9) {$\textcolor{tumOrange}{z\left(t, e\right)} + a \geq 1 + a$};
        \node[sectionblue] at (1, -2.9) {?};
        \node[white] at (1, -3.4) {$z\left(t, f\right) + a = a$};
        \node at (0, -3.9) {Total weight of adjacent edges:};
        \node at (0, -4.4) {$2n - 1 + \deg\left(v\right)a \leq 2n - 1 + Ka$};
        \node at (0, -5) {$\Prb{\textrm{exit along }e} \geq \frac{1 + a}{2n - 1 + Ka} = \Prb{Y_{n-1}'=1}$};
      \end{tikzpicture} ~~~~~~~~~~
      \begin{tikzpicture}
        \node at (0, 1.8) {First exit along $e$, $Y_0'=1$};
        \draw[tumGray] (0, 0) -- (1.4, 1.4);
        \draw[tumGray] (0, 0) -- (-1.4, 1.4);
        \draw[tumGray] (0, 0) -- (0, -2);
        \node[circle,fill=black,inner sep=0.7mm] (N0) at (0, 0) {};
        \node[circle,fill=black,inner sep=0.7mm] (N1) at (2, 0) {};
        \node[circle,fill=black,inner sep=0.7mm] (Nm1) at (-2, 0) {};
        \node[circle,fill=tumGray,inner sep=0.7mm] (N1b) at (1.4, 1.4) {};
        \node[circle,fill=tumGray,inner sep=0.7mm] (N1c) at (-1.4, 1.4) {};
        \node[circle,fill=tumGray,inner sep=0.7mm] (N1d) at (0, -2) {};
        \node[below=1mm] at (N0) {$\check{f} = v$};
        \node[below=1mm] at (N1) {$\hat{f}$};
        \draw (N0) edge[-{Latex[length=2mm,width=2mm]},sectionblue] node[above] {$f$} (N1);
        \draw (Nm1.north) edge[-{Latex[length=2mm,width=2mm]},tumOrange] node[above] {$f'$} (N0.north);
        \draw (N0.south) edge[-{Latex[length=2mm,width=2mm]},tumOrange] node[below] {$e$} (Nm1.south);
        \node at (0, -2.4) {Weights at the $n$-th time at which};
        \node at (0, -2.9) {$v$ is exited along $e$ or $f$:};
        \node[tumOrange] at (-1, -3.4) {?};
        \node[sectionblue] at (1, -3.4) {$z\left(t, f\right) + a = a$};
        \node at (0, -3.9) {Total weight of edges $e$ and $f$:};
        \node at (0, -4.4) {$z\left(t, e\right) + 2a \geq n + 2a$};
        \node at (0, -5) {$\Prb{\textrm{exit along }f} \leq \frac{a}{n + 2a} = \Prb{Y_{n-1}=1}$};
      \end{tikzpicture}
      \caption{Bounding $Q$ with the help of Bernoulli random variables}
      \label{fig:bounding_q}
    \end{center}
  \end{figure}
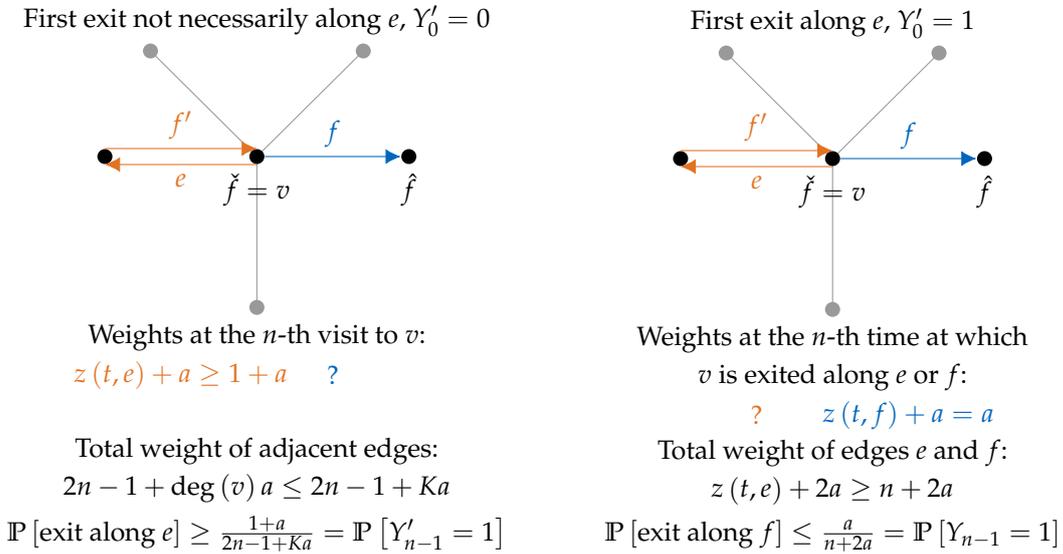
  The figure above demonstrates why the indicated coupling is possible. We are in the situation
  where $v$ has not yet been exited along both $e$ and $f$. On the left, the first exit from $v$
  was not necessarily along $e$, and we want to lower bound the probability to exit along $e$ on
  the $n$-th visit to $v$. $v$ was first reached by $f'$, so the weight of the edge $e$ must already
  be at least $1 + a$, which yields the desired lower bound of $\Prb{Y_{n-1}'=1}$. On the right, the
  first exit from $v$ was along $e$, and we want to upper bound the probability to exit along $f$,
  if we already know that we exit along one of the two edges $e$ and $f$. Since we look at the situation
  that not both edges have been used yet, $f$ cannot have been used yet (also not in reverse direction,
  at least not on the event $D_{\gamma}$), so its weight is still $a$. The weight of $e$, on the other
  hand must be at least $n + a$ if this is the $n$-th time at which we exit along one of these two
  edges (it could be higher, since $e$ might also have been used in the reverse direction). This
  yields the upper bound of $\Prb{Y_{n-1}=1}$.

  The preceding paragraphs hopefully gave enough details to intuitively understand how the random
  walk and the variables $Y$ and $Y'$ should be coupled, and why $\overline{Q}$ is an upper bound
  for $Q$. If not, the technical details will be given at the end of the proof. Assuming that we know
  that $Q\left(f\right) \leq \overline{Q}\left(f\right)$ for all $f \in \gamma, f \neq \gamma_1$
  on the event $D_{\gamma}$, we can now finish the prove as follows. Let $s \in \left(0, \frac{1}{2}\right)$.
  We have $\Ex{\overline{Q}\left(f\right)^s} = \Ex{\overline{Q}\left(f\right)^s \indic{Y_0' = 0}} + \Ex{\overline{Q}\left(f\right)^s \indic{Y_0' = 1}}$.

  If $Y_0' = 0$, then $\overline{M}_e = 1$. In addition,
  \begin{align*}
    \Prb{Y_0' = 0, \overline{M}_f = n} &= \Prb{Y_0' = 0} \cdot \Prb{Y_n' = 1} \cdot \prod_{j=1}^{n-1} \Prb{Y_j' = 0} \\
    &= \underbrace{\frac{\left(K-1\right)a}{1 + Ka}}_{\leq \frac{Ka}{1 + Ka}} \cdot \underbrace{\frac{1 + a}{2n + 1 + Ka}}_{\leq \frac{1 + Ka}{2n}} \cdot \prod_{j=1}^{n-1} \left(1 - \frac{1 + a}{2j + 1 + Ka}\right) \\
    &\leq \frac{Ka}{2n} \cdot \prod_{j=1}^{n-1} \left(1 - \frac{1 + a}{2j + 1 + Ka}\right)
  \end{align*}
  For the terms in the product, we have the following estimate, since $a > 0$:
  \begin{align*}
    \frac{1 + a}{2j + 1 + Ka} &\geq \min \left\{ \frac{1}{2j + 1}, \frac{a}{Ka} \right\} = \min \left\{ \frac{1}{2j + 1}, \frac{1}{K} \right\} \\
    \implies 1 - \frac{1 + a}{2j + 1 + Ka} &\leq \exp\left(- \frac{1 + a}{2j + 1 + Ka}\right) \leq \exp\left(- \min \left\{ \frac{1}{2j + 1}, \frac{1}{K} \right\}\right)
  \end{align*}
  Note that $- \frac{1}{2j + 1} = -\frac{1}{2j} + \frac{1}{j\left(4j + 2\right)} \leq -\frac{1}{2j} + \frac{1}{4j^2}$.
  Therefore,
  \begin{align*}
    \prod_{j=1}^{n-1} \left(1 - \frac{1 + a}{2j + 1 + Ka}\right) &\leq \exp\left(- \sum_{j = 1}^{n - 1} \frac{1}{2j} + \sum_{j = 1}^{n - 1} \frac{1}{4j^2} + C_1\left(K\right)\right)
  \end{align*}
  where $0 < C_1\left(K\right) < \infty$ accounts for the fact that for only finitely many $j \in \bbN$, we have that $\min \left\{ \frac{1}{2j + 1}, \frac{1}{K} \right\} = \frac{1}{K}$.
  We conclude, since $\sum_{j = 1}^{\infty} \frac{1}{4j^2} < \infty$,
  \begin{align*}
    \prod_{j=1}^{n-1} \left(1 - \frac{1 + a}{2j + 1 + Ka}\right) &\leq \exp\left(- \sum_{j = 1}^{n} \frac{1}{2j} + C_2\left(K\right)\right) \leq e^{C_2\left(K\right)} \exp\left(-\frac{1}{2} \ln\left(n\right)\right)
  \end{align*}
  We can therefore conclude
  \begin{align*}
    \Prb{Y_0' = 0, \overline{M}_f = n} &\leq \frac{Ka}{2n} \cdot e^{C_2\left(K\right)} \cdot n^{-\frac{1}{2}} = C\left(K\right) an^{-\frac{3}{2}}
  \end{align*}
  This yields, noting that $s \in \left(0, \frac{1}{2}\right)$ and therefore $s-\frac{3}{2} < -1$:
  \begin{align*}
    \Ex{\overline{Q}\left(f\right)^s \indic{Y_0' = 0}} &= \sum_{n \geq 1} n^s \Prb{Y_0' = 0, \overline{M}_f = n} \\
    &\leq \sum_{n \geq 1} C\left(K\right) an^{s-\frac{3}{2}} = C\left(K\right) a \underbrace{\sum_{n \geq 1} n^{s-\frac{3}{2}}}_{< \infty} \leq C\left(s, K\right)a
  \end{align*}
  If $Y_0' = 1$, then $\overline{M}_f = 1$. In addition,
  \begin{align*}
    \Prb{Y_0' = 1, \overline{M}_e = n} &\leq \Prb{Y_n = 1} = \frac{a}{n + 1 + 2a} \leq \frac{a}{n} \\
    \Ex{\overline{Q}\left(f\right)^s \indic{Y_0' = 1}} &= \sum_{n \geq 1} n^{-s} \Prb{Y_0' = 1, \overline{M}_e = n}
    \leq \sum_{n \geq 1} an^{-1-s} < C\left(s\right) a
  \end{align*}
  Thus (recall that the random variables $\overline{Q}$ were constructed to be iid),
  \begin{align*}
    \Exc{v_0}{\prod_{f \in \gamma, f \neq \gamma_1} Q\left(f\right)^{s} \cdot \indic{D_{\gamma}}}
    &\leq \Exc{v_0}{\prod_{f \in \gamma, f \neq \gamma_1} \overline{Q}\left(f\right)^{s} \cdot \indic{D_{\gamma}}}
    \leq \Exc{v_0}{\prod_{f \in \gamma, f \neq \gamma_1} \overline{Q}\left(f\right)^{s}} \\
    &= \prod_{f \in \gamma, f \neq \gamma_1} \Exc{v_0}{\overline{Q}\left(f\right)^{s}} \\
    &= \prod_{f \in \gamma, f \neq \gamma_1} \left(\Ex{\overline{Q}\left(f\right)^s \indic{Y_0' = 0}} + \Ex{\overline{Q}\left(f\right)^s \indic{Y_0' = 1}}\right) \\
    &\leq \left( C\left(s, K\right)a \right)^{\left|\gamma\right| - 1}
  \end{align*}
  This finishes the proof up to the technical details of the coupling of the random variables $Y, Y'$ and the domination of $Q$ by $\overline{Q}$.

  The intuition for how the random variables $Y, Y'$ should be coupled with the LERRW
  was already given in \autoref{fig:bounding_q}: the events $Y = 1$ and $Y' = 1$ upper and lower bound
  the events that a certain edge in the deterministic path $\gamma$ (which was fixed at the beginning of the proof) is crossed.
  We now give the details of the coupling by constructing the LERRW depending on the values
  taken by $Y, Y'$, and some additional randomness. If the random walker already did $t$ steps and is at
  vertex $v$, then the LERRW proceeds as follows ($v \in \gamma$ shall denote that $v$ is in the interior
  of the path $\gamma$, i.e.~not one of its endpoints):
  \begin{itemize}
    \item If $v \notin \gamma$, then the next edge is chosen as always in the LERRW: according to
    the current reinforced edge weights after the first $t$ steps. The choice of edge is independent of the
    variables $Y, Y'$.
    \item If $v \in \gamma$ and the first $t$ steps of the LERRW show that $D_{\gamma}$ does not occur,
    then, as in the previous case, the next edge is chosen according to the edge weights independently of $Y, Y'$.
    We can ignore the values of the variables $Y, Y'$ here because we only claimed $Q \leq \overline{Q}$ if
    $D_{\gamma}$ occurs. The first $t$ steps of the LERRW are inconsistent with $D_{\gamma}$ if
    an edge in $\gamma$ is traversed only after its inverse or if the first arrival to some node in $\gamma$
    was not through the preceding edge in $\gamma$.
    \item If $v \in \gamma$ and the value of $Q$ is already determined by the first $t$ steps of the LERRW,
    then we again ignore $Y, Y'$. By $Q$ being determined we mean that $v = \check{f}$ with $f \in \gamma$ and both
    $f$ as well as $e = \overleftarrow{f'}$ have already been traversed.
  \end{itemize}
  In all remaining cases, we thus have that the first $t$ steps of the LERRW are are consistent with
  $D_{\gamma}$ occurring, and that $v \in \gamma$, i.e.~$v = \check{f}$ with $\gamma_1 \neq f \in \gamma$. Again, call
  $e = \overleftarrow{f'}$. As we are not in the last case, we have that $f$ and $e$ have not both been
  traversed yet. The following cases remain:
  \begin{itemize}
    \item The walker is at $v$ for the first time. All incident edges thus have weight $a$, except for $e$, which has weight
    $1 + a$. The probability to take $e$ is $\frac{1 + a}{1 + \deg\left(v\right) a} \geq \frac{1 + a}{1 + Ka} = \Prb{Y_0' = 1}$.
    The LERRW will be coupled as follows: whenever $Y_0' = 1$, the walk exits along $e$, and sometimes also if $Y_0' = 0$.
    \item Later visits to $v$, $Y_0' = 0$ (the walk did not necessarily exit along $e$ after the first visit). If $v$ is visited
    for the $n$-th time, then the weight of $e$ is still at least $1 + a$, and the total weight of all incident edges is $2n - 1 + \deg\left(v\right)a$.
    The probability to take $e$ is therefore
    \begin{align*}
      \frac{z\left(t, e\right) + a}{2n - 1 + \deg\left(v\right)a} \geq \frac{1 + a}{2n - 1 + Ka} = \Prb{Y_{n - 1}' = 1}
    \end{align*}
    where $z\left(t, e\right)$ was the number of undirected traversals of $e$ up to time $t$. The
    LERRW will now be coupled such that it takes $e$ if $Y_{n - 1}' = 1$, and sometimes also if $Y_{n - 1}' = 0$.
    \item Later visits to $v$, $Y_0' = 1$ (the walk exited along $e$ after the first visit). $f$ has not been traversed yet
    since not both $f$ and $e$ have been traversed already. As the first $t$ steps are consistent with $D_{\gamma}$,
    $\overleftarrow{f}$ has also not been traversed yet, and therefore, the weight of $f$ is $a$. We now decide, independently
    of $Y, Y'$ and with probabilities corresponding to the current edge weights, whether one of the edges $f, e$ will be used or not.
    In the latter case, the variables $Y, Y'$ are ignored again. If $f$ or $e$ is used, and this happens for the $n$-th time,
    then $z\left(t, e\right) \geq n$ since $f$ was not traversed yet, and since $e$ was also traversed in the other
    direction at least once. The probability to choose $f$ from $f, e$ is
    \begin{align*}
      \frac{a}{z\left(t, e\right) + 2a} \leq \frac{a}{n + 2a} = \Prb{Y_{n - 1} = 1}
    \end{align*}
    We couple the LERRW such that if $Y_{n - 1} = 0$, then $e$ is chosen, and sometimes also if $Y_{n - 1} = 1$.
  \end{itemize}
  
  To finish the proof, we give some details on why $Q\left(f\right) \leq \overline{Q}\left(f\right)$ on $D_{\gamma}$ (fix some edge
  $\gamma_1 \neq f \in \gamma$ with the usual notation). Hence we assume that $D_{\gamma}$ occurs.
  \begin{itemize}
    \item If $Y_0' = 0$ (the walk did not necessarily exit along $e$ after the first visit), then only the variables $Y'$ are
    relevant. At every visit to $v$ (until $e$ is used), one variable $Y'$ is considered. By the coupling given above,
    $e$ is used at the latest when $Y_{n - 1}' = 1$ for the first time, during the $n$-th visit to $v$. At this point,
    $f$ can have been taken at most $n - 1$ times. We have
    \begin{align*}
      Q\left(f\right) \leq M_f \leq n - 1 = \overline{M}_f = \overline{Q}\left(f\right)
    \end{align*}
    \item If $Y_0' = 1$ (the walk exited along $e$ after the first visit), then only the variables $Y$ are relevant.
    At every visit to $v$ where $f$ or $e$ is used, one variable $Y$ is considered. $f$ is used at the earliest when $Y_{n - 1} = 1$
    for the first time, during the $n$-th visit where one of $f, e$ is used. At this point, $e$ will have been taken at least
    $n - 1$ times. We have
    \begin{align*}
      Q\left(f\right) &= \frac{1}{M_e} \leq \frac{1}{n - 1} = \frac{1}{\overline{M}_e} = \overline{Q}\left(f\right) \qedhere
    \end{align*}
  \end{itemize}
\end{tproof}

\subsection{Proof of Recurrence}

Before we arrive at the final result, we are now able to finish the last intermediate step:

\begin{theorem}[Exponential Conductance Decay \textnormal{(\cite[Theorem 2]{localizationerrw})}]
  \label{thm:exp_conduc_decay}
  Let $G$ be a graph with degree at most $K$ such that the LERRW on $G$ is recurrent for
  any choice of initial weights, let $v_0 \in V$ be the starting vertex, and let $s \in \left(0, \frac{1}{4}\right)$.
  Then, there is some $a_0 > 0$, depending on $s$ and $K$, such that for initial weights $a \in \left(0, a_0\right)$
  we have
  \begin{align*}
    \Exc{v_0}{c_e^s} \leq 2K \left( C\left(s, K\right)\sqrt{a} \right)^{\dist{v_0}{e}}
  \end{align*}
  where $C\left(s, K\right)$ solely depends on $K$ and $s$.
\end{theorem}

\begin{tproof}
  We have
  \begin{align*}
    \Exc{v_0}{c_e^s} &= \sum_{\gamma \in \Gamma_e} \Exc{v_0}{c_e^s \cdot \indic{D_{\gamma}}}
    \leq \sum_{\gamma \in \Gamma_e} \Exc{v_0}{\left(\frac{c_e}{c_{\gamma_1}}\right)^s \cdot \indic{D_{\gamma}}} \\
    &\overset{\textrm{CSI}}{\leq}
    \sum_{\gamma \in \Gamma_e} \Exc{v_0}{\prod_{f \in \gamma, f \neq \gamma_1} \left(\frac{R\left(f\right)}{Q\left(f\right)}\right)^{2s} \cdot \indic{D_{\gamma}}}^{\frac{1}{2}}
    \Exc{v_0}{\prod_{f \in \gamma, f \neq \gamma_1} Q\left(f\right)^{2s} \cdot \indic{D_{\gamma}}}^{\frac{1}{2}} \\
    &\overset{\circledast}{\leq} \sum_{\gamma \in \Gamma_e}
    \left(C\left(2s\right)^{\left|\gamma\right| - 1}\right)^{\frac{1}{2}}
    \left(\left(C\left(2s,K\right)a\right)^{\left|\gamma\right| - 1}\right)^{\frac{1}{2}} \\
    &= \sum_{\gamma \in \Gamma_e} \left(\vphantom{C\left(2s,K\right)^\frac{1}{2}\sqrt{a}}\right.
    \underbrace{C\left(2s\right)^\frac{1}{2} C\left(2s,K\right)^{\frac{1}{2}}}_{=: C_0}\left. \vphantom{C\left(2s,K\right)^\frac{1}{2}}\sqrt{a} \right)^{\left|\gamma\right| - 1}
  \end{align*}
  where $\circledast$ holds by \autoref{lem:error_bound} and \autoref{lem:estimate_bound}. Now, choose
  $a_0$ such that $KC_0\sqrt{a_0} = \frac{1}{2}$ (\circledchar{1}). Then, for $a < a_0$ (\circledchar{2}), we have, since any $\gamma \in \Gamma_e$
  has length at least $\left|\gamma\right| = \dist{v_0}{e} + 1$  (\circledchar{3}), and since there are at most $K^l$ paths
  of length $l$ (\circledchar{4}):
  \begin{align*}
    \Exc{v_0}{c_e^s} &\leq \sum_{\gamma \in \Gamma_e} \left( C_0 \sqrt{a} \right)^{\left|\gamma\right| - 1}
    \overset{\textrm{\circledchar{3}}\vphantom{X_{\left(k\right)}}}{=} \sum_{l \geq \dist{v_0}{e} + 1} \sum_{\gamma \in \Gamma_e, \left|\gamma\right| = l} \left( C_0 \sqrt{a} \right)^{l-1} \\
    &\overset{\textrm{\circledchar{4}}}{\leq} \sum_{l \geq \dist{v_0}{e} + 1} K \cdot \left( KC_0 \sqrt{a} \right)^{l-1}
    = K \cdot \left(KC_0 \sqrt{a}\right)^{\dist{v_0}{e}} \cdot \sum_{l \geq 0} \left(KC_0 \sqrt{a}\right)^l \\
    &\overset{\textrm{\circledchar{1}\circledchar{2}}}{<} K \cdot \left(KC_0 \sqrt{a}\right)^{\dist{v_0}{e}} \cdot \underbrace{\sum_{l \geq 0} \left(\frac{1}{2}\right)^l}_{2}
    = 2K \left(KC_0 \sqrt{a}\right)^{\dist{v_0}{e}} \qedhere
  \end{align*}
\end{tproof}

We can now finally prove the result we are actually interested in:

\begin{theorem}[Recurrence on Bounded Degree Graphs \textnormal{(\cite[Theorem 1]{localizationerrw})}]
  \label{thm:rec_bndd_deg}
  Let $K \in \bbN_{\geq 1}$. Then, there exists $a_0 > 0$ such that if the graph $G$ has degree at most $K$,
  then the LERRW on $G$ with initial weights set to $a \in \left(0, a_0\right)$ is recurrent.
  The corresponding RWRE is a.s.~positive recurrent, but this does not necessarily imply
  that the expected return time to the starting vertex is finite.
\end{theorem}

\begin{tproof}
  The idea is to apply \autoref{thm:exp_conduc_decay} to the LERRW on finite balls of the infinite
  graph $G$. Call the starting vertex $v_0$ and let $B_R\left(v_0\right) := \left\{ v \in V : \dist{v_0}{v} \leq R \right\}$
  be the set of vertices of distance at most $R$ to the starting vertex. Consider the LERRW $X^{\left(R\right)}$
  on the finite ball $B_R\left(v_0\right)$. Since $B_R\left(v_0\right)$ is finite, $X^{\left(R\right)}$ is recurrent for
  any choice of initial weights, and \autoref{thm:exp_conduc_decay} is applicable. By \autoref{thm:lerrw_is_rwre}, $X^{\left(R\right)}$ is equal in law to a RWRE with random
  conductances $\mathbf{c}^{\left(R\right)} = \left(c_e^{\left(R\right)}\right)_{e \in B_R\left(v_0\right)}$. Denote the mixing measure
  (giving the distribution of $\mathbf{c}^{\left(R\right)}$) by $\mu^{\left(R\right)}$. The measures $\mu^{\left(R\right)}$ are a sequence
  of measures on the set of possible conductance vectors.

  The following calculation will, next to other results, show that for fixed initial weights $a$, the measures
  $\mu^{\left(R\right)}$ are tight (see end of the proof). So, by
  \href{https://en.wikipedia.org/wiki/Prokhorov%27s_theorem}{Prokhorov's theorem}, there is a subsequence
  converging to some measure $\mu$. The first $R$ steps of the LERRW on the whole graph $G$ are
  equal in law to the first $R$ steps of the LERRW on $B_R\left(v_0\right)$. So the LERRW
  on $G$ has the same distribution as a RWRE with mixing measure $\mu$.

  Fix $s \in \left(0, \frac{1}{4}\right)$ and let $e \in B_R\left(v_0\right)$ be an edge in the finite ball of radius
  $R$. By Markov's inequality, and by \autoref{thm:exp_conduc_decay}, for initial weights $a$ small enough, we have
  \begin{align*}
    \mu^{\left(R\right)}\left(c_e^{\left(R\right)} > Q\right)
    = \mu^{\left(R\right)}\left(\left( c_e^{\left(R\right)} \right)^s > Q^s\right)
    \leq \frac{\Ex{\left( c_e^{\left(R\right)} \right)^s}}{Q^s}
    \leq \frac{2K \left( C\left(s, K\right)\sqrt{a} \right)^{\dist{v_0}{e}}}{Q^s}
  \end{align*}
  Choose $Q = \left(2K\right)^{-\dist{v_0}{e}}$ and note that there are at most $K^{l+1}$
  edges at distance $l$ from $v_0$. Then
  \begin{align}
    \label{eq:conduc_bound_q}
    \begin{split}
      &\mu^{\left(R\right)}\left(\exists e \in B_R\left(v_0\right): \dist{v_0}{e} = l \textrm{ and } c_e^{\left(R\right)} > \left(2K\right)^{-l}\right) \\
      &\leq \sum_{e \in B_R\left(v_0\right): \dist{v_0}{e} = l} \mu^{\left(R\right)}\left(c_e^{\left(R\right)} > \left(2K\right)^{-l}\right) \\
      &\leq K^{l+1} 2K \left( C\left(s, K\right)\sqrt{a} \right)^l \left(2K\right)^{sl}
      = 2K^2 \left(2^s K^{1+s} C\left(s, K\right)\sqrt{a}\right)^l
    \end{split}
  \end{align}
  If we now choose $a_0$ such that $2^s K^{1+s} C\left(s, K\right)\sqrt{a_0} \leq \frac{1}{2}$
  and such that $a_0$ is smaller than the bound given by \autoref{thm:exp_conduc_decay}, then, for
  intital weights $a \in \left(0, a_0\right)$,
  \begin{align*}
    \mu^{\left(R\right)}\left(\exists e \in B_R\left(v_0\right): \dist{v_0}{e} = l \textrm{ and } c_e^{\left(R\right)} > \left(2K\right)^{-l}\right)
    < 2K^2 \left(\frac{1}{2}\right)^l = K^2 2^{1 - l}
  \end{align*}
  This bound is uniform in $R$, so it holds for all measures $\mu^{\left(R\right)}$, and in consequence
  also for the limit $\mu$ (by the \href{https://en.wikipedia.org/wiki/Convergence_of_measures#Weak_convergence_of_measures}{Portmanteau theorem}).
  Since $K^2 2^{1 - l}$ is summable in $l$, Borel-Cantelli implies that the probability that infinitely many of the events
  $\left\{c_e > \left(2K\right)^{-\dist{v_0}{e}}\right\}$ occur is $0$.
  In other words, for all but a finite number of edges, $c_e \leq \left(2K\right)^{-\dist{v_0}{e}}$ a.s.
  But then,
  \begin{align*}
    \sum_{e \in E} c_e = \sum_{l \geq 0} \sum_{e \in E: \dist{v_0}{e} = l} c_e
    \leq A + \sum_{l \geq 0} K^{l+1} \left(2K\right)^{-l} = A + K \cdot \sum_{l \geq 0} 2^{-l} < \infty
  \end{align*}
  where $A$ accounts for the finitely many exceptions to the $c_e$ bound. So the RWRE is a.s.~positive recurrent.

  To finish the proof, we will show how tightness follows from the above calculation for fixed initial weights $a$.
  Let $\varepsilon > 0$. We have to show there is a compact set $K_\varepsilon \subseteq \bbR_{\geq 0}^E$ such that for all $R$, 
  the measure of the complement $\mu^{\left(R\right)}\left( K_\varepsilon^{\textrm{C}} \right)$ is smaller than $\varepsilon$.
  If we choose $Q = \left(\lambda K\right)^{-\dist{v_0}{e}}$ in \autoref{eq:conduc_bound_q}, we get
  \begin{align*}
    &\mu^{\left(R\right)}\left(\exists e \in B_R\left(v_0\right): \dist{v_0}{e} = l \textrm{ and } c_e^{\left(R\right)} > \left(\lambda K\right)^{-l}\right) \\
    &\leq 2K^2 \left(\lambda^s K^{1+s} C\left(s, K\right)\sqrt{a}\right)^l
  \end{align*}
  Thus,
  \begin{align*}
    &\mu^{\left(R\right)}\left(\exists e \in B_R\left(v_0\right): c_e^{\left(R\right)} > \left(\lambda K\right)^{-\dist{v_0}{e}}\right) \\
    &\leq \sum_{l \geq 0} \mu^{\left(R\right)}\left(\exists e \in B_R\left(v_0\right): \dist{v_0}{e} = l \textrm{ and } c_e^{\left(R\right)} > \left(\lambda K\right)^{-l}\right) \\
    &\overset{\circledast}{=} \sum_{l \geq 1} \mu^{\left(R\right)}\left(\exists e \in B_R\left(v_0\right): \dist{v_0}{e} = l \textrm{ and } c_e^{\left(R\right)} > \left(\lambda K\right)^{-l}\right) \\
    &\leq 2K^2 \sum_{l \geq 1} \left(\lambda^s K^{1+s} C\left(s, K\right)\sqrt{a}\right)^l
    = \frac{2K^2\lambda^s K^{1+s} C\left(s, K\right)\sqrt{a}}{1 - \lambda^s K^{1+s} C\left(s, K\right)\sqrt{a}}
  \end{align*}
  For $\circledast$, recall that we normalized the conductances to $c_{v_0} = 1$, so for any edge $e$ with distance $0$ to $v_0$, it holds that
  $c_e \leq 1$. We can now choose $\lambda$ such that
  \begin{align*}
    \lambda < \left(\frac{\varepsilon}{K^{1+s}C\left(s, K\right)\sqrt{a}\left(2K^2 + \varepsilon\right)}\right)^{\frac{1}{s}}
  \end{align*}
  Then
  \begin{align*}
    \mu^{\left(R\right)}\left(\exists e \in B_R\left(v_0\right): c_e^{\left(R\right)} > \left(\lambda K\right)^{-\dist{v_0}{e}}\right)
    < \varepsilon
  \end{align*}
  which holds independently of $R$. So we can choose
  \begin{align*}
    K_\varepsilon = \prod_{e \in E} \left[0, \left(\lambda K\right)^{-\dist{v_0}{e}}\right]
  \end{align*}
  which is compact. Note that $K_\varepsilon$ grows when $\varepsilon$ gets smaller, because $\lambda$ must also be chosen
  smaller in this case.
\end{tproof}

\clearpage
\newpage

\section{Reinforced Random Walk with Multiple Walkers}
\label{sec:multiple_walkers}

Here, we consider multiple walkers on a single environment influencing each other.

\subsection{A Two-Player Urn}

We start with a very simple model, the linearly edge-reinforced random walk with two walkers on
a segment of $\bbZ$ with three nodes.
\begin{figure}[H]
  \begin{center}
    \begin{tikzpicture}
      \node[circle,fill=black,inner sep=0.7mm] (N0) at (0, 0) {};
      \node[circle,fill=black,inner sep=0.7mm] (N1) at (2, 0) {};
      \node[circle,fill=black,inner sep=0.7mm] (Nm1) at (-2, 0) {};
      \node[below=5mm] at (N0) {$0$};
      \node[below=5mm] at (Nm1) {$-1$};
      \node[below=5mm] at (N1) {$1$};
      \draw (Nm1) -- node[above=2mm] {$w\left(n,0\right)$} (N0) -- node[above=2mm] {$w\left(n,1\right)$} (N1);
    \end{tikzpicture}
    \caption{Edge weights of the line segment at time $n$}
    \label{fig:errw_segment_multiple}
  \end{center}
\end{figure}
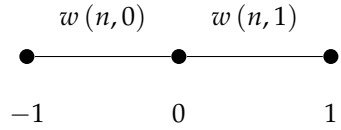

\subsubsection{Alternating Players}

We first define the following dynamics:
\begin{itemize}
  \item There are two walkers $X^{\left(1\right)}$ and $X^{\left(2\right)}$, which both start at the
  node in the center, i.e.~$X^{\left(1\right)}_0 = 0, X^{\left(2\right)}_0 = 0$.
  \item Initially, both edge weights are $1$. We denote the edge weight of the left edge at time $n$
  by $w\left(n,0\right)$, the weight of the right edge by $w\left(n,1\right)$.
  \item Whenever an edge is crossed by either of the walkers, its weight is increased by $1$, so we
  set $W\left(n\right) = 1 + n$ in terms of \autoref{def:errw}.
  \item The walkers move alternately, i.e.~at odd time steps, walker $1$ moves (in particular, walker
  $1$ moves first at step $1$) and at even time steps, walker $2$ moves. This implies that the walkers
  will meet at the node in the center every four steps.
  \item When a walker at the node in the center is about to move, he chooses the edge to traverse with
  probability proportional to the respective edge weight.
\end{itemize}

\begin{lemma}
  \label{lem:two_pl_urn_mart}
  The random variables $\frac{w\left(4n,0\right)}{w\left(4n,0\right) + w\left(4n,1\right)}$ for $n \geq 0$
  form a martingale.
\end{lemma}

\begin{tproof}
  We can calculate:
  \begin{align*}
    &\Ex{\frac{w\left(4n + 4,0\right)}{w\left(4n + 4,0\right) + w\left(4n + 4,1\right)} \midd| w\left(4n, 0\right) = a, w\left(4n, 1\right) = b} \\
    &= \frac{a}{a + b} \frac{a + 1}{a + b + 1} \frac{a + 4}{a + b + 4}
    + \frac{a}{a + b} \frac{b}{a + b + 1} \frac{a + 2}{a + b + 4} \\
    &\hphantom{\;=\;} + \frac{b}{a + b} \frac{a}{a + b + 1} \frac{a + 2}{a + b + 4}
    + \frac{b}{a + b} \frac{b + 1}{a + b + 1} \frac{a}{a + b + 4} \\
    &= \frac{a}{a + b} \cdot \frac{\left(a + 1\right)\left(a + 4\right) + 2b\left(a + 2\right) + b\left(b + 1\right)}{\left(a + b + 1\right)\left(a + b + 4\right)}
    = \frac{a}{a + b} \qedhere
  \end{align*}
\end{tproof}

\begin{corollary}
  \label{cor:two_pl_urn_conv}
  $\frac{w\left(n,0\right)}{w\left(n,0\right) + w\left(n,1\right)}$ converges almost surely for $n \to \infty$
\end{corollary}

\begin{tproof}
  We immediately get from \autoref{lem:two_pl_urn_mart} that $\frac{w\left(4n,0\right)}{w\left(4n,0\right) + w\left(4n,1\right)}$
  converges a.s.~and then, it is easy to see that also $\frac{w\left(n,0\right)}{w\left(n,0\right) + w\left(n,1\right)}$ converges.
\end{tproof}

\begin{conjecture}
  \label{conj:two_pl_urn_limit}
  Define the random variable $Y := \lim_{n \to \infty} \frac{w\left(n,0\right)}{w\left(n,0\right) + w\left(n,1\right)}$.
  Then $Y \in \left[0, 1\right]$ has a density w.r.t.~the Lebesgue measure on $\left[0, 1\right]$.
  $Y$ is not Beta-distributed.
\end{conjecture}

\subsubsection{Random Player Selection}

Consider next the case where at every step, we choose uniformly at random (independent of all other steps)
which of the two walkers moves. \autoref{lem:two_pl_urn_rand_meeting_time} shows that the expected time to meet again in the
middle, if both walkers start in the center, is $4$, just as in the previous case. Of course, the difference
now is that the next meeting time is random. We call $\tau_0 = 0$ and
$\tau_n = \inf\left\{ k > \tau_{n-1} : X^{\left(1\right)}_k = X^{\left(2\right)}_k = 0 \right\}$.

\begin{lemma}
  \label{lem:two_pl_urn_rand_meeting_time}
  For any $n \geq 0$ and any $l \geq 1$, it holds that $\Prb{\tau_{n+1} - \tau_n = 2l} = 2^{-l}$
  and $\Ex{\tau_{n+1} - \tau_n} = 4$.
\end{lemma}

\begin{tproof}
  We consider a MC consisting of three states and coupled with the edge-reinforced random walk.
  The MC is in state $s_{\textrm{center}}$ if both walkers are in the center, in state $s_{\textrm{mixed}}$
  if one walker is in the center and the other in either of the two outer nodes, and in state $s_{\textrm{none}}$
  if none of the walkers is in the center. It is easy to verify that this is indeed a MC
  with the following transition probabilities:
  \begin{center}
    \begin{tikzpicture}
      \node[circle,draw=black,inner sep=0.7mm] (cen) at (-2, 0) {$s_{\textrm{center}}$};
      \node[circle,draw=black,inner sep=0.7mm] (mix) at (0, 0) {$s_{\textrm{mixed}}$};
      \node[circle,draw=black,inner sep=0.7mm] (non) at (2, 0) {$s_{\textrm{none}}$};
      \draw (cen) edge[-{Latex[length=2mm,width=2mm]},bend left] node[above] {$1$} (mix);
      \draw (mix) edge[-{Latex[length=2mm,width=2mm]},bend left] node[below] {$\frac{1}{2}$} (cen);
      \draw (mix) edge[-{Latex[length=2mm,width=2mm]},bend left] node[above] {$\frac{1}{2}$} (non);
      \draw (non) edge[-{Latex[length=2mm,width=2mm]},bend left] node[below] {$1$} (mix);
    \end{tikzpicture}
  \end{center}
  Let $n \geq 0$. At time $\tau_n$, both walkers are in the center and the MC is therefore
  in state $s_{\textrm{center}}$. The time $\tau_{n+1} - \tau_n$ corresponds to the time needed to
  return again to the state $s_{\textrm{center}}$. It is now a standard calculation to show
  $\Prb{\tau_{n+1} - \tau_n = 2l} = 2^{-l}$. Consequently,
  \begin{align*}
    \Ex{\tau_{n+1} - \tau_n} &= \sum_{l \geq 1} 2^{-l} \cdot 2l = \sum_{l \geq 1} l \cdot \left(\frac{1}{2}\right)^{l-1}
    = \frac{1}{\left(1 - \frac{1}{2}\right)^2} = 4 \qedhere
  \end{align*}
\end{tproof}

We next want to prove, in \autoref{lem:two_pl_urn_mart_rand}, that the fraction of the left edge weight,
$\frac{w\left(\tau_n,0\right)}{w\left(\tau_n,0\right) + w\left(\tau_n,1\right)}$,
is a martingale, as it was in the previous case and as in the one-player urn with the difference that
we look at the fraction not at every time step, but at certain stopping times. We use two lemmas for
the proof and introduce the following notation:
\begin{enumerate}[(1)]
  \item \label{term:two_pl_urn_rand_expec} We look at the expectation of the proportion of the left edge weight
  multiplied by the indicator of the event that the walkers need $2l$ steps to meet again:
  \begin{align*}
    \bbE_{a,b,l} := \Ex{\frac{w\left(\tau_{n+1},0\right)}{w\left(\tau_{n+1},0\right) + w\left(\tau_{n+1},1\right)} \cdot \mathbbm{1}_{\left\{\tau_{n+1} - \tau_n = 2l\right\}} \midd| w\left(\tau_n, 0\right) = a, w\left(\tau_n, 1\right) = b}
  \end{align*}
  \item \label{term:two_pl_urn_rand_proba} We also consider the probability that the last walker which
  returns to the center comes from the left node, again intersected with the event that the walkers need $2l$ steps to meet again:
  \begin{align*}
    L_n &:= \left\{ \textrm{the walker returning to the center at time } \tau_{n+1} \textrm{ comes from the left node} \right\} \\
    q_{a,b,l} &:= \Prb{L_n \textrm{ occurs and } \tau_{n+1} - \tau_n = 2l \midd| w\left(\tau_n, 0\right) = a, w\left(\tau_n, 1\right) = b}
  \end{align*}
\end{enumerate}

\begin{lemma}
  \label{lem:two_pl_urn_rand_recursion}
  We have the following recursive equations:
  \begin{align*}
    \bbE_{a,b,l+1} &= \frac{1}{2} \bbE_{a,b,l} + \frac{1}{\left(a+b+2l-1\right)\left(a+b+2l+2\right)} \left(\bbE_{a,b,l} - q_{a,b,l}\right) \\
    q_{a,b,l+1} &= \frac{1}{2} q_{a,b,l} + \frac{1}{4} \cdot \frac{a+b+2l}{a+b+2l-1} \left(\bbE_{a,b,l} - q_{a,b,l}\right)
  \end{align*}
\end{lemma}

\begin{tproof}
  The following notation will be useful. We define a path of the two walkers as a sequence of the symbols
  $1_{\textrm{l}}, 1_{\textrm{r}}, 2_{\textrm{l}}, 2_{\textrm{r}}$ which correspond to the first (respectively
  second) walker moving left and right, where we assume that both walkers start in the center.
  The set $\textrm{Path}_{2l}$ contains all the paths of length $2l$ (a sequence of $2l$ symbols) such
  that the first time at which both walkers are in the center at the same time again is at the end of
  the path. Note that any such path must be of even length since each walker can only be in the center
  after having made an even number of movements.

  For, $\rho \in \textrm{Path}_{2l}$, we set $x\left(\rho\right)$ to be the number of traversals of the
  left edge when the path $\rho$ is taken and we can now write
  \begin{align*}
    \bbE_{a,b,l} &= \sum_{\rho \in \textrm{Path}_{2l}} \underbrace{\Prb{\textrm{the walkers move according to } \rho \midd| w\left(\tau_n, 0\right) = a, w\left(\tau_n, 1\right) = b}}_{=:\bbP_{\rho,a,b}} \cdot \frac{a + x\left(\rho\right)}{a + b + 2l}
  \end{align*}
  To prove the recurrence relation, we build paths of length $2\left(l + 1\right)$ out of paths of length
  $2l$.
  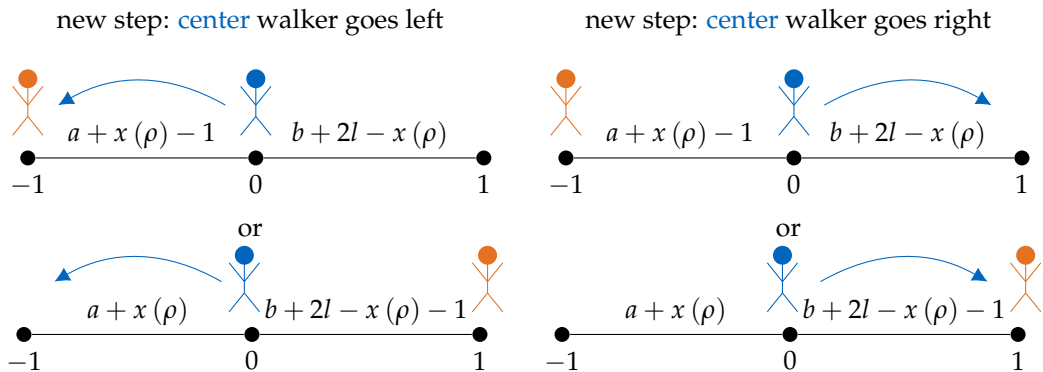
\begin{figure}[H]
    \begin{center}
      \begin{tabular}{cc}
        new step: \textcolor{sectionblue}{center} walker goes left & new step: \textcolor{sectionblue}{center} walker goes right \\
        ~ & ~ \\
        \begin{tikzpicture}
          \node[circle,fill=black,inner sep=0.7mm] (N0) at (0, 0) {};
          \node[circle,fill=black,inner sep=0.7mm] (N1) at (3, 0) {};
          \node[circle,fill=black,inner sep=0.7mm] (Nm1) at (-3, 0) {};
          \stickman{(-3,0.8)}{tumOrange}
          \stickman{(0,0.8)}{sectionblue}
          \node[below=1mm] at (N0) {$0$};
          \node[below=1mm] at (Nm1) {$-1$};
          \node[below=1mm] at (N1) {$1$};
          \draw (Nm1) -- node[above] {$a + x\left(\rho\right) - 1$} (N0) -- node[above] {$b + 2l - x\left(\rho\right)$} (N1);
          \draw[sectionblue] ([shift={(-0.4,0.7)}]N0.center) edge[bend right,-{Latex[length=2mm,width=2mm]}] ([shift={(0.4,0.7)}]Nm1.center);
        \end{tikzpicture}
        & \begin{tikzpicture}
          \node[circle,fill=black,inner sep=0.7mm] (N0) at (0, 0) {};
          \node[circle,fill=black,inner sep=0.7mm] (N1) at (3, 0) {};
          \node[circle,fill=black,inner sep=0.7mm] (Nm1) at (-3, 0) {};
          \stickman{(-3,0.8)}{tumOrange}
          \stickman{(0,0.8)}{sectionblue}
          \node[below=1mm] at (N0) {$0$};
          \node[below=1mm] at (Nm1) {$-1$};
          \node[below=1mm] at (N1) {$1$};
          \draw (Nm1) -- node[above] {$a + x\left(\rho\right) - 1$} (N0) -- node[above] {$b + 2l - x\left(\rho\right)$} (N1);
          \draw[sectionblue] ([shift={(0.4,0.7)}]N0.center) edge[bend left,-{Latex[length=2mm,width=2mm]}] ([shift={(-0.4,0.7)}]N1.center);
        \end{tikzpicture}
        \\ or & or \\
        \begin{tikzpicture}
          \node[circle,fill=black,inner sep=0.7mm] (N0) at (0, 0) {};
          \node[circle,fill=black,inner sep=0.7mm] (N1) at (3, 0) {};
          \node[circle,fill=black,inner sep=0.7mm] (Nm1) at (-3, 0) {};
          \stickman{(3.1,0.8)}{tumOrange}
          \stickman{(-0.1,0.8)}{sectionblue}
          \node[below=1mm] at (N0) {$0$};
          \node[below=1mm] at (Nm1) {$-1$};
          \node[below=1mm] at (N1) {$1$};
          \draw (Nm1) -- node[above] {$a + x\left(\rho\right)$} (N0) -- node[above] {$b + 2l - x\left(\rho\right) - 1$} (N1);
          \draw[sectionblue] ([shift={(-0.4,0.7)}]N0.center) edge[bend right,-{Latex[length=2mm,width=2mm]}] ([shift={(0.4,0.7)}]Nm1.center);
        \end{tikzpicture}
        & \begin{tikzpicture}
          \node[circle,fill=black,inner sep=0.7mm] (N0) at (0, 0) {};
          \node[circle,fill=black,inner sep=0.7mm] (N1) at (3, 0) {};
          \node[circle,fill=black,inner sep=0.7mm] (Nm1) at (-3, 0) {};
          \stickman{(3.1,0.8)}{tumOrange}
          \stickman{(-0.1,0.8)}{sectionblue}
          \node[below=1mm] at (N0) {$0$};
          \node[below=1mm] at (Nm1) {$-1$};
          \node[below=1mm] at (N1) {$1$};
          \draw (Nm1) -- node[above] {$a + x\left(\rho\right)$} (N0) -- node[above] {$b + 2l - x\left(\rho\right) - 1$} (N1);
          \draw[sectionblue] ([shift={(0.4,0.7)}]N0.center) edge[bend left,-{Latex[length=2mm,width=2mm]}] ([shift={(-0.4,0.7)}]N1.center);
        \end{tikzpicture}
      \end{tabular}
      \caption{Possible walker locations after $2l - 1$ steps of a path $\rho$ of length $2l$}
      \label{fig:new_paths}
    \end{center}
  \end{figure}
  For any possible path of the walkers, at any uneven time step, there will be one walker which is in
  the \textcolor{sectionblue}{center} and one which is in one of the \textcolor{tumOrange}{outer} nodes. The paths in $\textrm{Path}_{2l}$ have the additional
  condition that at any even time step (except for the beginning and end), both walkers have to be in
  the outer nodes, not necessarily the same one. We can therefore construct all paths in
  $\textrm{Path}_{2\left(l + 1\right)}$ as follows: for every $\rho \in \textrm{Path}_{2l}$, we take
  the path $\rho$ up to time $2l - 1$. The possible walker locations and edge weights after $2l - 1$
  steps are depicted in \autoref{fig:new_paths}. The path $\rho$ now continues with the walker in the
  \textcolor{tumOrange}{outer} node moving back to the center. To get a path of length $2\left(l + 1\right)$,
  the walker in the center node has to move instead in the next step, and then both walkers will return
  to the center in either order. Out of $\rho$, we can thus construct $4$ new paths of length $2\left(l + 1\right)$:
  we have two choices for which way the center walker moves in step $2l$, and then two choices for
  the order in which the walkers return. Note that we do not construct any path twice, since any two
  different paths in $\textrm{Path}_{2l}$ must already differ somewhere in the first $2l - 1$ steps.

  Now, we look at how this changes the outcome (the final edge weights) and the probability of the path.
  Let $\rho \in \textrm{Path}_{2l}$.
  \begin{itemize}
    \item If the \textcolor{tumOrange}{outer} walker is on the left after $2l - 1$ steps (this is a condition on $\rho$):
    \begin{itemize}
      \item If the \textcolor{sectionblue}{center} walker should move left in step $2l$ (this is one choice for creating the new path):
      the probability $\bbP_{\rho,a,b}$ of $\rho$ is a product over the probability that the walker indicated in $\rho$ is chosen
      at the respective step (which is always $\frac{1}{2}$) and the probability that the walker moves
      in the direction indicated by $\rho$ (this can either be $1$, if the walker moves back to the center,
      or a fraction depending on the edge weights). In our new modified path of length $2\left(l + 1\right)$,
      the first difference is that we choose a different walker to move in step $2l$. This event has
      probability $\frac{1}{2}$, but this is the same as chosing the original walker, hence this part
      is already included in the product giving the probability $\bbP_{\rho,a,b}$ of $\rho$.

      Next, the probability for the center walker to go left, if he is chosen to move at step $2l$,
      is given by $\frac{a + x\left(\rho\right) - 1}{a + b + 2l - 1}$ (compare with \autoref{fig:new_paths}).
      This factor is new and has to be added to the product. In the next step, one of the two walkers
      is chosen and will move back to the center. It is irrelevant which one moves, so we get only a
      new factor of $1$. Finally, in step $2\left(l + 1\right)$, the walker which is still in an
      outer node has to be chosen to move back into the center. This happens with probability $\frac{1}{2}$,
      and this is the final factor to be added to the product. Our new path therefore has probability
      \begin{align*}
        \bbP_{\rho,a,b} \cdot \frac{1}{2} \cdot \frac{a + x\left(\rho\right) - 1}{a + b + 2l - 1}
      \end{align*}
      The new path will end with the following ratio of the left edge weight divided by the total edge weights:
      \begin{align*}
        \frac{a + x\left(\rho\right) + 2}{a + b + 2l + 2}
      \end{align*}
      since the left edge will be traversed twice more by the center walker.
      \item If the \textcolor{sectionblue}{center} walker should move right in step $2l$:
      \begin{align*}
        \textrm{new path probability: }\bbP_{\rho,a,b} \cdot \frac{1}{2} \cdot \frac{b + 2l - x\left(\rho\right)}{a + b + 2l - 1},
        \quad \textrm{new outcome: }\frac{a + x\left(\rho\right)}{a + b + 2l + 2}
      \end{align*}
    \end{itemize}
    \item If the \textcolor{tumOrange}{outer} walker is on the right after $2l - 1$ steps:
    \begin{itemize}
      \item If the \textcolor{sectionblue}{center} walker should move left in step $2l$:
      \begin{align*}
        \textrm{new path probability: }\bbP_{\rho,a,b} \cdot \frac{1}{2} \cdot \frac{a + x\left(\rho\right)}{a + b + 2l - 1},
        \quad \textrm{new outcome: }\frac{a + x\left(\rho\right) + 2}{a + b + 2l + 2}
      \end{align*}
      \item If the \textcolor{sectionblue}{center} walker should move right in step $2l$:
      \begin{align*}
        \textrm{new path probability: }\bbP_{\rho,a,b} \cdot \frac{1}{2} \cdot \frac{b + 2l - x\left(\rho\right) - 1}{a + b + 2l - 1},
        \quad \textrm{new outcome: }\frac{a + x\left(\rho\right)}{a + b + 2l + 2}
      \end{align*}
    \end{itemize}
  \end{itemize}

  For $\rho \in \textrm{Path}_{2l}$, we write $\textrm{left}\left(\rho\right) = 1$ if the outer walker
  after step $2l - 1$ is on the left, and $\textrm{left}\left(\rho\right) = 0$ otherwise. Using the
  recursive path construction above, we get that
  \begin{align*}
    &\bbE_{a,b,l+1} = \sum_{\rho \in \textrm{Path}_{2\left(l+1\right)}} \bbP_{\rho,a,b} \cdot \frac{a + x\left(\rho\right)}{a + b + 2l + 2} =\\
    &\sum_{\rho \in \textrm{Path}_{2l}} \bbP_{\rho,a,b} \cdot \frac{1}{2} \cdot
    \left(
      \textrm{left}\left(\rho\right) \cdot \left(
        \frac{a + x\left(\rho\right) - 1}{a + b + 2l - 1} \cdot \frac{a + x\left(\rho\right) + 2}{a + b + 2l + 2}
        + \frac{b + 2l - x\left(\rho\right)}{a + b + 2l - 1} \cdot \frac{a + x\left(\rho\right)}{a + b + 2l + 2}
      \right) \right. \\
    & \qquad\;\left. +
      \left(1 - \textrm{left}\left(\rho\right)\right) \cdot \left(
        \frac{a + x\left(\rho\right)}{a + b + 2l - 1} \cdot \frac{a + x\left(\rho\right) + 2}{a + b + 2l + 2}
        + \frac{b + 2l - x\left(\rho\right) - 1}{a + b + 2l - 1} \cdot \frac{a + x\left(\rho\right)}{a + b + 2l + 2}
      \right)
    \right)
  \end{align*}
  Using this expression for $\bbE_{a,b,l+1}$, we can calculate (a simple, but longer calculation which
  we skip here, see
  \href{https://www.wolframalpha.com/input?i=simplify+%28a%2Bx-1%29%28a%2Bx%2B2%29%28a%2Bb%2B2l%29+%2B+%28b%2B2l-x%29%28a%2Bx%29%28a%2Bb%2B2l%29+-+%28a%2Bx%29%28a%2Bb%2B2l-1%29%28a%2Bb%2B2l%2B2%29}{https://bit.ly/3trqm2r}
  and
  \href{https://www.wolframalpha.com/input?i=simplify+%28a%2Bx%29%28a%2Bx%2B2%29%28a%2Bb%2B2l%29+%2B+%28b%2B2l-x-1%29%28a%2Bx%29%28a%2Bb%2B2l%29+-+%28a%2Bx%29%28a%2Bb%2B2l-1%29%28a%2Bb%2B2l%2B2%29}{https://bit.ly/391ma2o}):
  \begin{align*}
    \bbE_{a,b,l+1} - \frac{1}{2} \bbE_{a,b,l} &=
    \sum_{\rho \in \textrm{Path}_{2l}} \bbP_{\rho,a,b} \cdot
    \frac{1}{\left(a + b + 2l - 1\right)\left(a + b + 2l + 2\right)} \;\cdot \\
    &\hphantom{=\qquad} 
    \left(
      \textrm{left}\left(\rho\right) \cdot \frac{x\left(\rho\right) - b - 2l}{a + b + 2l}
      +  \left(1 - \textrm{left}\left(\rho\right)\right) \cdot \frac{a + x\left(\rho\right)}{a + b + 2l}
    \right) \\
    &= \frac{1}{\left(a + b + 2l - 1\right)\left(a + b + 2l + 2\right)} \;\cdot \\
    &\hphantom{=\qquad} \sum_{\rho \in \textrm{Path}_{2l}} \bbP_{\rho,a,b} \cdot
    \left(
      \textrm{left}\left(\rho\right) \cdot \frac{- a - b - 2l}{a + b + 2l}
      +  \frac{a + x\left(\rho\right)}{a + b + 2l}
    \right) \\
    &= \frac{1}{\left(a + b + 2l - 1\right)\left(a + b + 2l + 2\right)} \cdot \left(
      \bbE_{a,b,l} - \sum_{\rho \in \textrm{Path}_{2l}} \bbP_{\rho,a,b} \cdot \textrm{left}\left(\rho\right)
    \right) \\
    &= \frac{1}{\left(a + b + 2l - 1\right)\left(a + b + 2l + 2\right)} \cdot \left(
      \bbE_{a,b,l} - q_{a,b,l}
    \right)
  \end{align*}
  This proves the first equation. For the second equation, we use the same strategy. For our newly
  constructed paths, we have to calculate the probability that the last walker to move to the center
  comes from the left node. Let $\rho \in \textrm{Path}_{2l}$.
  \begin{itemize}
    \item If the \textcolor{tumOrange}{outer} walker is on the left after $2l - 1$ steps (this is a condition on $\rho$):
    \begin{itemize}
      \item If the \textcolor{sectionblue}{center} walker should move left in step $2l$ (this is one choice for creating the new path):
      \begin{align*}
        \textrm{new path probability: }\bbP_{\rho,a,b} \cdot \frac{1}{2} \cdot \frac{a + x\left(\rho\right) - 1}{a + b + 2l - 1},
        \quad \textrm{probability of }L_n\textrm{: }1
      \end{align*}
      \item If the \textcolor{sectionblue}{center} walker should move right in step $2l$:
      \begin{align*}
        \textrm{new path probability: }\bbP_{\rho,a,b} \cdot \frac{1}{2} \cdot \frac{b + 2l - x\left(\rho\right)}{a + b + 2l - 1},
        \quad \textrm{probability of }L_n\textrm{: }\frac{1}{2}
      \end{align*}
    \end{itemize}
    \item If the \textcolor{tumOrange}{outer} walker is on the right after $2l - 1$ steps:
    \begin{itemize}
      \item If the \textcolor{sectionblue}{center} walker should move left in step $2l$:
      \begin{align*}
        \textrm{new path probability: }\bbP_{\rho,a,b} \cdot \frac{1}{2} \cdot \frac{a + x\left(\rho\right)}{a + b + 2l - 1},
        \quad \textrm{probability of }L_n\textrm{: }\frac{1}{2}
      \end{align*}
      \item If the \textcolor{sectionblue}{center} walker should move right in step $2l$:
      \begin{align*}
        \textrm{new path probability: }\bbP_{\rho,a,b} \cdot \frac{1}{2} \cdot \frac{b + 2l - x\left(\rho\right) - 1}{a + b + 2l - 1},
        \quad \textrm{probability of }L_n\textrm{: }0
      \end{align*}
    \end{itemize}
  \end{itemize}
  Therefore
  \begin{align*}
    q_{a,b,l+1} - \frac{1}{2} q_{a,b,l} &=
    \sum_{\rho \in \textrm{Path}_{2l}} \bbP_{\rho,a,b} \cdot \frac{1}{2} \;\cdot \\
    &\hphantom{=\qquad}\left(
      \frac{a + x\left(\rho\right) - 1}{a + b + 2l - 1} \cdot \textrm{left}\left(\rho\right)
      + \frac{1}{2} \cdot \frac{b + 2l - x\left(\rho\right)}{a + b + 2l - 1} \cdot \textrm{left}\left(\rho\right) \right.\\
    &\hphantom{=\qquad(}\left.  + \;\frac{1}{2} \cdot \frac{a + x\left(\rho\right)}{a + b + 2l - 1} \cdot \left(1 - \textrm{left}\left(\rho\right)\right)
      - \textrm{left}\left(\rho\right)
    \right)\\
    &= \frac{1}{2} \cdot \sum_{\rho \in \textrm{Path}_{2l}} \bbP_{\rho,a,b} \;\cdot \\
    &\hphantom{=\qquad}\left(
      -\frac{1}{2} \cdot \frac{a + b + 2l}{a + b + 2l - 1} \cdot \textrm{left}\left(\rho\right)
      + \frac{1}{2} \cdot \frac{a + x\left(\rho\right)}{a + b + 2l - 1}
    \right) \\
    &= \frac{1}{4} \cdot \frac{a + b + 2l}{a + b + 2l - 1} \cdot \sum_{\rho \in \textrm{Path}_{2l}} \bbP_{\rho,a,b} \cdot
    \left(
      \frac{a + x\left(\rho\right)}{a + b + 2l} - \textrm{left}\left(\rho\right)
    \right) \\
    &= \frac{1}{4} \cdot \frac{a + b + 2l}{a + b + 2l - 1} \cdot \left(\bbE_{a,b,l} - q_{a,b,l}\right) \qedhere
  \end{align*}
\end{tproof}

\begin{corollary}
  \label{cor:two_pl_urn_rand_recursion_solved}
  The expectation of the proportion of the left edge weight and the probability that the last walker
  returning to the center comes from the left coincide, and:
  \begin{align*}
    \bbE_{a,b,l} = q_{a,b,l} = \frac{1}{2^l} \cdot \frac{a}{a + b}
  \end{align*}
\end{corollary}

\begin{tproof}
  We use \autoref{lem:two_pl_urn_rand_recursion}. Let us first calculate $\bbE_{a,b,1}$ and $q_{a,b,1}$.
  There are four possible paths of length $2$ which end again with both walkers in the center: first,
  we choose which of the two walker moves, and this walker can than either move left or right and then
  back to the center. Since it is irrelevant which walker we choose in the beginning, we can disregard
  which walker moves. We thus get:
  \begin{itemize}
    \item With probability $\frac{a}{a + b} \cdot \frac{1}{2}$ the walker which moves in the first step
    moves left and is then chosen again to move back to the center in the next step. In this case, $L_n$
    occurs and the resulting edge weight ratio is $\frac{a + 2}{a + b + 2}$.
    \item With probability $\frac{b}{a + b} \cdot \frac{1}{2}$ the walker which moves in the first step
    moves right and is then chosen again to move back to the center in the next step. In this case, $L_n$
    does not occur and the resulting edge weight ratio is $\frac{a}{a + b + 2}$.
  \end{itemize}
  We see directly that $q_{a,b,1} = \frac{1}{2} \cdot \frac{a}{a + b}$ and that
  \begin{align*}
    \bbE_{a,b,1} = \frac{1}{2} \cdot \left(\frac{a}{a + b} \cdot \frac{a + 2}{a + b + 2} + \frac{b}{a + b} \cdot \frac{a}{a + b + 2}\right)
    = \frac{1}{2} \cdot \frac{a}{a + b}
  \end{align*}
  The remaining proof is now a simple induction using \autoref{lem:two_pl_urn_rand_recursion}, where
  it should be noted that $\bbE_{a,b,l} - q_{a,b,l} = 0$ under the induction assumption.
\end{tproof}

\begin{lemma}
  \label{lem:two_pl_urn_mart_rand}
  The random variables $\frac{w\left(\tau_n,0\right)}{w\left(\tau_n,0\right) + w\left(\tau_n,1\right)}$ for $n \geq 0$
  form a martingale.
\end{lemma}

\begin{tproof}
  We have by \autoref{cor:two_pl_urn_rand_recursion_solved}:
  \begin{align*}
    &\Ex{\frac{w\left(\tau_{n+1},0\right)}{w\left(\tau_{n+1},0\right) + w\left(\tau_{n+1},1\right)} \cdot \mathbbm{1}_{\left\{\tau_{n+1} - \tau_n = 2l\right\}} \midd| w\left(\tau_n, 0\right) = a, w\left(\tau_n, 1\right) = b}
    = 2^{-l} \cdot \frac{a}{a + b} \\
    \implies
    &\Ex{\frac{w\left(\tau_{n+1},0\right)}{w\left(\tau_{n+1},0\right) + w\left(\tau_{n+1},1\right)} \midd| w\left(\tau_n, 0\right) = a, w\left(\tau_n, 1\right) = b}
    = \frac{a}{a+b} \cdot \sum_{l \geq 1} 2^{-l} = \frac{a}{a+b} \qedhere
  \end{align*}
\end{tproof}

\begin{lemma}
  \label{lem:two_pl_urn_vals_rand}
  Assume that $w\left(0,0\right) = a > 0, w\left(0,1\right) = b > 0$.
  Then, the random variables $M_n := \frac{w\left(\tau_n,0\right)}{w\left(\tau_n,0\right) + w\left(\tau_n,1\right)}$
  take values in the set
  \begin{align*}
    R := \left\{ \frac{a + 2x}{a + b + 2l} \midd| x \in \bbN_{\geq 0}, l \in \bbN_{\geq 1}, x \leq l \right\}
  \end{align*}
  and, for every $r \in R$ and every $n \in \bbN_{\geq 1}$, we have that $\Prb{M_n = r} > 0$.
  The set $R$ is dense in $\left[0, 1\right]$.
\end{lemma}

\begin{tproof}
  It is immediately clear that $M_n$ can take only values in the set $R$: $w\left(\tau_n,0\right)$
  will be equal to $a$ plus the number of crossings of the left edge until time $\tau_n$. Since both
  walkers are in the center at time $\tau_n$, this number of edge crossings has to be even, hence
  we can write $w\left(\tau_n,0\right) = a + 2x$ for some $x \in \bbN_{\geq 0}$. On the other hand,
  the total edge weight also must have increased by an even number since both walkers have to do an
  even number of steps before meeting again in the center.

  We now want to construct a sequence of walker steps leading to any possible outcome $\frac{a + 2x}{a + b + 2l}$
  with positive probability. First, note that $w\left(\tau_n,0\right) + w\left(\tau_n,1\right) = a + b + \tau_n$.
  Since $\tau_n \geq 2n$, we cannot have a sequence of length $2l$ leading to our desired outcome when
  $l < n$. In this case, we have to reach $\frac{a + 2x}{a + b + 2l}$ by a longer sequence. We have, for any $k$,
  \begin{align*}
    \frac{a + 2x}{a + b + 2l} = \frac{k\left(a + 2x\right)}{k\left(a + b + 2l\right)}
    = \frac{a + 2\left(kx + \frac{1}{2}\left(k-1\right)a\right)}{a + b + 2\left(kl + \frac{1}{2}\left(k-1\right)\left(a+b\right)\right)}
  \end{align*}
  To reach our desired fraction, we thus want a sequence of length $2\left(kl + \frac{1}{2}\left(k-1\right)\left(a+b\right)\right) \geq 2n$.
  In addition, $\left(k-1\right)\left(a+b\right)$ and $\left(k-1\right)a$
  should both be even. It thus suffices to choose an uneven $k$ which is large enough such that
  $kl + \frac{1}{2}\left(k-1\right)\left(a+b\right) \geq n$. This allows us to assume
  w.l.o.g.~that $l \geq n$.

  What remains to show is that we can construct a sequence of length $2l$ such that after all $2l$
  steps, the walkers meet again in the center for the $n$-the time and such that the left edge is
  traversed $2x$ times. Since any finite sequence of possible walker movements has positive probability,
  this concludes the proof. Finding such a sequence is easy: in the first $2\left(n-1\right)$ steps,
  one walker moves $\min\left\{x, n-1\right\}$ times to the left node and back, and in the remaining
  steps (if there are any), the same walker moves to the right node and back. Then, in the $\left(2n-1\right)$-th
  step, the same walker moves to the left node, if the left edge was not yet crossed $2x$ times, and otherwise
  to the right node. The walker then remains there until returning to the center in the $2l$-th step.
  In the remaining steps of the complete sequence of length $2l$, the other walker moves in such a way
  as to reach a total number of $2x$ left edge crossings. This leads to the desired fraction at time
  $\tau_n$.

  It is easy to see that $R$ is dense in $\left[0, 1\right]$. Indeed, choosing $x = k \cdot p$ and
  $l = k \cdot q$ for increasing $k$ shows that we can find sequences in $R$ converging to any rational number,
  and these are already dense in $\left[0, 1\right]$.
\end{tproof}

\begin{corollary}
  \label{cor:two_pl_urn_conv_rand}
  $\frac{w\left(n,0\right)}{w\left(n,0\right) + w\left(n,1\right)}$ converges almost surely for $n \to \infty$,
  and the limit is identical to the limit of $\frac{w\left(\tau_n,0\right)}{w\left(\tau_n,0\right) + w\left(\tau_n,1\right)}$.
\end{corollary}

\begin{tproof}
  We use the following notation:
  \begin{align*}
    F_n := \frac{w\left(n,0\right)}{w\left(n,0\right) + w\left(n,1\right)} \qquad
    M_n := F_{\tau_n} = \frac{w\left(\tau_n,0\right)}{w\left(\tau_n,0\right) + w\left(\tau_n,1\right)} \qquad
    M_{\infty} := \lim\limits_{n \to \infty} M_n
  \end{align*}
  Let $\varepsilon > 0$ and consider the events
  \begin{align*}
    A_n := \left\{ \left|F_k - M_{\infty}\right| > \varepsilon \textrm{ for some } k \in \left[\tau_n, \tau_{n+1}\right] \right\}
  \end{align*}
  It is sufficient to show that only finitely many of the events $A_n$ can occur a.s.
  Further set
  \begin{align*}
    B_n := \left\{ \left|F_k - M_n\right| > \frac{\varepsilon}{2} \textrm{ for some } k \in \left[\tau_n, \tau_{n+1}\right] \right\}
  \end{align*}
  Since $M_n$ converges a.s.~by \autoref{lem:two_pl_urn_mart_rand}, there is some (random) $N \in \bbN$ such that for all $n \geq N$,
  it holds that $\left|M_n - M_{\infty}\right| < \frac{\varepsilon}{2}$. For $n \geq N$, the occurrence of $A_n$
  implies that $B_n$ occurs as well, so it is sufficient to show that only finitely many of the events
  $B_n$ can occur.
  
  Now, at time $\tau_n$, the random walkers must have moved at least $2n$ times, so
  $w\left(\tau_n,0\right) + w\left(\tau_n,1\right) \geq 2n$. If $\left|F_k - M_n\right| = \left|F_k - F_{\tau_n}\right| > \frac{\varepsilon}{2}$
  for some $k \in \left[\tau_n, \tau_{n+1}\right]$, it is therefore necessary that at least $\varepsilon n$
  steps were made by the walkers between time $\tau_n$ and time $k$, since every step changes the value
  of $F_k$ by at most $\frac{1}{2n}$. Thus
  \begin{align*}
    B_n \subseteq \left\{ \tau_{n+1} - \tau_n \geq \varepsilon n \right\}
  \end{align*}
  We have
  \begin{align*}
    \Prb{\tau_{n+1} - \tau_n \geq \varepsilon n}
    &= \sum_{l \geq \left\lceil \frac{\varepsilon n}{2} \right\rceil} \Prb{\tau_{n+1} - \tau_n = 2l}
    \overset{\textrm{\autoref{lem:two_pl_urn_rand_meeting_time}}}{=} \sum_{l \geq \left\lceil \frac{\varepsilon n}{2} \right\rceil} 2^{-l} \\
    &\leq 2^{-\frac{\varepsilon n}{2}} \sum_{l \geq 0} 2^{-l} = 2^{1 - \frac{\varepsilon n}{2}} \\
    \implies \sum_{n \geq 1} \Prb{\tau_{n+1} - \tau_n \geq \varepsilon n}
    &\leq \sum_{n \geq 1} 2^{1 - \frac{\varepsilon n}{2}} = 2 \cdot \sum_{n \geq 1} \left(2^{-\frac{\varepsilon}{2}}\right)^n < \infty
  \end{align*}
  By Borel-Cantelli, it follows that only finitely many of the events $\left\{ \tau_{n+1} - \tau_n \geq \varepsilon n \right\}$
  can occur, and therefore also only finitely many of the events $B_n$. This concludes the proof.
\end{tproof}

The same model with only one random walker results in a standard P\'olya urn where two balls of the
drawn color are added after every draw. For the P\'olya urn, the fraction of balls of one of the colors
is a martingale and converges to a limit which is Beta-distributed. We have shown now for the two-player
urn that the fraction of the left edge weight also has a limit, even if it is a martingale only if looked
at at certain stopping times. Simulations and intuition suggest that the limit should also have a density.
However, we couldn't prove this so far, so we only give the following conjecture. Simulations also suggest
that the limiting distribution is not a Beta distribution.

\begin{conjecture}
  \label{conj:two_pl_urn_limit_rand}
  Define the random variable $Y := \lim_{n \to \infty} \frac{w\left(n,0\right)}{w\left(n,0\right) + w\left(n,1\right)}$.
  Then $Y \in \left[0, 1\right]$ has a density w.r.t.~the Lebesgue measure on $\left[0, 1\right]$.
  $Y$ is not Beta-distributed, and its distribution is different from the distribution of the limit
  in \autoref{conj:two_pl_urn_limit}.
\end{conjecture}

\subsection{Model on $\bbZ$}

So far, we have looked at the LERRW with multiple walkers only on a very simple graph and
only with $2$ walkers. We now look at the edge-reinforced random walk on $\bbZ$ with an arbitrary finite number of
random walkers. Formally, we have $k$ sequences (for $k$ walkers)
$\left(X^{\left(m\right)}_n\right)_{n \geq 0}, 1 \leq m \leq k$ of random variables with the following
dynamics. The transition probabilities depend on the edge weights $w\left(n, j\right) > 0$ for $n \geq 0, j \in \bbZ$
where $j$ corresponds to the edge from $j$ to $j+1$.
\begin{figure}[H]
  \begin{center}
    \begin{tikzpicture}
      \node[circle,fill=black,inner sep=0.7mm] (N0) at (0, 0) {};
      \node[circle,fill=black,inner sep=0.7mm] (N1) at (2, 0) {};
      \node[circle,fill=black,inner sep=0.7mm] (N2) at (4, 0) {};
      \node[circle,fill=black,inner sep=0.7mm] (N3) at (6, 0) {};
      \node[circle,fill=black,inner sep=0.7mm] (Nm1) at (-2, 0) {};
      \node[circle,fill=black,inner sep=0.7mm] (Nm2) at (-4, 0) {};
      \node[below=5mm] at (N0) {$0$};
      \node[below=5mm] at (N1) {$j-1$};
      \node[below=5mm] at (N2) {$j$};
      \node[below=5mm] at (N3) {$j+1$};
      \draw (Nm2) -- (Nm1) -- (N0) -- (N1) -- node[above=2mm] {$w\left(n,j-1\right)$} (N2) -- node[above=2mm] {$w\left(n,j\right)$} (N3);
      \draw[dashed] (N3) -- (7, 0);
      \draw[dashed] (Nm2) -- (-5, 0);
      \node (aN2l) at (3.9, 0.5) {};
      \node (aN2r) at (4.1, 0.5) {};
    \end{tikzpicture}
    \caption{Edge weights on $\bbZ$ at time $n$}
    \label{fig:errw_z_multiple}
  \end{center}
\end{figure}
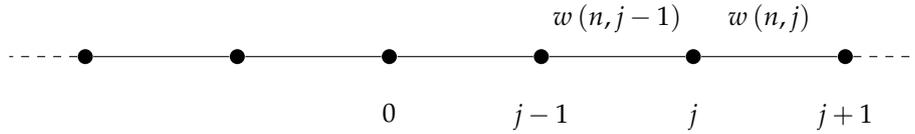
If $\calG_n$ denotes $\sigma\left(\left\{ X^{\left(m_2\right)}_{m_1} : 0 \leq m_1 \leq n, 1 \leq m_2 \leq k \right\} \cup \left\{w\left(m, j\right) : 0 \leq m \leq n, j \in \bbZ \right\}\right)$
i.e.~the history of the random walkers and edge weights up to and including time $n$, then we define,
conditioned on $\calG_n$, the following transition probabilities:
\begin{itemize}
  \item A random walker $m$ ($1 \leq m \leq k$) which is going to jump is chosen uniformly at random
  (independent of $\calG_n$) amongst the $k$ walkers.
  \item If the chosen random walker is at position $j$ (i.e.~$X^{\left(m\right)}_n = j$), he then jumps
  \begin{itemize}
    \item to the right (i.e.~$X^{\left(m\right)}_{n + 1} = j + 1$) with probability $\frac{w\left(n, j\right)}{w\left(n, j - 1\right) + w\left(n, j\right)}$
    \item to the left (i.e.~$X^{\left(m\right)}_{n + 1} = j - 1$) with probability $\frac{w\left(n, j - 1\right)}{w\left(n, j - 1\right) + w\left(n, j\right)}$
  \end{itemize}
  i.e.~the jump probabilities are proportional to the corresponding edge weights.
  \item If $j^{\ast}$ is the traversed edge ($j^{\ast} = j$ if the walker jumps to the right, $j^{\ast} = j - 1$
  if he jumps to the left), then for $i \neq j^\ast$, $w\left(n, i\right) = w\left(n + 1, i\right)$
  and $w\left(n, j^{\ast}\right) \leq w\left(n + 1, j^{\ast}\right)$, i.e.~the weight of the
  traversed edge may be increased according to some reinforcement scheme.
  
  We consider schemes where
  the increment $w\left(n + 1, j^{\ast}\right) - w\left(n, j^{\ast}\right)$ may solely depend on
  $j^{\ast}$, and the number of times the edge was crossed up to time $n$. In other words, $w\left(n, j\right)$
  can still be written in terms of the weight function $\ensuremath{W_e\left(k\right)}$ as defined in \autoref{def:errw}.
  Some generalizations are possible, such as making the weight also depend on $n$ (i.e.~at which times
  the edge was crossed), but will not be considered here.
  \item The initial edge weights can be chosen arbitrarily, but all of them must be positive.
  \item The initial positions of the $k$ walkers can be chosen arbitrarily.
\end{itemize}

\subsection{Recurrence or Finite Range on $\bbZ$}

\begin{definition}
  \label{def:rectrans_multiple}
  We say that one of the walkers $m$ ($1 \leq m \leq k$)
  \begin{itemize}
    \item is \textbf{transient}, if he visits every integer only finitely often, that is, every integer
    appears only finitely often in the sequence $\left(X^{\left(m\right)}_n\right)_{n \geq 1}$
    \item is \textbf{recurrent}, if he visits every integer infinitely often, that is, every integer
    appears infinitely often in the sequence $\left(X^{\left(m\right)}_n\right)_{n \geq 1}$
    \item has \textbf{finite range}, if he only visits finitely many integers, that is, the number of distinct
    integers appearing in the sequence $\left(X^{\left(m\right)}_n\right)_{n \geq 1}$ is finite \qedhere
  \end{itemize}
\end{definition}

\begin{lemma}
  \label{lem:retRootOrFinManyNew}
  Assume the edge-reinforced random walk with $k$ walkers starts with an initial configuration of the
  weights $w\left(0, j\right)$ such that all but finitely many of them are $1$.
  Then, for every random walker $m$ ($1 \leq m \leq k$) with $X^{\left(m\right)}_0 \geq 0$, we have the following:
  \begin{align*}
    &\Prb{X^{\left(m\right)}_n = 0 \textrm{ for some }n \geq 0} \\
    +\; &\bbP\left[ X^{\left(m\right)}_n \neq 0\textrm{ for all }n \geq 0 \textrm{ and } X^{\left(m\right)} \textrm{ only visits finitely many nodes}\right. \\
    &\hphantom{\bbP[}\left. \vphantom{X^{\left(m\right)}_n}\qquad\quad \textrm{which have not been visited before by any other walker} \right] = 1 \qedhere
  \end{align*}
\end{lemma}

\begin{tproof}
  We follow the proof of \cite[Lemma 3.0]{errwz}. Consider a fixed random walker $m$ with $X^{\left(m\right)}_0 \geq 0$. We now define:
  \begin{align*}
    F\left(n, j\right) &:= \begin{cases}
      \sum_{i=0}^{j-1} \frac{1}{w\left(n, i\right)} & \textrm{ if } j > 0 \\
      0 & \textrm{ if } j \leq 0
    \end{cases}
    \qquad \qquad \tau^{\left(m\right)} := \inf \left\{ n \geq 0: X^{\left(m\right)}_n \leq 0 \right\} \\
    M^{\left(m\right)}_n &:= F\left(n \land \tau^{\left(m\right)}, X^{\left(m\right)}_{n \land \tau^{\left(m\right)}}\right) \\
    H^{\left(m\right)}_n &:= M^{\left(m\right)}_n + \sum_{i=1}^n
    \left(\frac{1}{w\left(i-1,X^{\left(m\right)}_{i-1}\right)} - \frac{1}{w\left(i,X^{\left(m\right)}_{i-1}\right)}\right)
    \cdot \mathbbm{1}_{X^{\left(m\right)}_i > X^{\left(m\right)}_{i-1}, i\leq\tau^{\left(m\right)}} \\
    &\hphantom{\;:=\;} + \sum_{i = 1}^n \; \sum_{l = 1, l \neq m}^k \; \sum_{j = 0}^{\infty}
    \left(\frac{1}{w\left(i-1,j\right)} - \frac{1}{w\left(i,j\right)}\right)
    \cdot \underbrace{\mathbbm{1}_{\left\{X^{\left(l\right)}_{i-1}, X^{\left(l\right)}_i\right\} = \left\{j, j + 1\right\}, i\leq\tau^{\left(m\right)}, j < X^{\left(m\right)}_i}}_{= 1 \textrm{ for at most one pair of }l,j}
  \end{align*}
  $M^{\left(m\right)}_n$ is nonnegative by definition of $F$, and $H^{\left(m\right)}_n \geq M^{\left(m\right)}_n \geq 0$
  since edge weights can only increase and therefore, all terms in the sums in the definition of $H^{\left(m\right)}_n$
  are nonnegative. $H^{\left(m\right)}_n$ is a martingale: let
  \begin{align*}
    d^{\left(m\right)}_n \; &:= \\
    H^{\left(m\right)}_n - H^{\left(m\right)}_{n-1}
    &= \underbrace{M^{\left(m\right)}_n - M^{\left(m\right)}_{n-1}}_{:= e^{\left(m\right)}_n}
    + \underbrace{\left(\frac{1}{w\left(n-1,X^{\left(m\right)}_{n-1}\right)} - \frac{1}{w\left(n,X^{\left(m\right)}_{n-1}\right)}\right)
    \cdot \mathbbm{1}_{X^{\left(m\right)}_n > X^{\left(m\right)}_{n-1}, n\leq\tau^{\left(m\right)}}}_{:= f^{\left(m\right)}_n} \\
    &\hphantom{\;=\;}+ \underbrace{\sum_{l = 1, l \neq m}^k \; \sum_{j = 0}^{\infty}
    \left(\frac{1}{w\left(n-1,j\right)} - \frac{1}{w\left(n,j\right)}\right)
    \cdot \mathbbm{1}_{\left\{X^{\left(l\right)}_{n-1}, X^{\left(l\right)}_n\right\} = \left\{j, j + 1\right\}, n\leq\tau^{\left(m\right)}, j < X^{\left(m\right)}_n}}_{:= g^{\left(m\right)}_n}
  \end{align*}
  then we have to show that $\Ex{d^{\left(m\right)}_n \midd| \calG_{n-1}} = 0$. We have:
  \begin{itemize}
    \item if $n - 1 \geq \tau^{\left(m\right)}$, then $d^{\left(m\right)}_n = 0$. Hence, it suffices to consider the case $X^{\left(m\right)}_{n-1} = j > 0$
    and $\tau^{\left(m\right)} \geq n$.
    \item with probability $\frac{1}{k}$, the walker $m$ jumps at time $n-1$. In this case, $g^{\left(m\right)}_n = 0$
    since no other walker can jump and the indicator variable in $g^{\left(m\right)}_n$ is therefore $0$.
    If he jumps to the right (with probability
    $\frac{1}{k} \cdot \frac{w\left(n - 1, j\right)}{w\left(n - 1, j - 1\right) + w\left(n - 1, j\right)}$), then
    $e^{\left(m\right)}_n = w\left(n, j\right)^{-1}$ and $f^{\left(m\right)}_n = w\left(n-1,j\right)^{-1} - w\left(n,j\right)^{-1}$, hence $d^{\left(m\right)}_n = w\left(n-1,j\right)^{-1}$.
    If he jumps left (with probability
    $\frac{1}{k} \cdot \frac{w\left(n - 1, j - 1\right)}{w\left(n - 1, j - 1\right) + w\left(n - 1, j\right)}$), then
    $e^{\left(m\right)}_n = -w\left(n - 1, j - 1\right)^{-1}$ and $f^{\left(m\right)}_n = 0$, hence $d^{\left(m\right)}_n = -w\left(n - 1, j - 1\right)^{-1}$.
    \item with probability $\frac{k-1}{k}$, the walker $m$ does not jump. In this case, $f^{\left(m\right)}_n = 0$
    since the indicator variable in $f^{\left(m\right)}_n$ is therefore $0$.
    The value of $M^{\left(m\right)}_{n-1}$ now changes
    (that is, $e^{\left(m\right)}_n \neq 0$) if one of the other $k-1$ walkers crosses one of the edges between the nodes $0$ and $j$.
    At the same time, $g^{\left(m\right)}_n \neq 0$ only in this exact case. Now assume the walker $l \neq m$
    crosses the edge $i$ with $0 \leq i < j$. Then $e^{\left(m\right)}_n = \frac{1}{w\left(n,i\right)} - \frac{1}{w\left(n-1,i\right)}$ and
    $g^{\left(m\right)}_n = \frac{1}{w\left(n-1,i\right)} - \frac{1}{w\left(n,i\right)}$, hence $d^{\left(m\right)}_n = 0$.
    \item conditioned on $X^{\left(m\right)}_{n-1} = j > 0$ and $\tau^{\left(m\right)} \geq n$ (both events measurable w.r.t.~$\calG_{n-1}$), we can therefore conclude
    \begin{align*}
      \Ex{d^{\left(m\right)}_n \midd| \calG_{n-1}} &= \frac{1}{k} \cdot \frac{1}{w\left(n - 1, j - 1\right) + w\left(n - 1, j\right)} \cdot \left(
      \frac{w\left(n - 1, j\right)}{w\left(n - 1, j\right)} - \frac{w\left(n - 1, j - 1\right)}{w\left(n - 1, j - 1\right)}\right) \\
      &= \frac{1}{k} \cdot \frac{1}{w\left(n - 1, j - 1\right) + w\left(n - 1, j\right)} \cdot \left(1 - 1\right) = 0
    \end{align*}
    (By the same arguments, but only considering $e^{\left(m\right)}_n$, we can show that $M^{\left(m\right)}_n$
    is a supermartingale.)
  \end{itemize}
  As a nonnegative martingale, $H^{\left(m\right)}_n$ converges almost surely.
  
  \begin{figure}[H]
    \begin{center}
      \begin{tikzpicture}[scale=0.9]
  \draw[dashed] (-4.32,0) -- (-3.6,0);
  \draw[dashed] (10.8,0) -- (11.52,0);
  \draw (-3.6,0) -- node[above] {$1$} (-2.4,0) -- node[above] {$1$} (-1.2,0) -- node[above] {$1$} node[below=3mm] {$M^{\left(3\right)}_0\;=$} (0,0) node[above=1mm,sectionblue] {$0$} -- node[above] {$1$} node[below=3mm] {$\vphantom{M^{\left(m\right)}_n}1$} (1.2,0) -- node[above] {$1$} node[below=3mm] {$\vphantom{M^{\left(m\right)}_n}+\;1$} (2.4,0) -- node[above] {$1$} node[below=3mm] {$\vphantom{M^{\left(m\right)}_n}+\;1$} (3.6,0) -- node[above] {$1$} node[below=3mm] {$\vphantom{M^{\left(m\right)}_n}+\;1$} (4.8,0) -- node[above] {$1$} node[below=3mm] {$\vphantom{M^{\left(m\right)}_n}+\;1$} (6,0) -- node[above] {$1$} node[below=3mm] {$\vphantom{M^{\left(m\right)}_n}=\;5$} (7.2,0) -- node[above] {$1$} (8.4,0) -- node[above] {$1$} (9.6,0) -- node[above] {$1$} (10.8,0);
  \node[circle,fill=black,inner sep=0.7mm] at (-3.6,0) {};
  \node[circle,fill=black,inner sep=0.7mm] at (-2.4,0) {};
  \node[circle,draw=black,inner sep=1.1mm] at (-2.4,0) {};
  \node[circle,fill=black,inner sep=0.7mm] at (-1.2,0) {};
  \node[circle,fill=sectionblue,inner sep=0.7mm] at (0,0) {};
  \node[circle,fill=black,inner sep=0.7mm] at (1.2,0) {};
  \node[circle,fill=black,inner sep=0.7mm] at (2.4,0) {};
  \node[circle,draw=black,inner sep=1.1mm] at (2.4,0) {};
  \node[circle,fill=black,inner sep=0.7mm] at (3.6,0) {};
  \node[circle,fill=black,inner sep=0.7mm] at (4.8,0) {};
  \node[circle,fill=sectionblue,inner sep=0.7mm] at (6,0) {};
  \node[circle,draw=sectionblue,inner sep=1.1mm] at (6,0) {};
  \node[circle,fill=black,inner sep=0.7mm] at (7.2,0) {};
  \node[circle,fill=black,inner sep=0.7mm] at (8.4,0) {};
  \node[circle,draw=black,inner sep=1.1mm] at (8.4,0) {};
  \node[circle,fill=black,inner sep=0.7mm] at (9.6,0) {};
  \node[circle,fill=black,inner sep=0.7mm] at (10.8,0) {};
\end{tikzpicture}
 \\[1em]
      \begin{tikzpicture}[scale=0.9]
  \draw[line width=2pt] (6,0) -- (7.2,0);
  \draw[dashed] (-4.32,0) -- (-3.6,0);
  \draw[dashed] (10.8,0) -- (11.52,0);
  \draw (-3.6,0) -- node[above] {$1$} (-2.4,0) -- node[above] {$1$} (-1.2,0) -- node[above] {$1$} node[below=3mm] {$M^{\left(3\right)}_1\;=$} (0,0) node[above=1mm,sectionblue] {$0$} -- node[above] {$1$} node[below=3mm] {$\vphantom{M^{\left(m\right)}_n}1$} (1.2,0) -- node[above] {$1$} node[below=3mm] {$\vphantom{M^{\left(m\right)}_n}+\;1$} (2.4,0) -- node[above] {$1$} node[below=3mm] {$\vphantom{M^{\left(m\right)}_n}+\;1$} (3.6,0) -- node[above] {$1$} node[below=3mm] {$\vphantom{M^{\left(m\right)}_n}+\;1$} (4.8,0) -- node[above] {$1$} node[below=3mm] {$\vphantom{M^{\left(m\right)}_n}+\;1$} (6,0) -- node[above] {$2$} node[below=3mm] {$\vphantom{M^{\left(m\right)}_n}+\;\frac{1}{2}$} (7.2,0) -- node[above] {$1$} node[below=3mm] {$\vphantom{M^{\left(m\right)}_n}=\;\frac{11}{2}$} (8.4,0) -- node[above] {$1$} (9.6,0) -- node[above] {$1$} (10.8,0);
  \node[circle,fill=black,inner sep=0.7mm] at (-3.6,0) {};
  \node[circle,fill=black,inner sep=0.7mm] at (-2.4,0) {};
  \node[circle,draw=black,inner sep=1.1mm] at (-2.4,0) {};
  \node[circle,fill=black,inner sep=0.7mm] at (-1.2,0) {};
  \node[circle,fill=sectionblue,inner sep=0.7mm] at (0,0) {};
  \node[circle,fill=black,inner sep=0.7mm] at (1.2,0) {};
  \node[circle,fill=black,inner sep=0.7mm] at (2.4,0) {};
  \node[circle,draw=black,inner sep=1.1mm] at (2.4,0) {};
  \node[circle,fill=black,inner sep=0.7mm] at (3.6,0) {};
  \node[circle,fill=black,inner sep=0.7mm] at (4.8,0) {};
  \node[circle,fill=black,inner sep=0.7mm] at (6,0) {};
  \node[circle,fill=sectionblue,inner sep=0.7mm] at (7.2,0) {};
  \node[circle,draw=sectionblue,inner sep=1.1mm] at (7.2,0) {};
  \node[circle,fill=black,inner sep=0.7mm] at (8.4,0) {};
  \node[circle,draw=black,inner sep=1.1mm] at (8.4,0) {};
  \node[circle,fill=black,inner sep=0.7mm] at (9.6,0) {};
  \node[circle,fill=black,inner sep=0.7mm] at (10.8,0) {};
\end{tikzpicture}
 \\[1em]
      \begin{tikzpicture}[scale=0.9]
  \draw[line width=2pt] (1.2,0) -- (2.4,0);
  \draw[dashed] (-4.32,0) -- (-3.6,0);
  \draw[dashed] (10.8,0) -- (11.52,0);
  \draw (-3.6,0) -- node[above] {$1$} (-2.4,0) -- node[above] {$1$} (-1.2,0) -- node[above] {$1$} node[below=3mm] {$M^{\left(3\right)}_2\;=$} (0,0) node[above=1mm,sectionblue] {$0$} -- node[above] {$1$} node[below=3mm] {$\vphantom{M^{\left(m\right)}_n}1$} (1.2,0) -- node[above] {$2$} node[below=3mm] {$\vphantom{M^{\left(m\right)}_n}+\;\frac{1}{2}$} (2.4,0) -- node[above] {$1$} node[below=3mm] {$\vphantom{M^{\left(m\right)}_n}+\;1$} (3.6,0) -- node[above] {$1$} node[below=3mm] {$\vphantom{M^{\left(m\right)}_n}+\;1$} (4.8,0) -- node[above] {$1$} node[below=3mm] {$\vphantom{M^{\left(m\right)}_n}+\;1$} (6,0) -- node[above] {$2$} node[below=3mm] {$\vphantom{M^{\left(m\right)}_n}+\;\frac{1}{2}$} (7.2,0) -- node[above] {$1$} node[below=3mm] {$\vphantom{M^{\left(m\right)}_n}=\;5$} (8.4,0) -- node[above] {$1$} (9.6,0) -- node[above] {$1$} (10.8,0);
  \node[circle,fill=black,inner sep=0.7mm] at (-3.6,0) {};
  \node[circle,fill=black,inner sep=0.7mm] at (-2.4,0) {};
  \node[circle,draw=black,inner sep=1.1mm] at (-2.4,0) {};
  \node[circle,fill=black,inner sep=0.7mm] at (-1.2,0) {};
  \node[circle,fill=sectionblue,inner sep=0.7mm] at (0,0) {};
  \node[circle,fill=black,inner sep=0.7mm] at (1.2,0) {};
  \node[circle,draw=black,inner sep=1.1mm] at (1.2,0) {};
  \node[circle,fill=black,inner sep=0.7mm] at (2.4,0) {};
  \node[circle,fill=black,inner sep=0.7mm] at (3.6,0) {};
  \node[circle,fill=black,inner sep=0.7mm] at (4.8,0) {};
  \node[circle,fill=black,inner sep=0.7mm] at (6,0) {};
  \node[circle,fill=sectionblue,inner sep=0.7mm] at (7.2,0) {};
  \node[circle,draw=sectionblue,inner sep=1.1mm] at (7.2,0) {};
  \node[circle,fill=black,inner sep=0.7mm] at (8.4,0) {};
  \node[circle,draw=black,inner sep=1.1mm] at (8.4,0) {};
  \node[circle,fill=black,inner sep=0.7mm] at (9.6,0) {};
  \node[circle,fill=black,inner sep=0.7mm] at (10.8,0) {};
\end{tikzpicture}
 \\[1em]
      \begin{tikzpicture}[scale=0.9]
  \draw[line width=2pt] (7.2,0) -- (8.4,0);
  \draw[dashed] (-4.32,0) -- (-3.6,0);
  \draw[dashed] (10.8,0) -- (11.52,0);
  \draw (-3.6,0) -- node[above] {$1$} (-2.4,0) -- node[above] {$1$} (-1.2,0) -- node[above] {$1$} node[below=3mm] {$M^{\left(3\right)}_3\;=$} (0,0) node[above=1mm,sectionblue] {$0$} -- node[above] {$1$} node[below=3mm] {$\vphantom{M^{\left(m\right)}_n}1$} (1.2,0) -- node[above] {$2$} node[below=3mm] {$\vphantom{M^{\left(m\right)}_n}+\;\frac{1}{2}$} (2.4,0) -- node[above] {$1$} node[below=3mm] {$\vphantom{M^{\left(m\right)}_n}+\;1$} (3.6,0) -- node[above] {$1$} node[below=3mm] {$\vphantom{M^{\left(m\right)}_n}+\;1$} (4.8,0) -- node[above] {$1$} node[below=3mm] {$\vphantom{M^{\left(m\right)}_n}+\;1$} (6,0) -- node[above] {$2$} node[below=3mm] {$\vphantom{M^{\left(m\right)}_n}+\;\frac{1}{2}$} (7.2,0) -- node[above] {$2$} node[below=3mm] {$\vphantom{M^{\left(m\right)}_n}=\;5$} (8.4,0) -- node[above] {$1$} (9.6,0) -- node[above] {$1$} (10.8,0);
  \node[circle,fill=black,inner sep=0.7mm] at (-3.6,0) {};
  \node[circle,fill=black,inner sep=0.7mm] at (-2.4,0) {};
  \node[circle,draw=black,inner sep=1.1mm] at (-2.4,0) {};
  \node[circle,fill=black,inner sep=0.7mm] at (-1.2,0) {};
  \node[circle,fill=sectionblue,inner sep=0.7mm] at (0,0) {};
  \node[circle,fill=black,inner sep=0.7mm] at (1.2,0) {};
  \node[circle,draw=black,inner sep=1.1mm] at (1.2,0) {};
  \node[circle,fill=black,inner sep=0.7mm] at (2.4,0) {};
  \node[circle,fill=black,inner sep=0.7mm] at (3.6,0) {};
  \node[circle,fill=black,inner sep=0.7mm] at (4.8,0) {};
  \node[circle,fill=black,inner sep=0.7mm] at (6,0) {};
  \node[circle,fill=sectionblue,inner sep=0.7mm] at (7.2,0) {};
  \node[circle,draw=sectionblue,inner sep=1.1mm] at (7.2,0) {};
  \node[circle,draw=black,inner sep=1.5mm] at (7.2,0) {};
  \node[circle,fill=black,inner sep=0.7mm] at (8.4,0) {};
  \node[circle,fill=black,inner sep=0.7mm] at (9.6,0) {};
  \node[circle,fill=black,inner sep=0.7mm] at (10.8,0) {};
\end{tikzpicture}
 \\[1em]
      \begin{tikzpicture}[scale=0.9]
  \draw[line width=2pt] (6,0) -- (7.2,0);
  \draw[dashed] (-4.32,0) -- (-3.6,0);
  \draw[dashed] (10.8,0) -- (11.52,0);
  \draw (-3.6,0) -- node[above] {$1$} (-2.4,0) -- node[above] {$1$} (-1.2,0) -- node[above] {$1$} node[below=3mm] {$M^{\left(3\right)}_4\;=$} (0,0) node[above=1mm,sectionblue] {$0$} -- node[above] {$1$} node[below=3mm] {$\vphantom{M^{\left(m\right)}_n}1$} (1.2,0) -- node[above] {$2$} node[below=3mm] {$\vphantom{M^{\left(m\right)}_n}+\;\frac{1}{2}$} (2.4,0) -- node[above] {$1$} node[below=3mm] {$\vphantom{M^{\left(m\right)}_n}+\;1$} (3.6,0) -- node[above] {$1$} node[below=3mm] {$\vphantom{M^{\left(m\right)}_n}+\;1$} (4.8,0) -- node[above] {$1$} node[below=3mm] {$\vphantom{M^{\left(m\right)}_n}+\;1$} (6,0) -- node[above] {$3$} node[below=3mm] {$\vphantom{M^{\left(m\right)}_n}=\;\frac{9}{2}$} (7.2,0) -- node[above] {$2$} (8.4,0) -- node[above] {$1$} (9.6,0) -- node[above] {$1$} (10.8,0);
  \node[circle,fill=black,inner sep=0.7mm] at (-3.6,0) {};
  \node[circle,fill=black,inner sep=0.7mm] at (-2.4,0) {};
  \node[circle,draw=black,inner sep=1.1mm] at (-2.4,0) {};
  \node[circle,fill=black,inner sep=0.7mm] at (-1.2,0) {};
  \node[circle,fill=sectionblue,inner sep=0.7mm] at (0,0) {};
  \node[circle,fill=black,inner sep=0.7mm] at (1.2,0) {};
  \node[circle,draw=black,inner sep=1.1mm] at (1.2,0) {};
  \node[circle,fill=black,inner sep=0.7mm] at (2.4,0) {};
  \node[circle,fill=black,inner sep=0.7mm] at (3.6,0) {};
  \node[circle,fill=black,inner sep=0.7mm] at (4.8,0) {};
  \node[circle,fill=sectionblue,inner sep=0.7mm] at (6,0) {};
  \node[circle,draw=sectionblue,inner sep=1.1mm] at (6,0) {};
  \node[circle,fill=black,inner sep=0.7mm] at (7.2,0) {};
  \node[circle,draw=black,inner sep=1.1mm] at (7.2,0) {};
  \node[circle,fill=black,inner sep=0.7mm] at (8.4,0) {};
  \node[circle,fill=black,inner sep=0.7mm] at (9.6,0) {};
  \node[circle,fill=black,inner sep=0.7mm] at (10.8,0) {};
\end{tikzpicture}

      \caption{Evolution of the supermartingale during 4 steps of the ERRW with multiple walkers}
      \label{fig:errw_multiple_mart}
    \end{center}
    In the figure above, we consider $k = 4$ walkers whose positions are indicated by the circled nodes. The value
    of $M_n^{\left(3\right)}$ is shown for the first $4$ steps, and the corresponding walker $X^{\left(3\right)}$
    is highlighted in \textcolor{sectionblue}{blue}.
  \end{figure}
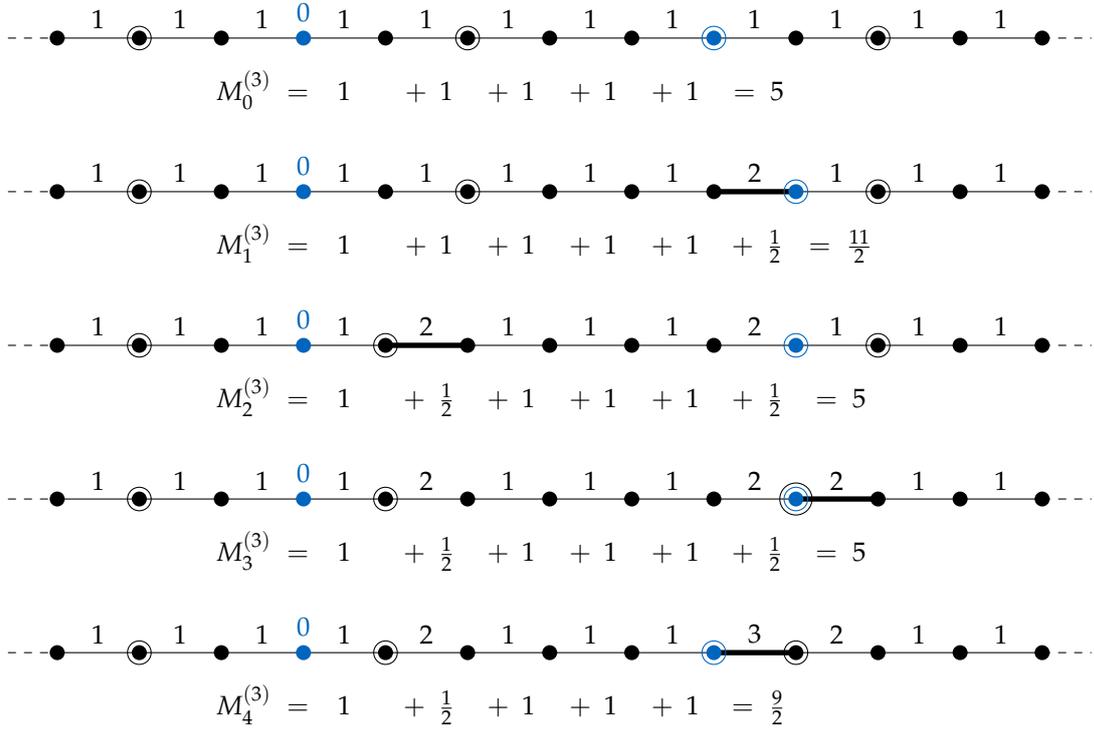

  We just showed this for all walkers $m$ ($1 \leq m \leq k$) with $X^{\left(m\right)}_0 \geq 0$.
  Further observe that for such a walker $m$, we have, on the event $B^{\left(m\right)}_n = \left\{X^{\left(m\right)}_n > X^{\left(m\right)}_{n-1}, n \leq \tau^{\left(m\right)}, w\left(n-1, X^{\left(m\right)}_{n-1}\right) = 1\right\}$,
  that $e^{\left(m\right)}_n = w\left(n, X^{\left(m\right)}_{n-1}\right)^{-1}$, $f^{\left(m\right)}_n = w\left(n - 1, X^{\left(m\right)}_{n-1}\right)^{-1} - w\left(n, X^{\left(m\right)}_{n-1}\right)^{-1}$, $g^{\left(m\right)}_n = 0$ and hence $d^{\left(m\right)}_n = 1$.
  Thus, by convergence, only a finite number of the events $B^{\left(m\right)}_n$ can occur
  for every such walker $m$.

  Now define $\Gamma$ to be the set of edges between two nonnegative integers to the right of the integer
  $\max\left\{ X^{\left(m\right)}_0 : 1 \leq m \leq k \right\}$ for which the initial weight
  was $1$ (all but finitely many edges meet the latter criterion), and further define the event
  \begin{align*}
    D_n &= \left\{ \exists l: X^{\left(l\right)}_0 \geq 0 \textrm{, an edge in } \Gamma \textrm{ is crossed between time } n-1 \textrm{ and } n \textrm{ for the first time by}\right.\\
    &\hphantom{\;=\;\{} \left.\vphantom{X^{\left(l\right)}_0}\textrm{any walker, the crossing walker is }l \textrm{ and } n \leq \tau^{\left(l\right)} \right\}
  \end{align*}
  Clearly, $D_n \subseteq B^{\left(m\right)}_n$ for some random walker $m$ with $X^{\left(m\right)}_0 \geq 0$, hence only a finite number of
  the events $D_n$ can occur.

  Now the proof cannot be continued along \cite[Lemma 3.0]{errwz} since the walkers starting to the left
  of $0$ and the walkers which reach $0$ can
  later cross edges to the right of $0$ without triggering $D_n$ and the other walkers can then follow
  them without triggering $D_n$. So, we only proved that walkers which never go to $0$ and start to
  the right of $0$ cannot
  visit infinitely many edges which have not been visited before by any other walker.
\end{tproof}

\begin{corollary}
  \label{cor:oneReturnsOrFiniteRange}
  Assume the edge-reinforced random walk with $k$ walkers starts with an initial configuration of the
  weights $w\left(0, j\right)$ such that all but finitely many of them are $1$.
  Then, we have the following:
  \begin{align*}
    &\Prb{\exists m: X^{\left(m\right)}_n = 0 \textrm{ for some }n \geq 0} \\
    +\; &\bbP\left[ \forall m: X^{\left(m\right)}_n \neq 0\textrm{ for all }n \geq 0 \textrm{ and } X^{\left(m\right)} \textrm{ has finite range} \right] = 1 \qedhere
  \end{align*}
\end{corollary}

\begin{tproof}
  It is sufficient to show that, conditional on
  $A^{\textrm{C}} = \left\{\nexists m: X^{\left(m\right)}_n = 0 \textrm{ for some }n \geq 0\right\}$, we have that the
  second event $B = \left\{ \forall m: X^{\left(m\right)}_n \neq 0\textrm{ for all }n \geq 0 \textrm{ and } X^{\left(m\right)} \textrm{ has finite range} \right\}$ occurs a.s.
  Now, by \autoref{lem:retRootOrFinManyNew}, every random walker $m$ separately either reaches $0$ or some other
  event $E^{\left(m\right)}$ occurs (by symmetry, \autoref{lem:retRootOrFinManyNew} can also be applied to
  random walkers which start to the left of $0$). But now, if $A^{\textrm{C}}$ occurs, then no random
  walker reaches $0$, hence the event $E^{\left(m\right)}$ occurs for every random walker $m$ where
  \begin{align*}
    E^{\left(m\right)} &= \left\{ X^{\left(m\right)}_n \neq 0\textrm{ for all }n \geq 0 \textrm{ and } X^{\left(m\right)} \textrm{ only visits finitely many nodes}\right. \\
    &\hphantom{\;=\;\{}\left.\vphantom{X^{\left(m\right)}_n}\textrm{which have not been visited before by any other walker} \right\}
  \end{align*}
  So every walker visits only finitely many nodes not visited before by any of the other walkers. But
  this implies that all walkers together can only visit finitely many nodes.
\end{tproof}

\begin{corollary}
  \label{cor:visitedInfOftenOrFiniteRange}
  Assume the edge-reinforced random walk with $k$ walkers starts with an initial configuration of the
  weights $w\left(0, j\right)$ such that all but finitely many of them are $1$.
  Then, we have the following:
  \begin{align*}
    &\Prb{\forall j \in \bbZ: j \textrm{ is visited }\infty\textrm{ often by at least one of the walkers}} \\
    +\; &\bbP\left[ \textrm{all walkers have finite range} \right] \qquad\qquad\qquad\qquad\qquad\qquad\;\; =\; 1 \qedhere
  \end{align*}
\end{corollary}

\begin{tproof}
  It suffices to show that for every $j \in \bbZ$, we have
  \begin{align}
    \begin{split}
    \Prb{j \textrm{ is visited }\infty\textrm{ often by at least one of the random walkers}} & \\
    + \;\bbP\left[j \textrm{ is visited only finitely often and all random walkers have finite range} \right] &= 1
    \end{split}
    \label{eq:singleJ}
  \end{align}
  To see this, assume that not all random walkers have finite range. Then, if we have proved
  \autoref{eq:singleJ}, we know that for every $j \in \bbZ$, $j$ is visited infinitely often by
  one of the random walkers. Thus, conditional on not all random walkers having finite range, all nodes
  are visited infinitely often, which is equivalent to \autoref{cor:visitedInfOftenOrFiniteRange}.

  To show \autoref{eq:singleJ}, it suffices in turn to show, for every $j \in \bbZ$ and $\forall n \geq 0$:
  \begin{align}
    \begin{split}
    \Prb{j \textrm{ is visited by at least one walker at a time }t \geq n}&\\
    + \;\Prb{\textrm{no walker visits }j\textrm{ at a time }t \geq n\textrm{ and all walkers have finite range}} &= 1
    \end{split}
    \label{eq:corReformulated}
  \end{align}
  since the given events are decreasing and increasing respectively, and their limits correspond to the
  events in \autoref{eq:singleJ}.

  But now, consider the random walkers $\left(X^{\left(m\right)}_{n+i}\right)_{i \geq 0}$. These form
  again an edge-reinforced random walk with $k$ walkers, and since until time $n$, only a finite number
  of edge weights can have changed, all but finitely many edges will still have weight $1$ at time $n$.
  Thus, we can apply \autoref{cor:oneReturnsOrFiniteRange} to the walkers $\left(X^{\left(m\right)}_{n+i}\right)_{i \geq 0}$,
  which directly proves \autoref{eq:corReformulated} (\autoref{cor:oneReturnsOrFiniteRange} only proves
  \autoref{eq:corReformulated} for $j = 0$ but of course we can relabel the nodes such that any other
  node gets the label $0$, hence \autoref{cor:oneReturnsOrFiniteRange} is valid for any choice of node).
\end{tproof}

\begin{lemma}
  \label{lem:labelExchange}
  Assume the edge-reinforced random walk with $k \geq 2$ walkers starts with an arbitrary initial configuration of the
  weights $w\left(0, j\right)$.
  Assume further that $X^{\left(1\right)}$ and $X^{\left(2\right)}$ meet infinitely often. Then (almost surely):
  \begin{enumerate}[(i)]
    \item \label{lem:labelExchangeBothInf} if at least one of the walkers $X^{\left(1\right)}$ and $X^{\left(2\right)}$ does
    not have finite range, then both $X^{\left(1\right)}$ and $X^{\left(2\right)}$ do not have finite range.
    \item \label{lem:labelExchangeBothVisInf} if some integer $z$ is visited infinitely often by at least one of the walkers $X^{\left(1\right)}$
    and $X^{\left(2\right)}$, then both $X^{\left(1\right)}$ and $X^{\left(2\right)}$ visit $z$ infinitely often.
    \item \label{lem:labelExchangeBothRec} if every integer is visited infinitely often by at least one of the walkers $X^{\left(1\right)}$
    and $X^{\left(2\right)}$, then both $X^{\left(1\right)}$ and $X^{\left(2\right)}$ are recurrent. \qedhere
  \end{enumerate}
\end{lemma}

\begin{tproof}
  The proof idea is the following: whenever $X^{\left(1\right)}$ and $X^{\left(2\right)}$ meet, we
  can randomly exchange their labels, i.e.~we can randomly decide whether we want to rename $X^{\left(1\right)}$
  to $X^{\left(2\right)}$ and vice versa, and the law of the edge-reinforced random walk with the
  two walkers is invariant under such relabelings because the only distinguishing feature of a random
  walker is his position. But now, to construct counterexamples to the two statements in \autoref{lem:labelExchange},
  we would have to choose a fixed labeling for infinitely many times at which the walkers meet. But
  if we randomize the labeling with a sequence of independent Bernoulli random variables, then the probability of
  choosing a certain fixed labeling at infinitely many points in the sequence is $0$, and since the
  law was invariant under random relabeling, it follows that the probability of any such counterexample
  is $0$. We continue with the formal proof.

  Set $\tau_1 := \inf\left\{ n \geq 0 : X^{\left(1\right)}_n = X^{\left(2\right)}_n \right\}$
  and $\tau_{i + 1} := \inf\left\{ n > \tau_i : X^{\left(1\right)}_n = X^{\left(2\right)}_n \right\}$.
  If $X^{\left(1\right)}, X^{\left(2\right)}$ meet infinitely
  often, then $\forall n: \tau_n < \infty$, but the construction also works if this is not the case. Let $\left(\omega_i\right)_{i \geq 1}$ be a sequence of
  iid random variables with $\Prb{\omega_i = 1} = \frac{1}{2} = \Prb{\omega_i = 0}$ (the $\omega_i$
  are also independent of $\calG_n$ for all $n$, i.e.~independent of the edge-reinforced random walk).
  Define $\widetilde{X}^{\left(1\right)}_n$ and $\widetilde{X}^{\left(2\right)}_n$ as follows (with $\omega_0 = 0$ and $\tau_0 = -1$):
  \begin{align*}
    \widetilde{X}^{\left(1\right)}_n := X^{\left(1,\omega\right)}_n &= \sum_{i \geq 0} \left( \left(1 - \omega_i\right) X^{\left(1\right)}_n + \omega_i X^{\left(2\right)}_n \right) \cdot
    \mathbbm{1}_{\tau_i < n \leq \tau_{i+1}} \\
    \widetilde{X}^{\left(2\right)}_n := X^{\left(2,\omega\right)}_n &= \sum_{i \geq 0} \left( \left(1 - \omega_i\right) X^{\left(2\right)}_n + \omega_i X^{\left(1\right)}_n \right) \cdot
    \mathbbm{1}_{\tau_i < n \leq \tau_{i+1}}
  \end{align*}
  Note that the sums collapse to a single term. $\omega_i = 1$ means that we switch the labels of $X^{\left(1\right)}$
  and $X^{\left(2\right)}$ during the time interval $\left(\tau_i,\tau_{i+1}\right]$. If we consider
  $\left(X^{\left(1\right)},X^{\left(2\right)}\right)$ and $\left(\widetilde{X}^{\left(1\right)},\widetilde{X}^{\left(2\right)}\right)$
  as sequences of pairs of integers, then we have
  \begin{align}
    \left(\widetilde{X}^{\left(i\right)}_n\right)_{1\leq i\leq 2, n\geq 0} &\overset{\textrm{d}}{=}
    \left(X^{\left(i\right)}_n\right)_{1\leq i\leq 2, n\geq 0}
    \label{eq:eqInDist}
  \end{align}
  The equality in distribution follows from the mentioned invariance of the law of the random walk
  under relabelings at meeting times which is quite intuitive, and could be proved formally by looking at cylinder events,
  for example.

  \begin{figure}[H]
    \begin{center}
      \begin{tikzpicture}[scale=0.8]
  \draw[-{Latex[length=1.8mm,width=1.8mm]},tumGray] (0,0) -- (15.5,0) node[right] {time};
  \draw[sectionblue,line width=1pt] (0,-1) -- (0.5,-1) -- (1,-1) -- (1.5,-1) -- (2,-1) -- (2.5,-1) -- (3,-1.5) -- (3.5,-2) -- (4,-2.5) -- (4.5,-2.5) -- (5,-2.5) -- (5.5,-2) -- (6,-2.5) -- (6.5,-2.5) -- (7,-2.5) -- (7.5,-2) -- (8,-1.5) -- (8.5,-1.5) -- (9,-1.5) -- (9.5,-1.5) -- (10,-1.5) -- (10.5,-1.5) -- (11,-1) -- (11.5,-1) -- (12,-1) -- (12.5,-1) -- (13,-0.5) -- (13.5,-1) -- (14,-1.5) -- (14.5,-1) -- (15,-0.5) node[right] {$X^{\left(1\right)}$};
  \draw[tumOrange,line width=1pt] (0,1) -- (0.5,0.5) -- (1,0) -- (1.5,-0.5) -- (2,-1) -- (2.5,-0.5) -- (3,-0.5) -- (3.5,-0.5) -- (4,-0.5) -- (4.5,-1) -- (5,-1.5) -- (5.5,-1.5) -- (6,-1.5) -- (6.5,-1) -- (7,-1.5) -- (7.5,-1.5) -- (8,-1.5) -- (8.5,-2) -- (9,-1.5) -- (9.5,-1) -- (10,-0.5) -- (10.5,-1) -- (11,-1) -- (11.5,-1.5) -- (12,-2) -- (12.5,-1.5) -- (13,-1.5) -- (13.5,-1.5) -- (14,-1.5) -- (14.5,-1.5) -- (15,-1.5) node[right] {$X^{\left(2\right)}$};
  \draw[-{Latex[length=1.8mm,width=1.8mm]}] (0,-3) -- (0,1.5) node[above] {$\bbZ$};
  \draw (0,0) -- (-0.2,0) node[left] {$0$};
  \draw (0,-2.5) -- (-0.2,-2.5);
  \draw (0,-2) -- (-0.2,-2);
  \draw (0,-1.5) -- (-0.2,-1.5);
  \draw (0,-1) -- (-0.2,-1);
  \draw (0,-0.5) -- (-0.2,-0.5);
  \draw (0,0) -- (-0.2,0);
  \draw (0,0.5) -- (-0.2,0.5);
  \draw (0,1) -- (-0.2,1);
  \node[circle,fill=black,inner sep=0.7mm] at (2,-1) {};
  \node[circle,fill=black,inner sep=0.7mm] at (8,-1.5) {};
  \node[circle,fill=black,inner sep=0.7mm] at (9,-1.5) {};
  \node[circle,fill=black,inner sep=0.7mm] at (11,-1) {};
  \node[circle,fill=black,inner sep=0.7mm] at (14,-1.5) {};
  \node[below=2mm] at (2,-1) {$\tau_1 = 4$};
  \node[above=2mm] at (8,-1.5) {$\tau_2 = 16$};
  \node[below=2mm] at (9,-1.5) {$\tau_3 = 18$};
  \node[above=2mm] at (11,-1) {$\tau_4 = 22$};
  \node[below=2mm] at (14,-1.5) {$\tau_5 = 28$};
\end{tikzpicture}
 \\[1em]
      \begin{tikzpicture}[scale=0.8]
  \fill[tumGray!20!white] (2,-3) rectangle ++(6,4.5);
  \fill[tumGray!20!white] (9,-3) rectangle ++(2,4.5);
  \fill[tumGray!20!white] (14,-3) rectangle ++(1.15,4.5);
  \draw[-{Latex[length=1.8mm,width=1.8mm]},tumGray] (0,0) -- (15.5,0) node[right] {time};
  \draw[sectionblue!50!white,line width=4pt] (0,-1) -- (0.5,-1) -- (1,-1) -- (1.5,-1) -- (2,-1) -- (2.5,-0.5) -- (3,-0.5) -- (3.5,-0.5) -- (4,-0.5) -- (4.5,-1) -- (5,-1.5) -- (5.5,-1.5) -- (6,-1.5) -- (6.5,-1) -- (7,-1.5) -- (7.5,-1.5) -- (8,-1.5) -- (8.5,-1.5) -- (9,-1.5) -- (9.5,-1) -- (10,-0.5) -- (10.5,-1) -- (11,-1) -- (11.5,-1) -- (12,-1) -- (12.5,-1) -- (13,-0.5) -- (13.5,-1) -- (14,-1.5) -- (14.5,-1.5) -- (15,-1.5) node[right] {$\widetilde{X}^{\left(1\right)}$};
  \draw[tumOrange!50!white,line width=4pt] (0,1) -- (0.5,0.5) -- (1,0) -- (1.5,-0.5) -- (2,-1) -- (2.5,-1) -- (3,-1.5) -- (3.5,-2) -- (4,-2.5) -- (4.5,-2.5) -- (5,-2.5) -- (5.5,-2) -- (6,-2.5) -- (6.5,-2.5) -- (7,-2.5) -- (7.5,-2) -- (8,-1.5) -- (8.5,-2) -- (9,-1.5) -- (9.5,-1.5) -- (10,-1.5) -- (10.5,-1.5) -- (11,-1) -- (11.5,-1.5) -- (12,-2) -- (12.5,-1.5) -- (13,-1.5) -- (13.5,-1.5) -- (14,-1.5) -- (14.5,-1) -- (15,-0.5) node[right] {$\widetilde{X}^{\left(2\right)}$};
  \draw[sectionblue,line width=1pt] (0,-1) -- (0.5,-1) -- (1,-1) -- (1.5,-1) -- (2,-1) -- (2.5,-1) -- (3,-1.5) -- (3.5,-2) -- (4,-2.5) -- (4.5,-2.5) -- (5,-2.5) -- (5.5,-2) -- (6,-2.5) -- (6.5,-2.5) -- (7,-2.5) -- (7.5,-2) -- (8,-1.5) -- (8.5,-1.5) -- (9,-1.5) -- (9.5,-1.5) -- (10,-1.5) -- (10.5,-1.5) -- (11,-1) -- (11.5,-1) -- (12,-1) -- (12.5,-1) -- (13,-0.5) -- (13.5,-1) -- (14,-1.5) -- (14.5,-1) -- (15,-0.5);
  \draw[tumOrange,line width=1pt] (0,1) -- (0.5,0.5) -- (1,0) -- (1.5,-0.5) -- (2,-1) -- (2.5,-0.5) -- (3,-0.5) -- (3.5,-0.5) -- (4,-0.5) -- (4.5,-1) -- (5,-1.5) -- (5.5,-1.5) -- (6,-1.5) -- (6.5,-1) -- (7,-1.5) -- (7.5,-1.5) -- (8,-1.5) -- (8.5,-2) -- (9,-1.5) -- (9.5,-1) -- (10,-0.5) -- (10.5,-1) -- (11,-1) -- (11.5,-1.5) -- (12,-2) -- (12.5,-1.5) -- (13,-1.5) -- (13.5,-1.5) -- (14,-1.5) -- (14.5,-1.5) -- (15,-1.5);
  \draw[-{Latex[length=1.8mm,width=1.8mm]}] (0,-3) -- (0,1.5) node[above] {$\bbZ$};
  \draw (0,0) -- (-0.2,0) node[left] {$0$};
  \draw (0,-2.5) -- (-0.2,-2.5);
  \draw (0,-2) -- (-0.2,-2);
  \draw (0,-1.5) -- (-0.2,-1.5);
  \draw (0,-1) -- (-0.2,-1);
  \draw (0,-0.5) -- (-0.2,-0.5);
  \draw (0,0) -- (-0.2,0);
  \draw (0,0.5) -- (-0.2,0.5);
  \draw (0,1) -- (-0.2,1);
  \node[circle,fill=black,inner sep=0.7mm] at (2,-1) {};
  \node[circle,fill=black,inner sep=0.7mm] at (8,-1.5) {};
  \node[circle,fill=black,inner sep=0.7mm] at (9,-1.5) {};
  \node[circle,fill=black,inner sep=0.7mm] at (11,-1) {};
  \node[circle,fill=black,inner sep=0.7mm] at (14,-1.5) {};
  \node[below=2mm] at (2,-1) {$\tau_1 = 4$};
  \node[above=2mm] at (8,-1.5) {$\tau_2 = 16$};
  \node[below=2mm] at (9,-1.5) {$\tau_3 = 18$};
  \node[above=2mm] at (11,-1) {$\tau_4 = 22$};
  \node[below=2mm] at (14,-1.5) {$\tau_5 = 28$};
  \node at (1,1.25) {$\omega_0 = 0$};
  \node at (5,1.25) {$\omega_1 = 1$};
  \node at (8.5,1.25) {$\omega_2 = 0$};
  \node at (10,1.25) {$\omega_3 = 1$};
  \node at (12.5,1.25) {$\omega_4 = 0$};
  \node at (14.58,1.25) {$\omega_5 = 1$};
\end{tikzpicture}

      \caption{Label exchange lemma}
      \label{fig:label_exchange}
    \end{center}
    The figure above illustrates how the ``label exchange'' of the two walkers works. Two sample paths
    for the two walkers are drawn, together with meeting points and the result of the label exchange.
  \end{figure}
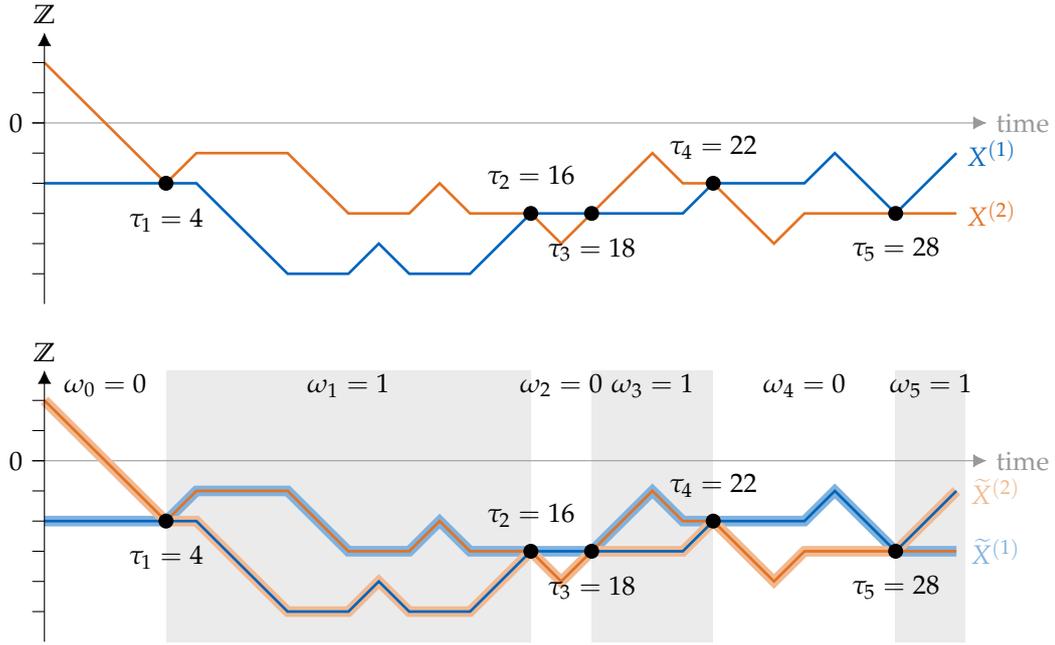
  
  We now show that any counterexamples to statements \ref{lem:labelExchangeBothInf} or \ref{lem:labelExchangeBothRec}
  have probability $0$:
  \begin{enumerate}[(i)]
    \item let $A$ be the event that one of the walkers $X^{\left(1\right)},X^{\left(2\right)}$ has finite range while the other one has infinite
    range, and that they meet infinitely often. It suffices to show that $\Prb{A} = 0$. Denote by $\bP$ the probability measure induced
    by the edge reinforced random walk alone and by $\bQ$ the probability measure induced by the
    sequence $\left(\omega_i\right)_{i \geq 1}$ alone. Then, by \autoref{eq:eqInDist}, we have
    \begin{align*}
      \Prb{A} = &\int \int \mathbbm{1}_{B} \dx{\bQ} \dx{\bP} \\
      \textrm{where } B := &\left\{\textrm{one of }\widetilde{X}^{\left(1\right)},\widetilde{X}^{\left(2\right)}\textrm{ has finite range while}\right.\\
      & \left.\hphantom{\;\{}\vphantom{\widetilde{X}^{\left(1\right)}}\textrm{the other has infinite range, they meet infinitely often}\right\}
    \end{align*}
    We have to show that the inner integral is $0$ almost surely with respect to $\bP$. Consider
    fixed walker sequences $X^{\left(1\right)}$ and $X^{\left(2\right)}$. If one of
    $\widetilde{X}^{\left(1\right)},\widetilde{X}^{\left(2\right)}$ should have finite range while the
    other has infinite range, then, by definition of $\widetilde{X}^{\left(1\right)}$ and $\widetilde{X}^{\left(2\right)}$, at
    least one of $X^{\left(1\right)},X^{\left(2\right)}$ must have infinite range. Of course, by definition,
    we also have that $\widetilde{X}^{\left(1\right)},\widetilde{X}^{\left(2\right)}$ meet infinitely often
    if, and only if, $X^{\left(1\right)},X^{\left(2\right)}$ meet infinitely often. Hence, the indicator variable in the
    integral above can only be $1$ in the case where one of the walkers $X^{\left(1\right)},X^{\left(2\right)}$ has infinite range and the two walkers
    meet infinitely often, so we only need to show that in this particular case, the inner integral is
    still $0$ almost surely.

    Assume $X^{\left(1\right)}$ does not have finite range (w.l.o.g.). Then, for every $n$, one can
    find $i$ such that between times $\tau_i$ and $\tau_{i+1}$ (all $\tau_i$ are finite if the two
    walkers meet infinitely often), $X^{\left(1\right)}$ visits a node at distance at least $n$ from
    the integer $0$. Call these times $\tau_{i_n}$ with $i_n$ strictly increasing in $n$ (w.l.o.g.).
    
    Now consider the walkers $\widetilde{X}^{\left(1\right)},\widetilde{X}^{\left(2\right)}$. One of
    them can have finite range only if the following holds. The same argument works for both walkers,
    we do it here for $\widetilde{X}^{\left(1\right)}$ w.l.o.g. $\widetilde{X}^{\left(1\right)}$
    can only have finite range if there exists $N$ such that for all $n \geq N$ we have $\omega_{i_n} = 1$.
    Assume to the contrary that no such $N$ exists. Then we can find arbitrarily large $n$ such that
    $\omega_{i_n} = 0$ which means that the labels of $X^{\left(1\right)}$ and $X^{\left(2\right)}$
    are \emph{not} exchanged in the interval $\left(\tau_{i_n},\tau_{i_n+1}\right]$. Since $X^{\left(1\right)}$ visits a node at distance at
    least $n$ from $0$ in this time interval, the same holds then for $\widetilde{X}^{\left(1\right)}$,
    so $\widetilde{X}^{\left(1\right)}$ would not have finite range.

    But the probability that the sequence $\omega_{i_n}$ is $1$ for all $n \geq N$ is $0$ for any $N$ (since the
    choice of $i_n$ only depends on the edge-reinforced random walk, i.e.~is independent of the $\omega_i$,
    and since the probability of $\omega$ being constantly $1$ on any fixed infinite subset of the
    integers is $0$ by the choice of $\omega$). Hence, the probability that such $N$ exists is $0$,
    and therefore the probability that $\widetilde{X}^{\left(1\right)}$ has finite range is $0$ as well,
    and the same arguments give that the probability for $\widetilde{X}^{\left(2\right)}$ having finite
    range is $0$ as well (both with respect to the measure $\bQ$).

    So the indicator variable in the integral above is $0$ almost surely w.r.t.~$\bQ$, and hence the
    inner integral is always $0$, which implies that the outer integral is also $0$ and hence $\Prb{A} = 0$.
    \item similar. Let $A$ now be the event that the integer $z$ is visited infinitely often by at least one of
    the walkers $X^{\left(1\right)},X^{\left(2\right)}$, that they meet infinitely often, and that one of
    them does \emph{not} visit $z$ infinitely often. Then, we have again:
    \begin{align*}
      \Prb{A} = &\int \int \mathbbm{1}_{B} \dx{\bQ} \dx{\bP} \\
      \textrm{where } B := &\left\{z\textrm{ visited }\infty\textrm{ often by at least one of }\widetilde{X}^{\left(1\right)},\widetilde{X}^{\left(2\right)}\textrm{,}\right.\\
      & \left.\hphantom{\;\{}\vphantom{\widetilde{X}^{\left(1\right)}}\textrm{they meet }\infty\textrm{ often, one of them visits }z\textrm{ only finitely often}\right\}
    \end{align*}
    We see that the indicator variable can be $1$ only if at least one of $X^{\left(1\right)},X^{\left(2\right)}$
    visits $z$ infinitely often, and w.l.o.g.~assume that this holds for $X^{\left(1\right)}$. As before,
    we can construct a stricly increasing sequence $i_n$ such that in the time interval $\left(\tau_{i_n},\tau_{i_n+1}\right]$,
    $X^{\left(1\right)}$ visits $z$. Again as before, one of $\widetilde{X}^{\left(1\right)},\widetilde{X}^{\left(2\right)}$,
    take $\widetilde{X}^{\left(1\right)}$ w.l.o.g., can visit $z$ only finitely often only if
    $\omega_{i_n} = 1$ for all $n \geq N$ for some $N$, an event which has again probability $0$ w.r.t.~$\bQ$.
    \item apply \ref{lem:labelExchangeBothVisInf} to every integer. \qedhere
  \end{enumerate}
\end{tproof}

\begin{definition}
  \label{def:limInfSup}
  Consider the edge-reinforced random walk with $k$ walkers. Define, for $1 \leq m \leq k$:
  \begin{itemize}
    \item $\overline{X}^{\left(m\right)} := \limsup_{n \to \infty} X^{\left(m\right)}_n$
    \item $\underline{X}^{\left(m\right)} := \liminf_{n \to \infty} X^{\left(m\right)}_n$ \qedhere
  \end{itemize}
\end{definition}

\begin{theorem}
  \label{thm:allRecurrentOrAllFiniteRange}
  Assume the edge-reinforced random walk with $k$ walkers starts with an initial configuration of the
  weights $w\left(0, j\right)$ such that all but finitely many of them are $1$.
  Then, we have the following:
  \begin{align*}
    \Prb{\forall m: X^{\left(m\right)} \textrm{ is recurrent}}
    + \Prb{\forall m: X^{\left(m\right)} \textrm{ has finite range}} &= 1 \qedhere
  \end{align*}
\end{theorem}

\begin{tproof}
  We present two variants of the proof. The first is a little less formal than the second, but hopefully
  easier to understand. The second variant is as formal as possible without becoming totally incomprehensible.

  \textbf{Proof Variant 1} (less formal):

  We have to show the following: if at least one of the walkers does not have finite range, then,
  almost surely, all of them are recurrent. Showing recurrence of all walkers is equivalent to showing
  that $\overline{X}^{\left(m\right)} = \infty, \underline{X}^{\left(m\right)} = -\infty$ for all walkers
  $m$. Now, if at least one walker does not have finite range, then we know by \autoref{cor:visitedInfOftenOrFiniteRange}
  that every integer is visited infinitely often by at least one of the walkers. Hence, there must
  be walkers $m_1$ and $m_2$ with $\overline{X}^{\left(m_1\right)} = \infty, \underline{X}^{\left(m_2\right)} = -\infty$.
  For a contradiction, assume there is some walker which is not recurrent. w.l.o.g.~we assume that there
  is some walker $m_3$ with $\underline{X}^{\left(m_3\right)} > -\infty$ (the proof is the same if the
  other condition on the $\limsup$ is not met). Since we have $\underline{X}^{\left(m_2\right)} = -\infty$,
  we can partition the set of walkers into two non-empty sets:
  \begin{align*}
    P_1 &:= \left\{ m: 1 \leq m \leq k \textrm{ and } \underline{X}^{\left(m\right)} = -\infty \right\} \\
    P_2 &:= \left\{ m: 1 \leq m \leq k \textrm{ and } \underline{X}^{\left(m\right)} > -\infty \right\}
  \end{align*}
  Choose $m_4 \in \argmax_{m \in P_1} \overline{X}^{\left(m\right)}$ and $m_5 \in \argmin_{m \in P_2} \underline{X}^{\left(m\right)}$.
  Let $y := \min\left\{1, \underline{X}^{\left(m_5\right)}\right\} - 1 \in \bbZ$.
  By choice of $m_5$, all walkers in $P_2$ visit $y$ only finitely often. Therefore, by \autoref{lem:retRootOrFinManyNew},
  each of them can only visit finitely many nodes not visited before by any other walker. We know that
  all nodes are visited infinitely often by at least one of the walkers, so we must have
  $\overline{X}^{\left(m_4\right)} = \infty$. If this was not the case, the walkers in $P_1$ would
  only visit nodes to the left of some fixed integer, and since the walkers in $P_2$ together only
  visit finitely many new nodes, this would imply that there is some largest visited integer, a
  contradiction to the fact that every integer is visited infinitely often (since at least one walker
  does not have finite range).

  Now consider $P_3 := \left\{ m \in P_1 : \overline{X}^{\left(m\right)} = \infty \right\} \neq \varnothing$
  (since $m_4 \in P_3$) and the walker $m_5$. If $m_5$ has finite range, then it is clear that $m_5$
  will meet any walker in $P_3$ infinitely often, but this is a contradiction to \autoref{lem:labelExchange}
  \ref{lem:labelExchangeBothInf}. So $m_5$ must have infinite range. This is only possible if
  $\overline{X}^{\left(m_5\right)} = \infty$. But since the walkers in $P_2$, including $m_5$, only
  visit finitely many nodes not visited before by any other walker, at least one walker $m_6$ in $P_3$ must be to
  the right of $m_5$ infinitely often in order to ``free the path'' for $m_5$. As the walkers in $P_3$
  are all recurrent, and as $m_5$ only visits nodes to the right of $y$, this implies that $m_5$ meets
  this walker $m_6$ infinitely often, which is a contradiction to \autoref{lem:labelExchange} \ref{lem:labelExchangeBothRec}.
  Hence, we almost surely arrive at a contradiction.

  To summarize: assuming that at least one walker does not have finite range and that least one walker
  is not recurrent at the same time leads (almost surely) to a contradiction. Therefore, the assumption
  that this can happen must be wrong (almost surely). This implies that if at least one walker does
  not have finite range, then (almost surely) all walkers must be recurrent.

  \begin{figure}[H]
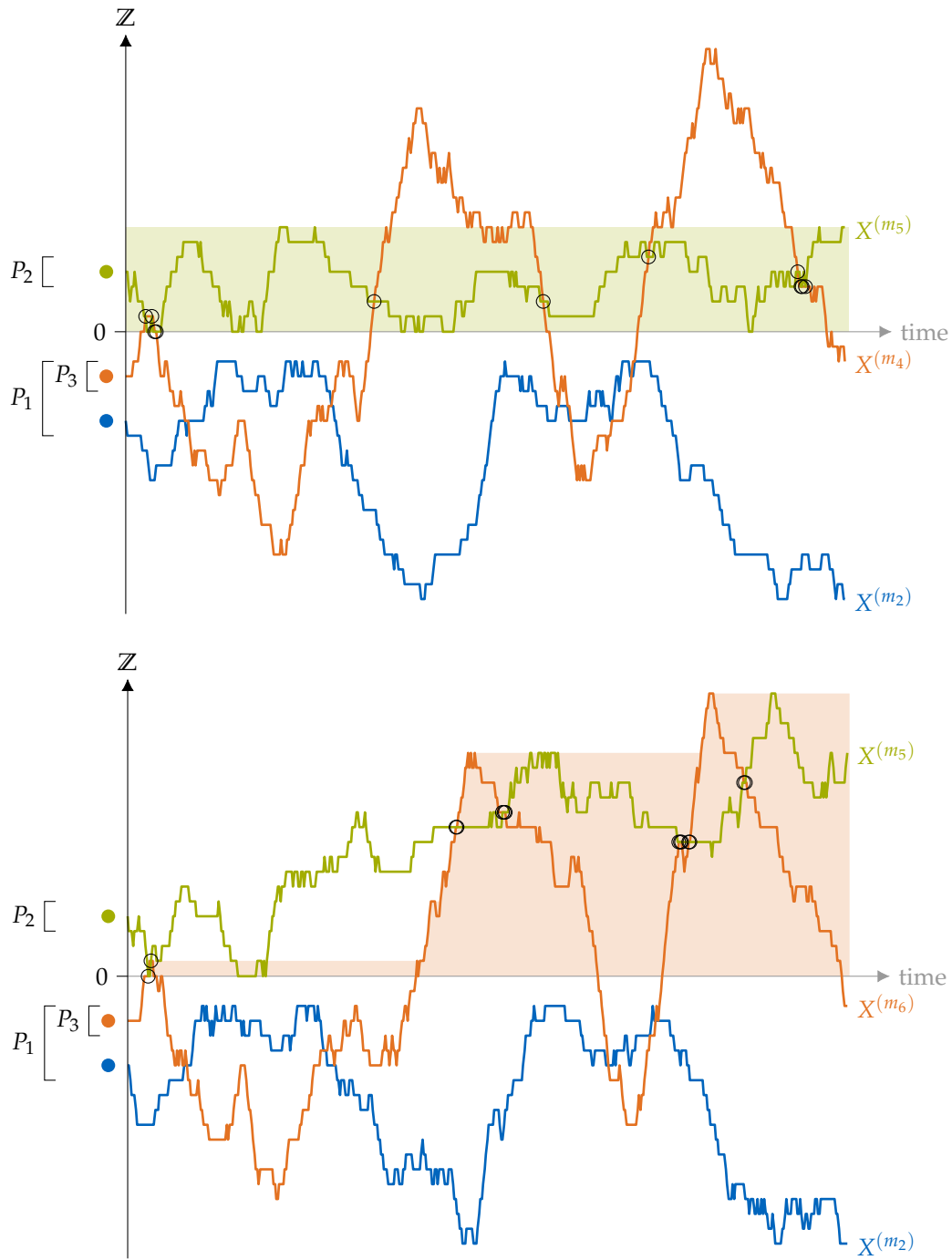

    \begin{center}

      \caption{Illustration: proof that all walkers recurrent or all have finite range}
      \label{fig:rec_proof}
    \end{center}
    Above, the two possible behaviors of the walker $m_5$ (which are both, in fact, almost surely
    impossible) are shown: either $m_5$ has finite range (indicated by the \textcolor{tumGreen}{light green} background),
    but then it would meet $m_4$ infinitely often,
    or $m_5$ does not have finite range, but then $m_6$ would have to ``free the path'' for $m_5$
    (indicated by the \textcolor{tumOrange}{light orange} background) and $m_5$ would meet $m_6$ infinitely often. The relevant
    meeting points are circled.
  \end{figure}
  
  \textbf{Proof Variant 2} (more formal):

  Our goal is to show the following:
  \begin{align*}
    \bbP\left[\vphantom{X^{\left(m_1\right)}}\right.
    \underbrace{\exists m_0, m_3: X^{\left(m_0\right)} \textrm{ does not have finite range and } X^{\left(m_3\right)} \textrm{ is not recurrent}}_{=: A}
    \left.\vphantom{X^{\left(m_1\right)}}\right]
    = 0
  \end{align*}
  The indices of $m$ are chosen in such a way that they agree with proof variant 1, and they therefore
  do not appear in any logical order in this proof variant.
  The proof proceeds by showing that the above event $A$ is subset of a null set. We first apply
  \autoref{cor:visitedInfOftenOrFiniteRange}. Since the intersection of $A$ and the event $\left\{\textrm{all walkers have finite range}\right\}$
  is empty, we can conclude that
  \begin{align*}
    A \subseteq &
    \underbrace{\left(A \cap \left\{\textrm{every integer is visited }\infty\textrm{ often by at least one of the walkers}\right\}\right)}_{B} \cup\; N \\
    &\textrm{where } N \textrm{ is a null set}
  \end{align*}
  Hence, it suffices to show that $\Prb{B} = 0$. $B$ can only occur if there are walkers $m_1$ and $m_2$
  with $\overline{X}^{\left(m_1\right)} = \infty, \underline{X}^{\left(m_2\right)} = -\infty$ since every
  integer is visited infinitely often (when $B$ occurs). Furthermore, for the non-recurrent walker $m_3$,
  we necessarily have $\overline{X}^{\left(m_3\right)} \neq \infty$ or $\underline{X}^{\left(m_3\right)} \neq -\infty$
  (when $B$ occurs). Hence,
  \begin{align*}
  B \subseteq & \underbrace{\left\{ \exists m_1, m_2, m_3: \overline{X}^{\left(m_3\right)} < \infty,
    \overline{X}^{\left(m_1\right)} = \infty, \underline{X}^{\left(m_2\right)} = -\infty \right\}}_{=: C} \\
    &\cup \underbrace{\left\{ \exists m_1, m_2, m_3: \underline{X}^{\left(m_3\right)} > -\infty,
    \overline{X}^{\left(m_1\right)} = \infty, \underline{X}^{\left(m_2\right)} = -\infty \right\}}_{=: D}
  \end{align*}
  We show w.l.o.g.~that $\Prb{D} = 0$ since the proof for $\Prb{C} = 0$ is the same. When $D$ occurs,
  the two sets $P_1 := \left\{ m: 1 \leq m \leq k \textrm{ and } \underline{X}^{\left(m\right)} = -\infty \right\}$
  and $P_2 := \left\{ m: 1 \leq m \leq k \textrm{ and } \underline{X}^{\left(m\right)} > -\infty \right\}$ are
  both non-empty (these are both random set). Choosing $m_5 \in \argmin_{m \in P_2} \underline{X}^{\left(m\right)}$ as well as
  $y := \min\left\{1, \underline{X}^{\left(m_5\right)}\right\} - 1 \in \bbZ$, we see that there is
  a random integer $y$ which is only visited finitely often by the walkers in the set $P_2$. Hence,
  \begin{align*}
    D \subseteq
    \bigcup_{y \in \bbZ} \bigcup_{n \in \bbN}
    &\overbrace{\left\{ \vphantom{\underline{X}^{\left(m_2\right)}} P_1 \neq \varnothing, P_2 \neq \varnothing, y \textrm{ not visited by walkers in } P_2 \textrm{ after time } n,\right.}^{=: E_{y,n}}\\
    &\left. \hphantom{\{} \qquad\qquad\qquad\qquad\qquad\;\;
    \exists m_1, m_2: \overline{X}^{\left(m_1\right)} = \infty, \underline{X}^{\left(m_2\right)} = -\infty \right\}
  \end{align*}
  It thus suffices to show $\Prb{E_{y,n}} = 0$. If we set $\left[k\right] := \left\{1, \ldots, k\right\}$ and $\frP := \calP\left(\left[k\right]\right) \setminus \left\{\varnothing, \left[k\right]\right\}$
  where $\calP$ denotes the power set, then we can further write
  \begin{align*}
    E_{y,n} \subseteq \bigcup_{P \in \frP} &\overbrace{\left\{ \vphantom{\underline{X}^{\left(m_2\right)}} P_1 = P = \left[k\right] \setminus P_2, y \textrm{ not visited by walkers in } P_2 \textrm{ after time } n,\right.}^{=: F_{y,n,P}}\\
    &\left.\hphantom{\{} \qquad\qquad\qquad\qquad\qquad\;\;\;\;
    \exists m_1, m_2: \overline{X}^{\left(m_1\right)} = \infty, \underline{X}^{\left(m_2\right)} = -\infty \right\}
  \end{align*}
  and it suffices to show $\Prb{F_{y,n,P}} = 0$ for all $y \in \bbZ, n \in \bbN, P \in \frP$. We are
  now ready to apply \autoref{lem:retRootOrFinManyNew} for every walker in the fixed set $\left[k\right] \setminus P = P_2$
  (equality if $F_{y,n,P}$ occurs). Since the intersection of $F_{y,n,P}$ and the event that any walker
  in $\left[k\right] \setminus P = P_2$ returns to $y$ at a time larger than $n$ is empty, we
  can conclude by \autoref{lem:retRootOrFinManyNew} (as in the proof of \autoref{cor:oneReturnsOrFiniteRange},
  only with the random walk started at time $n + 1$, and with a subset of the walkers) that
  \begin{align*}
    F_{y,n,P} \subseteq &
    \overbrace{\left\{ \vphantom{\overline{X}^{\left(m_1\right)}} P_1 = P, P_2 = \left[k\right] \setminus P, y \textrm{ not visited by walkers in } P_2 \textrm{ after time } n \textrm{ and the walkers}\right.}^{=: G_{y,n,P}}\\
    & \hphantom{\{\;} \left.\vphantom{\overline{X}^{\left(m_1\right)}} \quad
    \textrm{in }P_2\textrm{ only visit finitely many nodes not visited before by any other walker,} \right.\\
    & \hphantom{\{\;} \left.\vphantom{\overline{X}^{\left(m_1\right)}} \qquad\qquad\qquad\qquad\qquad\qquad\qquad\qquad\qquad\;\;\;
    \exists m_1, m_2: \overline{X}^{\left(m_1\right)} = \infty, \underline{X}^{\left(m_2\right)} = -\infty \right\} \\
    & \cup \; \widehat{N} \qquad \textrm{ where } \widehat{N} \textrm{ is a null set}
  \end{align*}
  and it suffices to show $\Prb{G_{y,n,P}} = 0$. Choose $m_4 \in \argmax_{m \in P_1} \overline{X}^{\left(m\right)}$. If $G_{y,n,P}$
  occurs, then we necessarily have $\overline{X}^{\left(m_4\right)} = \infty$. If this was not the case, the walkers in $P_1$ would
  only visit nodes to the left of some fixed integer, and since the walkers in $P_2$ together only
  visit finitely many new nodes (when $G_{y,n,P}$ occurs), this would imply that there is some largest visited integer, a
  contradiction to the fact that $\overline{X}^{\left(m_1\right)} = \infty$ when $G_{y,n,P}$ occurs.

  Now consider the random set of recurrent walkers $P_3 := \left\{ m \in P_1 : \overline{X}^{\left(m\right)} = \infty \right\} \neq \varnothing$
  (since $m_4 \in P_3$ when $G_{y,n,P}$ occurs) and the walker $m_5 \in \argmin_{m \in P_2} \underline{X}^{\left(m\right)}$.
  If $m_5$ has finite range, then it is clear that $m_5$ will meet any walker in $P_3$ infinitely often.
  On the other hand, if $m_5$ has infinite range, then $\overline{X}^{\left(m_5\right)} = \infty$
  since $m_5 \in P_2$. But since the walkers in $P_2$, including $m_5$, only
  visit finitely many nodes not visited before by any other walker when $G_{y,n,P}$ occurs, at least one walker $m_6$ in $P_3$ must be to
  the right of $m_5$ infinitely often in order to ``free the path'' for $m_5$. As the walkers in $P_3$
  are all recurrent, and as $m_5$ only visits nodes to the right of $y$, this implies that $m_5$ meets
  this walker $m_6$ infinitely often (when $G_{y,n,P}$ occurs). Therefore, regardless of whether $m_5$ has finite range or not,
  there is always a recurrent walker $m_6 \in P_3$ which meets $m_5$ infinitely often when $G_{y,n,P}$ occurs.

  We now apply \autoref{lem:labelExchange} \ref{lem:labelExchangeBothRec} to see that $G_{y,n,P}$
  is a null set (which concludes the proof): as done repeatedly in this proof, we find that $G_{y,n,P}$
  is a subset of a countable union of events of the form ``some fixed walker is not recurrent
  (corresponding to $m_5$) and meets some other fixed walker which is recurrent (corresponding to $m_6$)
  infinitely often''. However, any event of this form is a null set by \autoref{lem:labelExchange} \ref{lem:labelExchangeBothRec}.
\end{tproof}

\clearpage
\newpage

\section{Biased Reinforced Random Walk}
\label{sec:biased}

Here, we consider biased random walks with reinforcement. We could prove some results, but there are also still open
questions. This section therefore presents proven results next to open questions and conjectures, and sometimes also
just possible ways to tackle a question, which would have to be pursued further. In order to limit the length of this
thesis, this section will often only provide proof sketches or skip some proof altogether.

\subsection{$\lambda^\ast$-Biased Edge-Reinforced Random Walk}

The $\lambda^\ast$-biased edge-reinforced random walk uses the same linearly reinforced edge weights as the
LERRW, but introduces an additional bias in a certain direction. Here, we consider a biased reinforced
random walk on $\bbZ$, so the bias can be either to the left or to the right. Formally,
define the $\lambda^\ast$-biased edge-reinforced random walk on $\bbZ$ as the sequence $X_n$ of random variables
with
\begin{align*}
  X_0 &= 0 \\
  \Prb{X_{n+1} = x_n + 1 \midd| X_n = x_n, \ldots, X_0 = x_0} &= \frac{\lambda \cdot w\left(n, x_n\right)}{w\left(n, x_n-1\right) + \lambda \cdot w\left(n, x_n\right)} \\
  w\left(n, z\right) &= \underbrace{w\left(0, z\right)}_{= 1 \textrm{ here}} + \sum_{i=0}^{n-1} \mathbbm{1}_{X_i = z, X_{i + 1} = z + 1} + \mathbbm{1}_{X_i = z + 1, X_{i + 1} = z}
\end{align*}

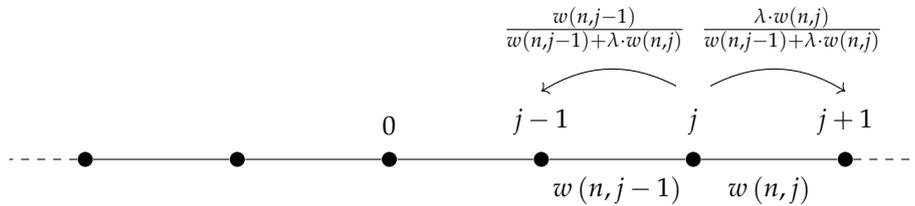
\begin{figure}[H]
  \begin{center}
    \begin{tikzpicture}
      \node[circle,fill=black,inner sep=0.7mm] (N0) at (0, 0) {};
      \node[circle,fill=black,inner sep=0.7mm] (N1) at (2, 0) {};
      \node[circle,fill=black,inner sep=0.7mm] (N2) at (4, 0) {};
      \node[circle,fill=black,inner sep=0.7mm] (N3) at (6, 0) {};
      \node[circle,fill=black,inner sep=0.7mm] (Nm1) at (-2, 0) {};
      \node[circle,fill=black,inner sep=0.7mm] (Nm2) at (-4, 0) {};
      \node[above=2mm] at (N0) {$0$};
      \node[above=2mm] at (N1) {$j-1$};
      \node[above=2mm] at (N2) {$j$};
      \node[above=2mm] at (N3) {$j+1$};
      \draw (Nm2) -- (Nm1) -- (N0) -- (N1) -- node[below=1mm] {$w\left(n,j-1\right)$} (N2) -- node[below=1mm] {$w\left(n,j\right)$} (N3);
      \draw[dashed] (N3) -- (7, 0);
      \draw[dashed] (Nm2) -- (-5, 0);
      \node (aN2l) at (3.9, 0.9) {};
      \node (aN2r) at (4.1, 0.9) {};
      \draw[->] (aN2l) edge[bend right] node[above=1mm] {$\frac{w\left(n, j-1\right)}{w\left(n, j-1\right) + \lambda \cdot w\left(n, j\right)} \quad$} (2, 0.9);
      \draw[->] (aN2r) edge[bend left] node[above=1mm] {$\quad \frac{\lambda \cdot w\left(n, j\right)}{w\left(n, j-1\right) + \lambda \cdot w\left(n, j\right)}$} (6, 0.9);
    \end{tikzpicture}
    \caption{Edge weights and transition probabilities for the biased walk on $\bbZ$ at time $n$}
    \label{fig:biased_errw_z}
  \end{center}
\end{figure}

In terms of \autoref{def:errw}, we take $W\left(n\right) = 1 + n$, so the edge weights
are being linearly reinforced, but then we introduce an additional bias by multiplying the right edge weight
with the parameter $\lambda$.

Now, consider $\lambda \in \bbQ_{>0}$ (the result can later be extended to all $\lambda \in \bbR_{>0}$)
and a single node, say $0$, the root. Initially, the weights on both adjacent edges are $1$. Whenever
the random walk leaves $0$ along one edge, it can only return to $0$ by
the same edge. The weight of this edge will then have increased by $2$. If $\lambda = 1$, then this
process is equivalent to drawing balls from an urn which initially contains $1$ black and $1$ white ball
and where a ball is replaced with $3$ balls of the same color as the ball which was drawn, i.e.~the urn
will contain two more balls of the respective color than before (white balls
corresponding to the left edge, black balls to the right edge). In general, the process is biased, and
the number of balls which are put into the urn has to be adapted:

\begin{lemma}
  \label{lem:urn}
  Let $\lambda = \frac{p}{q}$.
  Consider a node $z \in \bbZ$ with initial edge weights $w\left(0, z-1\right) = l_0$ and $w\left(0, z\right) = r_0$.
  If the $\lambda^\ast$-biased edge-reinforced random walk is started in $z$,
  then the sequence of left turns and right turns of the random walk at $z$ has the same distribution
  as the sequence of white and black balls drawn from the following urn:
  \begin{itemize}
    \item Initially, the urn contains $a_0 = q \cdot l_0$ white and $b_0 = p \cdot r_0$ black balls
    \item If a white ball is drawn, it is replaced with $2q + 1$ white balls (i.e.~$2q$ more white balls than before)
    \item If a black ball is drawn, it is replaced with $2p + 1$ black balls (i.e.~$2p$ more black balls than before)
  \end{itemize}
  The urn and the edge weights can be coupled such that, writing $\tau_n$ for the time
  at which $z$ is visited for the $n$-th time, $w\left(\tau_n, z-1\right) = \frac{a_n}{q}$
  and $w\left(\tau_n, z\right) = \frac{b_n}{p}$.
\end{lemma}

\begin{sproof}
  Note that, after the $n$-th draw and the $n$-th visit to the node $z$, we have
  with the coupling given in \autoref{lem:urn}:
  \begin{align*}
    \Prb{X_{\tau_n+1} = z + 1 \midd| X_{\tau_n} = z, \ldots, X_0 = x_0} &= \frac{\lambda \cdot w\left(\tau_n, z\right)}{w\left(\tau_n, z-1\right) + \lambda \cdot w\left(\tau_n, z\right)}
    = \frac{\frac{p}{q} \cdot \frac{b_n}{p}}{\frac{a_n}{q} + \frac{p}{q} \cdot \frac{b_n}{p}} \\
    &= \frac{b_n}{q} \cdot \frac{q}{a_n + b_n}
    = \frac{b_n}{a_n + b_n} \qedhere
  \end{align*}
\end{sproof}

\begin{lemma}
  \label{lem:asconv}
  Consider the $\lambda^\ast$-biased edge-reinforced random walk at node $z \in \bbZ$ with $\lambda \in \bbQ_{>0}$.
  If $z$ is visited infinitely often, then
  \begin{enumerate}[(i)]
    \item If $\lambda > 1$, then $\frac{\lambda w\left(n, z\right)}{w\left(n, z - 1\right) + \lambda w\left(n, z\right)}$
    converges almost surely to $1$
    \item If $\lambda < 1$, then $\frac{w\left(n, z - 1\right)}{w\left(n, z - 1\right) + \lambda w\left(n, z\right)}$
    converges almost surely to $1$ \qedhere
  \end{enumerate}
\end{lemma}

\begin{sproof}
  Use \autoref{lem:urn} and \autoref{thm:urnconv} with $\mu\left(\left\{p\right\}\right) = 1$ and $\nu\left(\left\{q\right\}\right) = 1$
  (or vice versa if $p < q$).
\end{sproof}

\begin{definition}
  The $\lambda^\ast$-biased edge-reinforced random walk is said to be
  \begin{itemize}
    \item \textbf{transient} if it a.s.~visits every node only finitely often
    \item \textbf{recurrent} otherwise, i.e.~if there is at least one node which is visited infinitely often \qedhere
  \end{itemize}
\end{definition}

\begin{lemma}
  \label{lem:rec}
  If the $\lambda^\ast$-biased edge-reinforced random walk is recurrent, then all nodes are visited infinitely often almost surely.
\end{lemma}

\begin{tproof}
  If the random walk is recurrent, then there is at least one node which is visited infinitely often, say $z$. To show that all nodes are visited
  infinitely often, it suffices to show that $z$ being visited infinitely often implies that both neighbors of $z$ are visited infinitely often
  (one can then continue by induction). Assume
  for a contradiction that one of the neighbors, say $y$, is visited only finitely often, and let $t$
  be the time of the last visit to $y$ (at time $t + 1$, the random walk then necessarily is at $z$).
  We can assume w.l.o.g.~$y = z + 1$ (in the other case, reflect $\bbZ$ at $z$, set a new $\lambda' = \frac{1}{\lambda}$ and continue as below).
  
  Let $\tau_1 = t + 1, \tau_2, \tau_3, \ldots$ be the times at which $z$ is visited after the last visit
  to $y$. Then, at each time $\tau_n$, we have
  \begin{align*}
    \Prb{X_{\tau_n+1} = z + 1 = y \midd| X_{\tau_n} = z, \ldots, X_0 = x_0} &= \frac{\lambda \cdot w\left(\tau_n, z\right)}{w\left(\tau_n, z-1\right) + \lambda \cdot w\left(\tau_n, z\right)}\\
    &= \frac{\lambda \cdot \overbrace{w\left(t+1, z\right)}^{w_r}}{\lambda \cdot w\left(t+1, z\right) + \underbrace{w\left(t+1, z-1\right)}_{w_l} + n-1} \\
    \implies \Prb{\forall n: X_{\tau_n+1} \neq y} = \prod_{n \geq 1} \frac{w_l + n - 1}{\lambda w_r + w_l + n - 1}
    &\overset{\circledast}{\leq} \exp\left(\sum_{n \geq 1} \left(\frac{w_l + n - 1}{\lambda w_r + w_l + n - 1} - 1\right)\right) \\
    &\hphantom{\;=\;} = \exp\left(-\sum_{n \geq 1} \frac{\lambda w_r}{\lambda w_r + w_l + n - 1}\right) = 0 \\
    \textrm{where }\circledast\textrm{ holds since } \forall x \in \bbR: x &\leq \exp\left(x - 1\right)
  \end{align*}
  Hence, the event $y$ being visited only finitely often has probability $0$.
\end{tproof}

Based on \autoref{lem:asconv}, it seems reasonable to conjecture that the $\lambda^\ast$-biased edge-reinforced random walk is transient
whenever $\lambda \neq 1$, since at every node, the probability to go in the direction of bias would converge to $1$ if the node was
visited infinitely often. This intuitively seems to be a contradiction: if the probability to go in one direction goes to $1$, then every
node should be visited only finitely often. However, we did not manage to prove the following conjecture.

\begin{conjecture}
  \label{conj:trans}
  The $\lambda^\ast$-biased edge-reinforced random walk is transient for $\lambda \neq 1$.
\end{conjecture}

Note that for $\lambda = 1$, the walk is recurrent since this case corresponds to the standard LERRW
on $\bbZ$.

\subsubsection{Some Simulations}

Below, the result of some simulations is presented. The $\lambda^\ast$-biased edge-reinforced random walk
was simulated in 100 simulations for 1 000 000 steps each, in 1 000 simulations for 100 000 steps each
and in 10 000 simulations for 10 000 steps each.

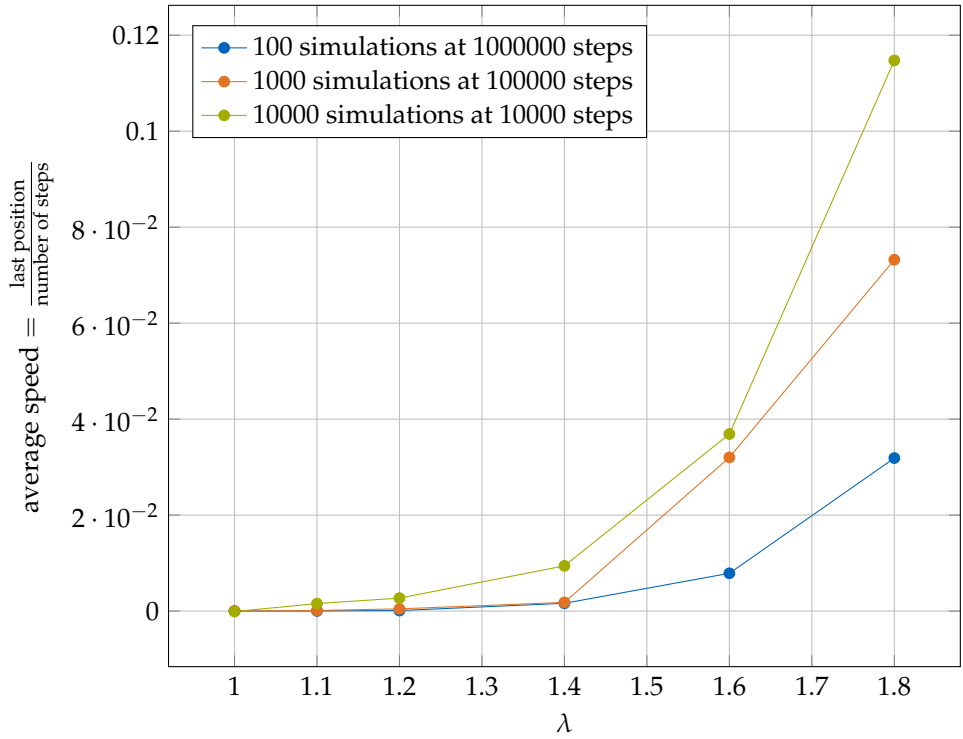
\begin{figure}[H]
  \begin{center}
    \begin{tikzpicture}
  \begin{axis}[
    xlabel={$\lambda$},
    ylabel={$\textrm{average speed} = \frac{\textrm{last position}}{\textrm{number of steps}}$},
    width=12cm,
    grid=major,
    legend pos=north west
  ]
    \addplot[color=tumBlue,mark=*] coordinates {
      (1.00, 0.00000) (1.10, 0.00002) (1.20, 0.00011) (1.40, 0.00162) (1.60, 0.00787) (1.80, 0.03188)
    };
  \addlegendentry{$100$ simulations at $1000000$ steps};
    \addplot[color=tumOrange,mark=*] coordinates {
      (1.00, 0.00000) (1.10, 0.00012) (1.20, 0.00046) (1.40, 0.00180) (1.60, 0.03202) (1.80, 0.07321)
    };
  \addlegendentry{$1000$ simulations at $100000$ steps};
    \addplot[color=tumGreen,mark=*] coordinates {
      (1.00, -0.00010) (1.10, 0.00156) (1.20, 0.00270) (1.40, 0.00942) (1.60, 0.03689) (1.80, 0.11474)
    };
  \addlegendentry{$10000$ simulations at $10000$ steps};
  \end{axis}
\end{tikzpicture}

    \caption{Average speed of the random walk}
    \label{fig:avg_speed}
  \end{center}
\end{figure}

As can be seen in \autoref{fig:avg_speed}, the simulated speed strongly depends on the number of steps
for which the random walk was simulated. In this case, the speed decreased with an increasing number of
simulated steps. This could either be an indication that the speed has not yet converged and longer simulations
would be necessary or that the random walk is actually always recurrent and hence the speed is actually
zero (however, it seems intuitively very likely that it will be transient for large values of $\lambda$).
For $\lambda = 1.1$, the speed was almost indistinguishable from $0$. To see if the walk is really recurrent, it
is however also of interest when the root was visited for the last time, since a speed of $0$ does not
necessarily imply recurrence.

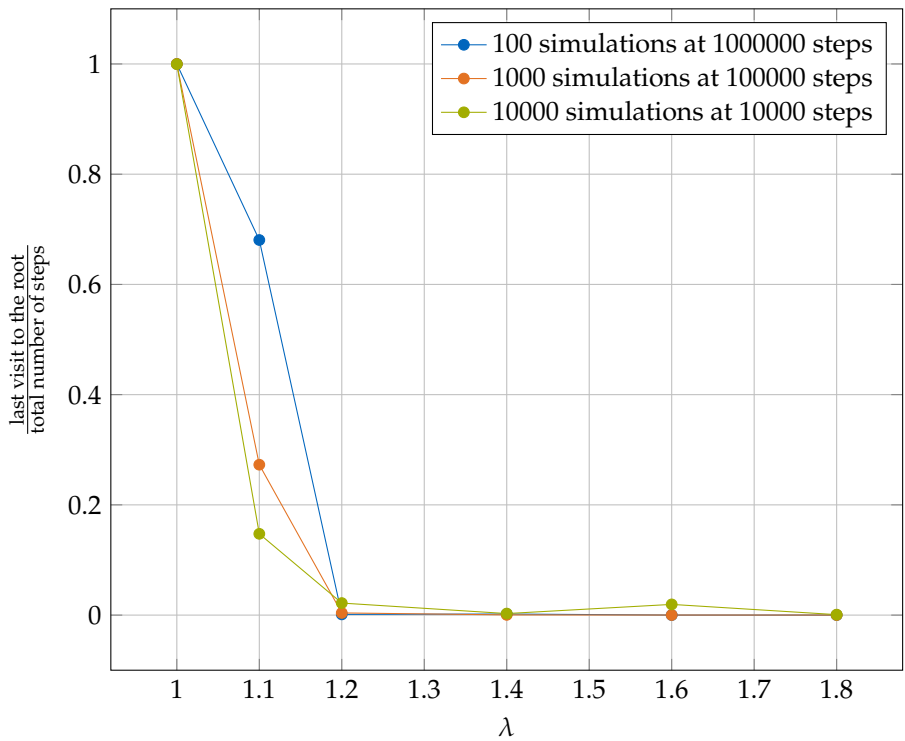
\begin{figure}[H]
  \begin{center}
    \begin{tikzpicture}
  \begin{axis}[
    xlabel={$\lambda$},
    ylabel={$\frac{\textrm{last visit to the root}}{\textrm{total number of steps}}$},
    width=12cm,
    grid=major
  ]
    \addplot[color=tumBlue,mark=*] coordinates {
      (1.00, 0.99994) (1.10, 0.68057) (1.20, 0.00127) (1.40, 0.00199) (1.60, 0.00002) (1.80, 0.00001)
    };
  \addlegendentry{$100$ simulations at $1000000$ steps};
    \addplot[color=tumOrange,mark=*] coordinates {
      (1.00, 0.99999) (1.10, 0.27286) (1.20, 0.00395) (1.40, 0.00038) (1.60, 0.00012) (1.80, 0.00007)
    };
  \addlegendentry{$1000$ simulations at $100000$ steps};
    \addplot[color=tumGreen,mark=*] coordinates {
      (1.00, 0.99972) (1.10, 0.14747) (1.20, 0.02183) (1.40, 0.00279) (1.60, 0.01942) (1.80, 0.00083)
    };
  \addlegendentry{$10000$ simulations at $10000$ steps};
  \end{axis}
\end{tikzpicture}

    \caption{Average last visit to the root by the random walk}
    \label{fig:avg_lastvis}
  \end{center}
\end{figure}

Longer or more simulations are necessary to get reliable results (note for example that in \autoref{fig:avg_lastvis},
the green curve suddenly goes up for $\lambda = 1.6$ which is likely not a property of the model, but
a case of variation in the simulation). However, it is already visible that for $\lambda = 1.1$, the last
visit to the root occurred on average relatively late, i.e.~at a noticeable fraction of the total number
of simulated steps. Indeed, the last visit to the root seemed to occur later (compared to the total number
of steps) when simulating the walk for longer, which we would not have expected. Of course, this could
also just be a case of variation in the simulation.

The simulations do unfortunately not present any evidence which would support \autoref{conj:trans},
but to argue about properties like recurrence and transience, which only make sense on an infinite time
horizon, is very difficult with the help of only finite simulations.

\subsubsection{Stochastic Approximation}

To better understand the $\lambda^\ast$-biased edge-reinforced random walk on $\bbZ$, it might be helpful
to approximate it by reinforced walks on finite circular graphs. We superficially present a basic idea how
this could be done. However, it is completely unclear if the analysis of the walk on a finite circular graph
will provide new insights for the walk on $\bbZ$. Since the easiest circular graph is the triangle, we
consider the $\lambda^\ast$-biased edge-reinforced random walk on the triangle, defined as follows:
\begin{align*}
  r\left(x\right) &= \left(x \mmod 3\right) + 1 \textrm{ the right neighbor of }x \in \left\{1, 2, 3\right\} \\
  l\left(x\right) &= \left(\left(x + 1\right) \mmod 3\right) + 1 \textrm{ the left neighbor of }x \\
  X_0 &= 1 \\
  \Prb{X_{n+1} = r\left(x_n\right) \midd| X_n = x_n, \ldots, X_0 = x_0} &= \frac{\lambda \cdot w\left(n,x_n\right)}{\lambda \cdot w\left(n,x_n\right) + w\left(n,l\left(x_n\right)\right)} \\
  w\left(n,i\right) &= \underbrace{w\left(0,i\right)}_{= 1 \textrm{ here}} + \sum_{k=0}^{n-1} \mathbbm{1}_{X_k = i, X_{k + 1} = r\left(i\right)} + \mathbbm{1}_{X_k = r\left(i\right), X_{k + 1} = i} \\
  c_{n,i} &= \frac{w\left(n,i\right)}{\sum_{j=1}^3 w\left(n,j\right)}
\end{align*}
For simplicity, we will assume $\lambda > 1$.

We call $\mathbf{c}_n$ the vector of normalized edge weights at time $n$, which is located on the
unit simplex and has entries $c_{n,i}$ where $i=1,2,3$.
If we assume that the normalized edge weights were fixed to some vector $c$ (on the unit simplex),
then we get a Markov chain with the following stationary distribution $\pi_c$ (only the value at the node
$1$ is shown, the others follow from the symmetry of the model):
\begin{align*}
  \pi_c\left(1\right) &=
  \left(
    \lambda^3 c_1c_2c_3 +
    \lambda^2 \left(c_1^2c_3 + c_2c_3^2\right) +
    \lambda \left(c_1^2c_2 + c_1c_3^2\right) +
    c_1c_2c_3
  \right) Z^{-1}
\end{align*}
where $Z$ is the appropriate normalizing constant.

\begin{figure}[H]
  \begin{center}
    \begin{tikzpicture}[scale=1.2]
      \node[circle,fill=black,inner sep=0.7mm] (N1) at (0,0) {};
			\node[circle,fill=black,inner sep=0.7mm] (N2) at (4,0) {};
			\node[circle,fill=black,inner sep=0.7mm] (N3) at (2,3.46) {};
			\node[below left=2mm] at (N1) {$1$};
			\node[below right=2mm] at (N2) {$2$};
			\node[above=3mm] at (N3) {$3$};
			\draw (N1) -- (N2) -- (N3) -- (N1);
			\draw ([shift={(0.3,-0.3)}]N1.center) edge[bend right,-{Latex[length=2mm,width=2mm]}] node[below] {$\frac{\lambda c_{n,1}}{\lambda c_{n,1}+c_{n,3}}$} ([shift={(-0.3,-0.3)}]N2.center);
			\draw ([shift={(0.11,0.41)}]N2.center) edge[bend right,-{Latex[length=2mm,width=2mm]}] node[right=1mm] {$\frac{\lambda c_{n,2}}{\lambda c_{n,2}+c_{n,1}}$} ([shift={(0.41,-0.11)}]N3.center);
			\draw ([shift={(-0.41,-0.11)}]N3.center) edge[bend right,-{Latex[length=2mm,width=2mm]}] node[left=1mm] {$\frac{\lambda c_{n,3}}{\lambda c_{n,3}+c_{n,2}}$} ([shift={(-0.11,0.41)}]N1.center);
			\draw ([shift={(0.51,0.28)}]N1.center) edge[bend right,-{Latex[length=2mm,width=2mm]}] node[right=1mm] {$\frac{c_{n,3}}{\lambda c_{n,1}+c_{n,3}}$} ([shift={(0.01,-0.58)}]N3.center);
    \end{tikzpicture}
    \caption{Transition probabilities of the $\lambda^\ast$-biased ERRW on the triangle at time $n$}
    \label{fig:biased_errw_triangle}
		(two arrows are missing for clarity)
  \end{center}
\end{figure}
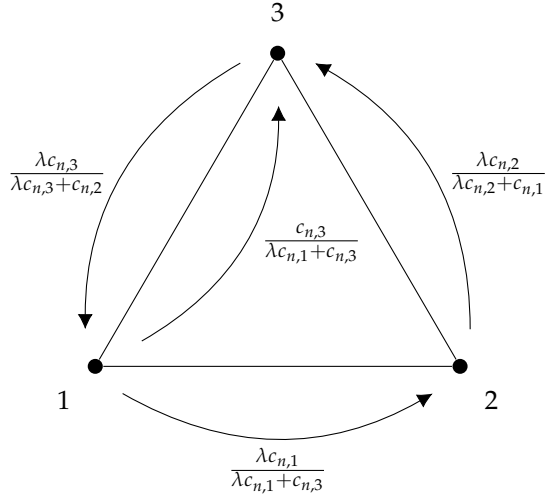

We further define the stationary distribution on
the edges $\pi_c^{\textrm{edge}}$ simply by looking at the time spent on each edge when the chain is
run starting from the stationary distribution $\pi$. For example,
\begin{align*}
  \pi_c^{\textrm{edge}}\left(1\right) &=
  \frac{\lambda c_1}{\lambda c_1+c_3} \cdot \pi_c\left(1\right) +
  \frac{c_1}{\lambda c_2+c_1} \cdot \pi_c\left(2\right)
\end{align*}
We now approximate the evolution of the time-dependent
vector $\mathbf{c}_n$ of the edge weights:
\begin{align*}
  1 \ll k \ll n \implies \left(n + k\right) \mathbf{c}_{n+k} &\approx n \mathbf{c}_n + k \pi_{\mathbf{c}_n}^{\textrm{edge}} \\
  \iff \mathbf{c}_{n+k} - \mathbf{c}_{n} &\approx \frac{k}{n+k} \left(\pi_{\mathbf{c}_n}^{\textrm{edge}} - \mathbf{c}_{n}\right) \\
  \textrm{therefore approximate with } \qquad \frac{\textrm{d}}{\textrm{d}t}\widetilde{\mathbf{c}}\left(t\right) &= \frac{1}{t}\left(\pi_{\widetilde{\mathbf{c}}\left(t\right)}^{\textrm{edge}} - \widetilde{\mathbf{c}}\left(t\right)\right) \\
  \textrm{exponential time change to get } \qquad \frac{\textrm{d}\mathbf{c}}{\textrm{d}t} &= \pi_{\mathbf{c}}^{\textrm{edge}} - \mathbf{c} \qquad \textrm{ with }
  \widetilde{\mathbf{c}}\left(t\right) = \mathbf{c}\left(\ln\left(t\right)\right)
\end{align*}

This is called a stochastic approximation: we approximate the random evolution of the vector of normalized
edge weights by a differential equation. This is possible because the changes in the vector of normalized
edge weights get ever smaller as time increases, and the randomness gets less noticeable by virtue of the
law of large numbers. To argue formally, many additional steps would be necessary, but this approximation
can already give a good intuition on what is happening.

\begin{figure}[H]
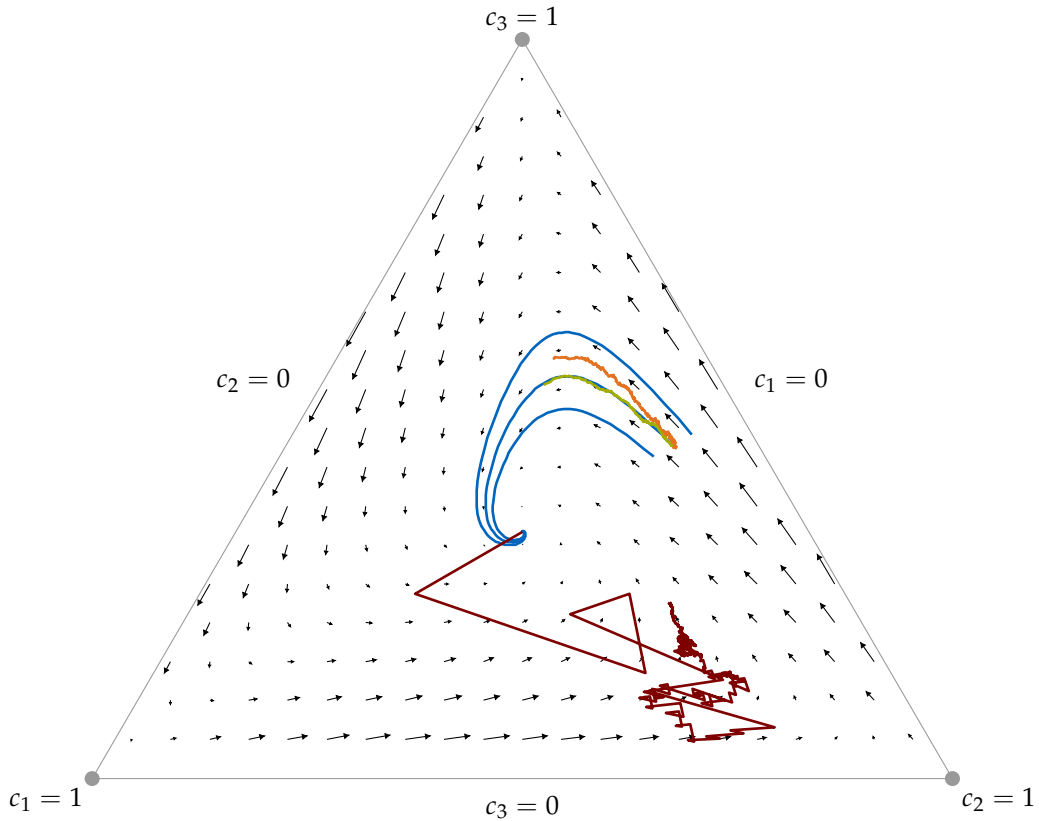

  \begin{center}

    \caption{The vector field of the differential equation on the unit simplex}
    \label{fig:diff_vfield_simplex}
  \end{center}
	In \textcolor{sectionblue}{blue}, three segments of solutions to the differential equation are shown.
	In \textcolor{linkred}{red}, \textcolor{tumOrange}{orange} and \textcolor{tumGreen}{green}, three actual trajectories obtained
  by simulating the $\lambda^\ast$-biased edge-reinforced random walk are shown. The \textcolor{linkred}{red}
  trajectory results from setting the initial weights to $1$ and running the simulation for 10 000
  steps. For the \textcolor{tumOrange}{orange} trajectory, the initial weights were $100, 450, 450$
  and the simulation was run for 100 000 steps; for the \textcolor{tumGreen}{green} one, they
  were $500, 2250, 2250$ and the simulation ran for 500 000 steps.
\end{figure}

Without any sort of formal proof by simply looking at \autoref{fig:diff_vfield_simplex}, we conclude
that the normalized edge weights will eventually converge to the uniform distribution. This probably
also holds for circles with more than three nodes, showing that the long-run behavior of the biased
reinforced random walk on circles is fundamentally different from the long-run behavior of the
walk on the integers (where at every vertex the quotient of the right edge weight divided by the left
edge weight converges to $1$).

This already gives an indication that approximating the wlak on $\bbZ$ with that on a finite circular graph
could be very hard. Indeed, \autoref{lem:asconv} shows that in the long run, the ratio of the adjacent edge weights
at every node behaves fundamentally different for the two models. It could however be possible to analyze the
reinforced walk on finite circular graphs after a finite amount of time, increasing with the size of the graph,
to make the transition from the circular graphs to $\bbZ$.

\subsection{$\lambda^+$-Biased Edge-Reinforced Random Walk}

Similarly to the multiplicative bias studied above, we can also consider an additive bias. Define the $\lambda^+$-biased
edge-reinforced random walk on $\bbZ$ as the sequence $X_n$ of random variables with
\begin{align*}
  X_0 &= 0 \\
  \Prb{X_{n+1} = x_n + 1 \midd| X_n = x_n, \ldots, X_0 = x_0} &= \frac{\lambda + w\left(n, x_n\right)}{\lambda + w\left(n, x_n\right) + w\left(n, x_n - 1\right)} \\
  w\left(n, z\right) &= \underbrace{w\left(0, z\right)}_{= 1 \textrm{ here}} + \sum_{i=0}^{n-1} \mathbbm{1}_{X_i = z, X_{i + 1} = z + 1} + \mathbbm{1}_{X_i = z + 1, X_{i + 1} = z}
\end{align*}
(where we require $\lambda \geq 0$). The problem with the definition of the $\lambda^\ast$--biased edge-reinforced
random walk was that it could no longer be represented as a mixture of Markov chains since the probability
of the occurrence of a certain edge sequence depended on the order in which the edges appeared in the sequence.
This is not the case with the new definition given above. In other words, when the edge weights are
represented by urns at every node, then the sequence of draws is now exchangeable, while this was not
the case before.

Following \cite[Lemma 1 and Lemma 2]{errwpemantle}, the $\lambda^+$-biased edge-reinforced random walk is
equivalent to a mixture of MCs (or a RWRE). For shorter notation, we will adapt the
notation of \autoref{def:rwre} as follows: we call $\Prb{X_{n+1} = x+1 \midd| X_n = x} = \ensuremath{P_{\mathbf{c}}}\left(x, x+1\right) =: \omega_x$.
Instead of looking at the distribution of the random conductances, we look directly at the distribution of the
$\omega_x$.

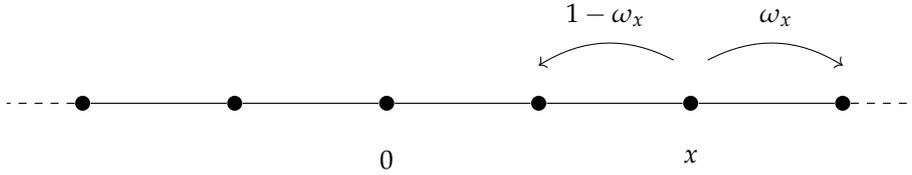
\begin{figure}[H]
  \begin{center}
    \begin{tikzpicture}
      \node[circle,fill=black,inner sep=0.7mm] (N0) at (0, 0) {};
      \node[circle,fill=black,inner sep=0.7mm] (N1) at (2, 0) {};
      \node[circle,fill=black,inner sep=0.7mm] (N2) at (4, 0) {};
      \node[circle,fill=black,inner sep=0.7mm] (N3) at (6, 0) {};
      \node[circle,fill=black,inner sep=0.7mm] (Nm1) at (-2, 0) {};
      \node[circle,fill=black,inner sep=0.7mm] (Nm2) at (-4, 0) {};
      \node[below=5mm] at (N0) {$0$};
      \node[below=5mm] at (N2) {$x$};
      \draw (Nm2) -- (Nm1) -- (N0) -- (N1) -- (N2) -- (N3);
      \draw[dashed] (N3) -- (7, 0);
      \draw[dashed] (Nm2) -- (-5, 0);
      \node (aN2l) at (3.9, 0.5) {};
      \node (aN2r) at (4.1, 0.5) {};
      \draw[->] (aN2l) edge[bend right] node[above=1mm] {$1 - \omega_x$} (2, 0.5);
      \draw[->] (aN2r) edge[bend left] node[above=1mm] {$\omega_x$} (6, 0.5);
    \end{tikzpicture}
    \caption{A random walk in a random environment on $\bbZ$}
    \label{fig:errw_z}
  \end{center}
\end{figure}

The mixture of MCs equivalent to the $\lambda^+$-biased walk is constructed by placing independent urns at every node which are coupled
with which edges the random walk takes. Initially, we start with $1$ black and $1 + \lambda$ (we can
also have half a ball or any positive real number of balls) white balls, the black balls representing
the edge to the left of the node, the white balls representing the edge to the right. Now, if a ball
is drawn, we take the edge corresponding to the color and in addition to the drawn ball add two more balls
of the same color (when the random walk returns to the node, the edge weight will have increased by two).
Now, for any nodes to the right of $0$, the first time the random walk gets there, there will already
be $2$ black balls in the urn, since the edge to the left was already traversed once. For nodes to the
left of $0$, there will already be $2 + \lambda$ white balls. So we have the following urns:
\begin{itemize}
  \item at nodes $>0$, the urn initially contains $2$ black and $1 + \lambda$ white balls
  \item at $0$, the urn initially contains $1$ black and $1 + \lambda$ white balls
  \item at nodes $<0$, the urn initially contains $1$ black and $2 + \lambda$ white balls
  \item whenever a ball of one color is drawn, it is put back together with two more balls of the same color 
\end{itemize}

With \cite[Lemma 1 and Lemma 2]{errwpemantle}, we get that the reinforced random walk with initial bias
is equivalent to the mixture of Markov Chains where the $\omega_x$ are independent with the following distributions
(also compare with \cite[Section 5]{lerrwbsc}):
\begin{align*}
  \textrm{for }x > 0\textrm{: }& \qquad \omega_x \sim \Betad{\frac{1 + \lambda}{2}}{\frac{2}{2}} \\
  \textrm{for }x = 0\textrm{: }& \qquad \omega_x \sim \Betad{\frac{1 + \lambda}{2}}{\frac{1}{2}} \\
  \textrm{for }x < 0\textrm{: }& \qquad \omega_x \sim \Betad{\frac{2 + \lambda}{2}}{\frac{1}{2}}
\end{align*}
where $\BETAD$ is the beta distribution. We use the following:

\begin{lemma}
  \label{lem:beta_distr}
  Let $A \sim \Betad{\alpha}{\beta}$ with $\alpha,\beta > 0$ and $t \in \bbR$. Then:
  \begin{align*}
    \Ex{\left(\frac{1 - A}{A}\right)^t} &= \begin{cases}
      \frac{\Gamma\left(\alpha - t\right)\Gamma\left(\beta + t\right)}{\Gamma\left(\alpha\right)\Gamma\left(\beta\right)}
      & \textrm{ if } -\beta < t < \alpha \\
      \infty & \textrm{ otherwise}
    \end{cases}
    \qquad \Ex{\frac{1 - A}{A}} = \begin{cases}
      \frac{\beta}{\alpha - 1}
      & \textrm{ if } \alpha > 1 \\
      \infty & \textrm{ otherwise}
    \end{cases} \\
    \Ex{\left(\frac{A}{1 - A}\right)^t} &= \begin{cases}
      \frac{\Gamma\left(\alpha + t\right)\Gamma\left(\beta - t\right)}{\Gamma\left(\alpha\right)\Gamma\left(\beta\right)}
      & \textrm{ if } -\alpha < t < \beta \\
      \infty & \textrm{ otherwise}
    \end{cases}
    \qquad \Ex{\frac{A}{1 - A}} = \begin{cases}
      \frac{\alpha}{\beta - 1}
      & \textrm{ if } \beta > 1 \\
      \infty & \textrm{ otherwise}
    \end{cases} \qedhere
  \end{align*}
\end{lemma}

\begin{tproof}
  We have:
  \begin{align*}
    \Ex{\left(\frac{1 - A}{A}\right)^t}
    &= \frac{\Gamma\left(\alpha + \beta\right)}{\Gamma\left(\alpha\right)\Gamma\left(\beta\right)}
    \int_0^1 \frac{\left(1 - x\right)^t}{x^t} \cdot x^{\alpha - 1} \left(1 - x\right)^{\beta - 1} \dx{x} \\
    &= \frac{\Gamma\left(\alpha + \beta\right)}{\Gamma\left(\alpha\right)\Gamma\left(\beta\right)}
    \int_0^1 x^{\alpha - t - 1} \left(1 - x\right)^{\beta + t - 1} \dx{x} \\
    &= \begin{cases}
      \frac{\Gamma\left(\alpha + \beta\right)}{\Gamma\left(\alpha\right)\Gamma\left(\beta\right)} \cdot
      \frac{\Gamma\left(\alpha - t\right)\Gamma\left(\beta + t\right)}{\Gamma\left(\alpha + \beta\right)}
      = \frac{\Gamma\left(\alpha - t\right)\Gamma\left(\beta + t\right)}{\Gamma\left(\alpha\right)\Gamma\left(\beta\right)}
      & \textrm{ if } -\beta < t < \alpha \\
      \infty & \textrm{ otherwise}
    \end{cases}
  \end{align*}
  Hence,
  \begin{align*}
    \Ex{\frac{1 - A}{A}}
    &= \begin{cases}
      \frac{\Gamma\left(\alpha - 1\right)\Gamma\left(\beta + 1\right)}{\Gamma\left(\alpha\right)\Gamma\left(\beta\right)}
      = \frac{\beta}{\alpha - 1}
      & \textrm{ if } -\beta < 1 < \alpha \\
      \infty & \textrm{ otherwise}
    \end{cases}
  \end{align*}
  The other statements follow from the fact that $1 - A \sim \Betad{\beta}{\alpha}$.
\end{tproof}

We can use \cite[Theorem 1.7]{rwresolomon} to analyze transience and recurrence. It only applies to the situation where the $\omega_x$
are iid. We have only iid variables on the positive half-line as well as iid variables on the negative
half-line, but for a moment we just assume that all $\omega_x$ are distributed as the variables on the positive half-line.
The $\sigma$ from \cite[Theorem 1.7]{rwresolomon} is $\frac{1-\omega_x}{\omega_x}$ in our case. We have,
with $C_{\lambda} > 0$ a constant depending on $\lambda$:
\begin{align*}
  \Ex{\ln\left(\frac{1-\omega_x}{\omega_x}\right)}
  = C_{\lambda} \int_0^1 \ln\left(\frac{1-x}{x}\right) x^{\frac{\lambda-1}{2}} \dx{x}
  \begin{cases}
    > 0 & \textrm{ if } 0 \leq \lambda < 1 \\
    = 0 & \textrm{ if } \lambda = 1 \\
    < 0 & \textrm{ if } \lambda > 1
  \end{cases}
\end{align*}
Unfortunately, this integral cannot be evaluated easily and in a nice way (in the general
case; for $\lambda = 1$ there is a nice representation and probably also for $\lambda \in \bbQ$).
However, WolframAlpha shows how the integral behaves. To formally conclude here, a bit more work
in analyzing the integral would be needed.
\begin{center}
  \begin{tabular}{ccc}
    $\lambda < 1$ & $\lambda = 1$ & $\lambda > 1$ \\
    \href{https://www.wolframalpha.com/input/?i=integral+from+0+to+1+of+ln%28%281-x%29%2Fx%29+*+x%5E%28%280.9-1%29%2F2%29+dx}{https://bit.ly/3yrjGCG} &
    \href{https://www.wolframalpha.com/input/?i=integral+from+0+to+1+of+ln%28%281-x%29%2Fx%29+*+x%5E%28%281-1%29%2F2%29+dx}{https://bit.ly/3F00yy3} &
    \href{https://www.wolframalpha.com/input/?i=integral+from+0+to+1+of+ln%28%281-x%29%2Fx%29+*+x%5E%28%281.1-1%29%2F2%29+dx}{https://bit.ly/3IK0jJQ}
  \end{tabular}
\end{center}

Hence, \cite[Theorem 1.7]{rwresolomon} implies that the random walk in the environment where all $\omega_x$
are iid is transient to the right when $\lambda > 1$, recurrent when $\lambda = 1$, and transient to
the left when $\lambda < 1$. But now, we are actually considering the random environment for the $\lambda^+$-biased
edge-reinforced random walk. The $\omega_x$ with $x \leq 0$ all stochastically dominate the $\omega_x$
with $x > 0$ (the probability to go to the right is higher). Hence, for the case $\lambda > 1$, it
is immediately (intuitively) clear that the random walk will still be transient to the right for this
random environment, and still recurrent for $\lambda = 1$. For the case $\lambda < 1$ we see that
the $\lambda^+$-biased edge-reinforced random walk always returns to the root from the positive half-line.
Now, we already know (e.g.~by \cite[Section 5]{lerrwbsc}) that the edge-reinforced walk without initial
bias is recurrent. The initial bias makes it only more likely to go to the right, hence the walk with
initial bias will also be recurrent on the negative half-line for $\lambda < 1$ and therefore recurrent
as a whole. We therefore get, even though a formal proof would still be needed:

\begin{theorem}
  \label{thm:add_bias_rectrans}
  The edge-reinforced random walk on $\bbZ$ with initial bias $\lambda$ is recurrent for $0 \leq \lambda \leq 1$
  and transient for $\lambda > 1$.
\end{theorem}

Next, we can look at the speed in the transient regime. Since the walk is only transient to the right
if it is transient at all, it suffices to look at the distributions of $\omega_x$ for $x > 0$ (the
finite number of steps spent to the left of the root $0$ is irrelevant for the speed). We have by \autoref{lem:beta_distr}:
\begin{align*}
  \Ex{\frac{1 - \omega_x}{\omega_x}} &=
  \begin{cases}
    \infty & \textrm{ if } \lambda \leq 1 \\
    \frac{2}{\lambda - 1} & \textrm{ if } \lambda > 1
  \end{cases}
\end{align*}
Clearly, $\frac{2}{\lambda - 1}$ is (strictly) decreasing in $\lambda$ for $\lambda > 1$, and hence,
we have $\Ex{\frac{1 - \omega_x}{\omega_x}} < 1 \iff \lambda > 3$. We now want to use \cite[Theorem 1.16]{rwresolomon}.
To get a statement of the type ``positive speed $\iff$ some condition'', we therefore also need to look at
the following, where we use again \autoref{lem:beta_distr}:
\begin{align*}
  \Ex{\frac{\omega_x}{1 - \omega_x}}
  &= \infty \textrm{ for all choices of } \lambda
\end{align*}
Hence, by \cite[Theorem 1.16]{rwresolomon}:

\begin{theorem}
  \label{thm:add_bias_speed}
  The edge-reinforced random walk on $\bbZ$ with initial bias $\lambda$ has positive speed, if, and only if,
  $\lambda > 3$. If this is the case, then $\frac{X_n}{n} \to \frac{\lambda - 3}{\lambda + 1}$ almost surely.
  Otherwise, $\frac{X_n}{n} \to 0$ almost surely.
\end{theorem}

\subsection{Reinforced Random Walk on Transient Environment}
\label{ssec:trans_env}

We now look at a special case of the LERRW where the initial edge weights
are not all $1$. Define the edge-reinforced random walk on $\bbZ$ with initially $\lambda$-biased environment
as the sequence $X_n$ of random variables with
\begin{align*}
  X_0 &= 0 \\
  \Prb{X_{n+1} = x_n + 1 \midd| X_n = x_n, \ldots, X_0 = x_0} &= \frac{w\left(n, x_n\right)}{w\left(n, x_n\right) + w\left(n, x_n - 1\right)} \\
  w\left(n, z\right) &= \underbrace{w\left(0, z\right)}_{:= \lambda^z} + \sum_{i=0}^{n-1} \mathbbm{1}_{X_i = z, X_{i + 1} = z + 1} + \mathbbm{1}_{X_i = z + 1, X_{i + 1} = z}
\end{align*}
where $\lambda > 0$.
This random walk can again be represented by a mixture of MCs, but the $\omega_x$ are no
longer identically distributed (they are still independent, though).

\begin{figure}[H]
  \begin{center}
    \begin{tikzpicture}
      \node[circle,fill=black,inner sep=0.7mm] (N0) at (0, 0) {};
      \node[circle,fill=black,inner sep=0.7mm] (N1) at (2, 0) {};
      \node[circle,fill=black,inner sep=0.7mm] (N2) at (4, 0) {};
      \node[circle,fill=black,inner sep=0.7mm] (N3) at (6, 0) {};
      \node[circle,fill=black,inner sep=0.7mm] (Nm1) at (-2, 0) {};
      \node[circle,fill=black,inner sep=0.7mm] (Nm2) at (-4, 0) {};
      \node[below=5mm] at (N0) {$0$};
      \node[below=5mm] at (N2) {$z$};
      \draw (Nm2) -- (Nm1) -- (N0) -- (N1) -- node[above=1mm] {$1 + \lambda^{z-1}$} (N2) -- node[above=1mm] {$\lambda^z$} (N3);
      \draw[dashed] (N3) -- (7, 0);
      \draw[dashed] (Nm2) -- (-5, 0);
      \node (aN2l) at (3.9, 0.5) {};
      \node (aN2r) at (4.1, 0.5) {};
    \end{tikzpicture}
    \caption{Edge weights of the walk on transient environment at time of first visit to $z > 0$}
    \label{fig:errw_z_initially_b_env}
  \end{center}
\end{figure}
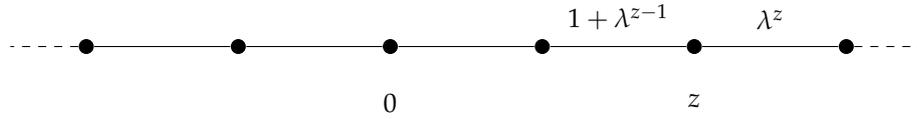

By \autoref{fig:errw_z_initially_b_env}, we get the following distributions for $\omega_x$:
\begin{align*}
  \textrm{for }x > 0\textrm{: }& \qquad \omega_x \sim \Betad{\frac{\lambda^x}{2}}{\frac{1 + \lambda^{x-1}}{2}} \\
  \textrm{for }x = 0\textrm{: }& \qquad \omega_x \sim \Betad{\frac{1}{2}}{\frac{1}{2 \lambda}} \\
  \textrm{for }x < 0\textrm{: }& \qquad \omega_x \sim \Betad{\frac{1 + \lambda^x}{2}}{\frac{\lambda^{x-1}}{2}}
\end{align*}
The electrical network corresponding to the random environment puts the following conductances $\mathbf{c}$
and resistances $\mathbf{r}$ on the edges:
\begin{align*}
  \textrm{for }z &> 0\textrm{: } \qquad c_{\left\{z, z + 1\right\}} = \prod_{x=1}^z \frac{\omega_x}{1-\omega_x}
  \qquad r_{\left\{z, z + 1\right\}} = \prod_{x=1}^z \frac{1-\omega_x}{\omega_x} \\
  \textrm{for }z &= 0\textrm{: } \qquad c_{\left\{0, 1\right\}} = 1
  \qquad r_{\left\{0, 1\right\}} = 1 \\
  \textrm{for }z-1 &< 0\textrm{: } \qquad c_{\left\{z - 1, z\right\}} = \prod_{x=z}^{0} \frac{1-\omega_x}{\omega_x}
  \qquad r_{\left\{z - 1, z\right\}} = \prod_{x=z}^{0} \frac{\omega_x}{1-\omega_x}
\end{align*}
The effective resistance between $0$ and $\infty$ is therefore (using the series law for the negative and
positive half-lines and the parallel law to merge these):
\begin{align*}
  \mathbf{R}_{\textrm{eff}} = \left(\left(1 + \sum_{z \geq 1} \prod_{x=1}^z \frac{1-\omega_x}{\omega_x}\right)^{-1} + \left(\sum_{z \leq 0} \prod_{x=z}^{0} \frac{\omega_x}{1-\omega_x}\right)^{-1}\right)^{-1}
\end{align*}
Hence, the random walk is recurrent if, and only if:
\begin{align*}
  \sum_{z \geq 1} \prod_{x=1}^z \frac{1-\omega_x}{\omega_x} = \infty = \sum_{z \leq 0} \prod_{x=z}^{0} \frac{\omega_x}{1-\omega_x}
\end{align*}
Otherwise, the random walk is transient. To get a better understanding, we calculate expectations.
First, for $x > 0$ using \autoref{lem:beta_distr}:
\begin{align*}
  \Ex{\frac{1 - \omega_x}{\omega_x}}
  &= \begin{cases}
    \frac{1 + \lambda^{x-1}}{2} \cdot \left(\frac{\lambda^x}{2} - 1\right)^{-1}
    = \frac{\lambda^{x-1} + 1}{\lambda^x - 2} & \textrm{ if } \frac{\lambda^x}{2} > 1 \iff \lambda^x > 2 \\
    \infty & \textrm{ otherwise}
  \end{cases}
\end{align*}
And for $x < 0$:
\begin{align*}
  \Ex{\frac{\omega_x}{1 - \omega_x}}
  &= \begin{cases}
    \frac{\lambda^x + 1}{\lambda^{x-1} - 2} & \textrm{ if } \lambda^{x-1} > 2 \\
    \infty & \textrm{ otherwise}
  \end{cases}
\end{align*}
Hence, if $\lambda > 2$:
\begin{align*}
  \Ex{\prod_{x=1}^z \frac{1-\omega_x}{\omega_x}}
  = \prod_{x=1}^z \frac{\lambda^{x-1}}{\lambda^{x-1}} \cdot \frac{1 + \frac{1}{\lambda^{x-1}}}{\lambda - \frac{2}{\lambda^{x-1}}}
  = \prod_{x=1}^z \frac{1 + \frac{1}{\lambda^{x-1}}}{\lambda - \frac{2}{\lambda^{x-1}}}
\end{align*}
Therefore, the expectation decreases exponentially for $\lambda > 2$. This seems to indicate transience for $\lambda > 2$
but a few calculations still remain to formally conclude. Indeed, for any $\lambda > 1$, the expectation
$\Ex{\frac{1 - \omega_x}{\omega_x}}$ will eventually be finite and of order $\frac{1}{\lambda}$ for $x$ large enough,
which can probably be used to show that the walk is transient for all $\lambda > 1$. The final proof is still
open, but will hopefully follow in the near future.

\clearpage
\newpage

\section{Conclusion}
\label{sec:conclusion}

We quickly look at how the considered modifications changed the behavior of the ERRW.

\begin{table}[H]
  \begin{center}
    \def\arraystretch{1.5}
    \begin{tabularx}{\textwidth}{l|X}
      Type of walk & Behavior \\ \hline
      $2$ walkers on $3$-node segment & \parbox{8.5cm}{
        ~\newline
        \textcolor{sectionblue}{Proven:} the proportion of the left edge weight compared to the total edge weight
        forms a martingale at certain stopping times, and converges to a limit (for linear reinforcement), similar to
        the case with $1$ walker.

        \textcolor{linkred}{Conjectured:} the limit of the left edge weight proportion is a random variable which has
        a density w.r.t.~the Lebesgue measure on $\left[0, 1\right]$, it is not Beta-distributed.
        \newline~
      } \\
      $k$ walkers on $\bbZ$ & \parbox{8.5cm}{
        ~\newline
        \textcolor{sectionblue}{Proven:} if all but finitely many initial weights are $1$, then, for almost arbitrary
        reinforcement, either all walkers are recurrent or all walkers have finite range a.s.

        \textcolor{linkred}{Conjectured:} for a certain reinforcement scheme, the $k$ walkers are recurrent if, and only if,
        $1$ walker is recurrent with the same reinforcement scheme.
        \newline~
      } \\
      multiplicative bias on $\bbZ$ & \parbox{8.5cm}{
        ~\newline
        \textcolor{sectionblue}{Proven:} if a node is visited infinitely often, then the probability to go in the direction
        of bias at that node converges to $1$ (for linear reinforcement).

        \textcolor{linkred}{Conjectured:} the walk is transient for any bias $\lambda \neq 1$.
        \newline~
      } \\
      multiplicative bias on triangle & \parbox{8.5cm}{
        ~\newline
        \textcolor{tumOrange}{Unfinished calculation:} the edge weights will eventually all be
        of the same order, and the walker will move around the circle with a constant speed
        (for linear reinforcement).
        \newline~
      } \\
      additive bias on $\bbZ$ & \parbox{8.5cm}{
        ~\newline
        \textcolor{sectionblue}{Proven (up to details):} recurrent for $0 \leq \lambda \leq 1$, transient for $\lambda > 1$.
        Positive speed for $\lambda > 3$ (for linear reinforcement).
        \newline~
      } \\
      initial transient edge weights on $\bbZ$ & \parbox{8.5cm}{
        ~\newline
        \textcolor{tumOrange}{Unfinished calculation:} transient for $\lambda > 1$ (for linear reinforcement).
      }
    \end{tabularx}
  \end{center}
  \caption{Overview of results and conjectures on modified reinforced walks}
  \label{tab:resultoverview}
\end{table}

There are two main conclusions which can be drawn from this overview. First, for reinforced random walks, having multiple random
walkers moving in the same environment which influence each other doesn't seem to fundamentally change the behavior in comparison
to a single walker. Indeed, after working on this topic for such a long period of time, the following (vague) conjecture seems reasonable:
if, for a given graph and given reinforcement scheme, the reinforced random walk is recurrent for a single walker, then it is also
recurrent for any finite number of walkers influencing each other, and vice versa. There is still a lot of work to do to get in any
way closer to prove some form of this conjecture, but analyzing various toy models seems to indicate that this could be true. Of course,
one would also have to formally define a reinforcement scheme to make a true mathematical claim, but the intuitive meaning should be clear.

The second conclusion is that introducing a bias can change the behavior. The bias, which makes the walk more transient, is sometimes
competing with the reinforcement, which makes the walk more recurrent. This result is insofar expected as on regular trees, the reinforcement
is also competing with the transient nature of a random walk on a tree. Indeed, depending on the strength of the reinforcement, there is a
phase transition between recurrence and transience on trees, a result which we find again for the additive bias on $\bbZ$.

All in all, the methods to prove the results given here are not new. Adding multiple walkers and a bias complicates an already complicated
model even more. Therefore, the random walks were only analyzed on very simple graphs, because the methods used here do not work anymore for
more general classes of graphs. Adding a bias or multiple walkers often destroys the property of exchangeability which the basic linearly
reinforced random walk with a single walker possesses (the probability of taking a certain path only depends on how often each edge in the
path is traversed, but not on the order of traversals). This makes the analysis harder and it is no longer possible to represent the reinforced
walk as a RWRE, which was often the tool of choice to prove previous results. One goal of this thesis was to better understand
how reinforced walks react to variations in the model. Some steps in the right direction were made by analyzing the behavior on simple graphs.
At the same time, better tools to analzye these complicated models are still missing. The mathematics here was often adapted from the basic
linear reinforcement case and is thus limited to a very restricted set of graphs.

\subsection{Outlook}

Even though some new results were obtained, many open questions remain. On the one hand, some proofs and calculations in this thesis still
lack some details or formal precision. A natural next step would be to fill these gaps to be sure that the unfinished calculations and proof
sketches actually show what they are supposed to show. However, the answers to more interesting questions are often only conjectured, and no
proof idea has been found yet. The following four main points would be very interesting for future research:
\begin{itemize}
  \item For the $2$ walkers on a $3$-node segment, which can also be seen as a kind of modified, $2$-player P\'olya urn, does the limit of
  the proportion of the left edge weight really have a density w.r.t.~the Lebesgue measure? What is the distribution of the limit?
  \item For $k$ walkers on $\bbZ$, is the behavior for $k$ walkers and for a single walker identical if the same reinforcement scheme
  is used?
  \item For linear reinforcement and additional multiplicative bias on $\bbZ$, is a single walker transient for any bias $\lambda \neq 1$?
  \item For general graphs, does the behavior of the reinforced random walk only depend on the graph and the reinforcement scheme, but not
  on the (finite) number of walkers?
\end{itemize}
In answering these question, another goal would be to find more general techniques to analyze these types of random walks.

This thesis was really only a starting point in better understanding how changes to the model of the reinforced random walk will
affect its behavior. While a definitive answer to the last question listed above still seems a long way off, the other three questions
seem easier to deal with and can hopefully be answered in the near future.

\newpage

\ifthenelse{\boolean{print}}{\thispagestyle{empty} { ~ } \newpage}{}

\clearpage
\newpage
\thispagestyle{empty}
{
	\makeatletter
	\centering
	~\\[19em]
	\def\svgwidth{96.05263069425209pt}\input{graphics/TUM_svg-tex.pdf_tex}\\[10em]
	\def\svgwidth{59.99796366911307pt}\input{graphics/tum_mathematik_svg-tex.pdf_tex}\\

	\makeatother
}

\end{document}

%% file: graphics/tum_mathematik_svg-tex.pdf_tex
\begingroup%
  \makeatletter%
  \providecommand\color[2][]{%
    \errmessage{(Inkscape) Color is used for the text in Inkscape, but the package 'color.sty' is not loaded}%
    \renewcommand\color[2][]{}%
  }%
  \providecommand\transparent[1]{%
    \errmessage{(Inkscape) Transparency is used (non-zero) for the text in Inkscape, but the package 'transparent.sty' is not loaded}%
    \renewcommand\transparent[1]{}%
  }%
  \providecommand\rotatebox[2]{#2}%
  \newcommand*\fsize{\dimexpr\f@size pt\relax}%
  \newcommand*\lineheight[1]{\fontsize{\fsize}{#1\fsize}\selectfont}%
  \ifx\svgwidth\undefined%
    \setlength{\unitlength}{29.76198926bp}%
    \ifx\svgscale\undefined%
      \relax%
    \else%
      \setlength{\unitlength}{\unitlength * \real{\svgscale}}%
    \fi%
  \else%
    \setlength{\unitlength}{\svgwidth}%
  \fi%
  \global\let\svgwidth\undefined%
  \global\let\svgscale\undefined%
  \makeatother%
  \begin{picture}(1,1.00003394)%
    \lineheight{1}%
    \setlength\tabcolsep{0pt}%
    \put(0,0){\includegraphics[width=\unitlength,page=1]{graphics/tum_mathematik_svg-tex.pdf}}%
  \end{picture}%
\endgroup

%% file: graphics/TUM_svg-tex.pdf_tex
\begingroup%
  \makeatletter%
  \providecommand\color[2][]{%
    \errmessage{(Inkscape) Color is used for the text in Inkscape, but the package 'color.sty' is not loaded}%
    \renewcommand\color[2][]{}%
  }%
  \providecommand\transparent[1]{%
    \errmessage{(Inkscape) Transparency is used (non-zero) for the text in Inkscape, but the package 'transparent.sty' is not loaded}%
    \renewcommand\transparent[1]{}%
  }%
  \providecommand\rotatebox[2]{#2}%
  \newcommand*\fsize{\dimexpr\f@size pt\relax}%
  \newcommand*\lineheight[1]{\fontsize{\fsize}{#1\fsize}\selectfont}%
  \ifx\svgwidth\undefined%
    \setlength{\unitlength}{54.75bp}%
    \ifx\svgscale\undefined%
      \relax%
    \else%
      \setlength{\unitlength}{\unitlength * \real{\svgscale}}%
    \fi%
  \else%
    \setlength{\unitlength}{\svgwidth}%
  \fi%
  \global\let\svgwidth\undefined%
  \global\let\svgscale\undefined%
  \makeatother%
  \begin{picture}(1,0.52054795)%
    \lineheight{1}%
    \setlength\tabcolsep{0pt}%
    \put(0,0){\includegraphics[width=\unitlength,page=1]{graphics/TUM_svg-tex.pdf}}%
  \end{picture}%
\endgroup